\documentclass[11pt,letterpaper,reqno]{amsart}%%12
\usepackage{amssymb,a4wide,amsthm,amsmath,amsfonts,xcolor}
\usepackage{mathrsfs}
\def\mathcal{\mathscr}
\usepackage{newsymbol}
\usepackage{verbatim}
\usepackage[active]{srcltx}
\usepackage[colorlinks,linkcolor={blue},
citecolor={blue},urlcolor={red},]{hyperref}
\usepackage{hyperref}
\usepackage{fouridx}
\newcommand{\domain}{\mathcal{O}}

\newcommand{\esssup}{\operatornamewithlimits{ess\,sup}}

\newcommand{\todown}{\searrow}
\newcommand{\rK}{\mathrm{K}}
\newcommand{\D}{\mathrm{D}}

\newcommand{\embed}{\hookrightarrow}

\newcommand\newe{e}
\newcommand\delc[1]{}
\newcommand\comc[1]{{\color{red} #1}}

\definecolor{ForestGreen}{rgb}{0.1,1,0.3}

\newcommand\delj[1]{}

\newcommand\delfeb[1]{}

\newcommand\deld[1]{}% December 30 2015
\newcommand\della[1]{}%% 20 February 2014
\newcommand{\eps}{\varepsilon}

\newcommand\addf[1]{{\color{blue} #1}} %14April 2014
\newcommand\delf[1]{}%14April 2014
\newcommand\delg[1]{}%10 June 2014

\newcommand\addi[1]{{\color{blue} #1}} %03 March 2015
\newcommand\deli[1]{}%03 March 2015

\newcommand{\rH}{\mathrm{ H}}
\newcommand{\rX}{\mathrm{ X}}
\newcommand{\rD}{\mathrm{ D}}
\newcommand{\rA}{\mathrm{ A}}

\newcommand{\lb}{\langle}
\newcommand{\rb}{\rangle}

\let\emptyset \undefined
\let\ge       \undefined
\let\le       \undefined
\newsymbol\le          1336  \let\leq\le
\newsymbol\ge          133E  \let\geq\ge
\newsymbol\emptyset    203F
\newsymbol\notle       230A
\newsymbol\notge       230B
\theoremstyle{plain}
\newtheorem{theorem}{Theorem}[section]
\theoremstyle{remark}
\newtheorem{remark}[theorem]{Remark}

\theoremstyle{plain}
\newtheorem{corollary}[theorem]{Corollary}
\newtheorem{lemma}[theorem]{Lemma}
\newtheorem{proposition}[theorem]{Proposition}
\newtheorem{definition}[theorem]{Definition}

\numberwithin{equation}{section}

\newcommand\del[1]{}
\def\la{\left(}
\def\ra{\right)}
\def\las{\left\langle}
\def\ras{\right\rangle}
\def\E{\mathbb E}
\def\P{\mathbb P}
\def\R{\mathbb R}

\def\HH{\mathbb H}
\newcommand{\imp}{\int\limits}
\newcommand{\LL}{\mathbb L}

\begin{document}
\title[Large deviations for stochastic LLG equation]
{Large deviations and transitions between equilibria for stochastic Landau-Lifshitz-Gilbert equation}
\author[Z. Brze\'zniak]{Zdzis{\l}aw Brze\'zniak}
\address{Department of Mathematics\\
The University of York\\
Heslington, York YO10 5DD, UK} \email{zb500@york.ac.uk}
\author[B. Goldys]{Ben Goldys}
\address{School of Mathematics and Statistics, The University of Sydney, Sydney 2006, Australia}
\email{Beniamin.Goldys@sydney.edu.au}
\author[T. Jegaraj]{Terence Jegaraj}
\address{School of Mathematics and Statistics, UNSW, Sydney 2052, Australia}
\thanks{The  work of Zdzis{\l}aw Brze{\'z}niak and of Ben Goldys was partially supported by the ARC
Discovery Grant DP120101886. The research on which we report in this paper was begun  at the Newton
Institute for Mathematical Sciences in Cambridge (UK) during the
    program "Stochastic Partial Differential Equations". The INI support and excellent working conditions are gratefully acknowledged by all three authors.    The first named author wishes to thank Clare Hall (Cambridge)  and the School of Mathematics, UNSW,  Sydney  for hospitality.
}
\begin{abstract}
We study a stochastic Landau-Lifshitz equation on a bounded interval and with finite dimensional noise. We first show that there exists a pathwise unique  solution to this equation and that this solution enjoys the  maximal regularity property. Next, we prove the large deviations principle for small noise asymptotic of solutions using the weak convergence method. An essential  ingredient   of the proof is  compactness, or weak to strong continuity,   of the solution map for a deterministic Landau-Lifschitz equation, when considered as a transformation of external fields. We then apply this large deviations principle to show that small noise can cause magnetisation reversal. We  also  show the importance of the shape anisotropy parameter for reducing the disturbance of the solution caused by small noise. The problem is motivated by applications of ferromagnetic nanowires to the fabrication of magnetic memories.
\end{abstract}
\date{\today}
\maketitle
\tableofcontents

\newtheorem{sta}{Statement}
\newtheorem{pro}{Proposition}
\newtheorem{thm}[pro]{Theorem}
\newtheorem{lem}[pro]{Lemma}
\newtheorem{cor}[pro]{Corollary}

\newtheorem{ass}[pro]{Assumption}

\newtheorem{defn}[pro]{Definition}

\newtheorem{rem}{Remark}
\newcommand{\lfunctional}[1]{{}_{{#1}^*}\!\langle}
\newcommand{\rfunctional}[1]{\rangle_{#1}}

\newcommand{\xone}{D(A)}
\newcommand{\yn}{y_n}%{y_{w_n}}
\newcommand{\yw}{y}%{y_{w}}

\newcommand\dela[1]{}

\newcommand{\bH}{\mathbb{H}}
\newcommand{\yp}{y^\prime}
\newcommand{\Ep}{E^\prime}
\newcommand{\Fp}{\mathcal{F}^\prime}
\newcommand{\Wp}{W^\prime}
\newcommand{\ydn}{y_n(s)\times\Delta y_n(s)}
\newcommand{\yydn}{y_n(s)\times(y_n(s)\times\Delta y_n(s))}
\newcommand{\yhn}{y_n(s)\times h}
\newcommand{\yyhn}{y_n(s)\times(y_n(s)\times h)}
\newcommand{\psionen}{\psi (|y_n(s)|_{L^\infty})}
\newcommand{\psitwon}{\psi (|\pi_n(y_n(s)\times h)|_{L^\infty})}
\newcommand{\psithreen}{\psi (|\pi_n(y_n(s)\times(y_n(s)\times h))|_{L^\infty})}
\newcommand{\psisqonen}{\psi^2(|y_n(s)|_{L^\infty})}
\newcommand{\psisqtwon}{\psi^2(|\pi_n(y_n(s)\times h)|_{L^\infty})}
\newcommand{\psisqthreen}{\psi^2(|\pi_n(y_n(s)\times(y_n(s)\times h))|_{L^\infty})}

\newcommand{\yhnp}{y^\prime_n(s)\times h}
\newcommand{\yyhnp}{y^\prime_n(s)\times(y^\prime_n(s)\times h)}
\newcommand{\psionenp}{\psi (|y^\prime_n(s)|_{L^\infty})}
\newcommand{\psitwonp}{\psi (|\pi_n(y^\prime_n(s)\times h)|_{L^\infty})}
\newcommand{\psithreenp}{\psi (|\pi_n(y^\prime_n(s)\times(y^\prime_n(s)\times h))|_{L^\infty})}
\newcommand{\psisqonenp}{\psi^2(|y^\prime_n(s)|_{L^\infty})}
\newcommand{\psisqtwonp}{\psi^2(|\pi_n(y^\prime_n(s)\times h)|_{L^\infty})}
\newcommand{\psisqthreenp}{\psi^2(|\pi_n(y^\prime_n(s)\times(y^\prime_n(s)\times h))|_{L^\infty})}
\newcommand{\rT}{\mathrm{T}}

\def\R{\mathbb R}
\def\la{\left(}
\def\ra{\right)}
\def\E{\mathbb E}

\newcommand{\Rthree}[0]{{\mathbb{R}^3}}
\newcommand{\vthree}[3]{\mbox{\scriptsize{$\left(\begin{array}{r}#1\\#2\\#3\end{array}\right)$}\normalsize}}

\newcommand{\dummyarray}{\begin{array}{c}{ }\\{ }\end{array}}
\newcommand{\dummyarraybig}{\begin{array}{c}{ }\\{ }\\{ }\end{array}}

\section{Introduction}\label{sectionIntro}

\noindent
Stochastic PDEs for manifold-valued processes have been first studied  by Funaki \cite{funaki} and  by Carroll and the first named authour in \cite{Carroll_1999,Brz+Carroll},    and recently studied
in \cite{bgo,bo1,bo2}. In this paper we consider a particular example of such an equation known as
the stochastic Landau-Lifshitz-Gilbert (LLG) equation with solutions taking values in the two-dimensional sphere $\mathbb S^2$, see \cite{llg-book} or \cite{Brz+Gold+Jeg_2013}. To introduce this equation, we will need the Sobolev space  $H^{1,2}\la\domain,\mathbb R^3\ra$ of functions defined on a bounded interval $\domain$ of the real line. To every $u\in  H^{1,2}\la\domain,\mathbb R^3\ra$ we associate its energy
\[\mathcal{E}(u)=\frac{a}{2}\int\limits_\domain|\nabla u|^2\,dx + \int\limits_\domain f(u)\,dx,\]
where, for instance,
\[f(y)=\frac{\beta}{2}\la y_2^2+y_3^2\ra,\quad u=\la y_1,y_2,y_3\ra\in\mathbb S^2\,.\]
Let
\[\mathcal H(u)=-\nabla\mathcal{E}(u)=a\Delta u -f^\prime(u)\]
denote the the $L^2$-gradient of the energy functional $\mathcal E$. We will consider
the following  Stratonovitch type stochastic PDE satisfied by a random function $M:[0,\infty)\times \domain\to\R^3$:
 \begin{equation}\label{sllg}
 \left\{\begin{array}{ll}
  dM=\left[ M\times\mathcal H(M)-\alpha M\times(M\times\mathcal H(M))\right]dt+\sqrt{\eps}g(M)\circ d\xi,&t\geq 0,\\
  &\\
   \left.\dfrac{\partial M}{\partial x}\right|_{\partial\domain}=0,\;\; t\geq 0,& \\
   &\\
  M(0)=M_0\,,&
  \end{array}\right.
  \end{equation}
where $\xi$ is a certain $L^2\la\domain,\R^3\ra$-valued Wiener process and the initial data satisfies  $|M_0(x)|=1$ for all $x\in\domain$.
The precise definitions of the noise  and the function $g:\R^3\to\mathcal L\la\R^3\ra$
are provided in Sections \ref{prelim} and \ref{sec_exist}, see \eqref{eqndef-g}.
The definition of solutions to \eqref{sllg} is implicitly given in Theorem \ref{thm-weak solution}.
Formal application of the It\^o formula easily shows that $|M(t,x)|=1$ for all times and all $x\in\domain$ so that \eqref{sllg} is indeed an example of a stochastic PDE for  an $\mathbb S^2$-valued process $M$.
\par\medskip\noindent
Equation \eqref{sllg} with $\eps=0$ is a relatively simple version of the general  Landau-Lifschitz-Gilbert equation that provides a basis for  the theory and applications of ferromagnetic materials, and fabrication of magnetic memories in particular, see for example \cite{llg-book,berkow,gilbert,landau}.
Let us recall that according to the Landau and Lifschitz  theory of ferrormagnetizm \cite{landau}, modified later by Gilbert \cite{gilbert}, the deterministic LLG equation
 \begin{equation}\label{sllgd}
 \left\{\begin{array}{ll}
\dfrac{dm}{dt}=m\times\mathcal H(m)-\alpha m\times(m\times\mathcal H(m)),&t\geq 0,\\
  &\\
   \left.\dfrac{\partial m}{\partial x}\right|_{\partial\domain}=0,& t\geq 0,\\
   &\\
  m(0)=m_0\,,&
  \end{array}\right.
  \end{equation}
describes the evolution of the magnetisation vector $m$ of a ferromagnet occupying the region $\domain$. For the derivation of  equations \eqref{sllg} and \eqref{sllgd} from physical principles and for the physical motivation to add a stochastic term to equation \eqref{sllgd}, see \cite{llg-book,berkow,BGJ2011,KRV,landau}. Here we mention only that the Landau-Lifschitz theory of ferromagnetizm requires coupling of equations \eqref{sllg} and \eqref{sllgd} with the Maxwell equations in the whole space. They need not be introduced in this paper because in one-dimensional domain the effect of coupling is incorporated in the term $\beta g(M)f(M)$, see \cite{carbou1} for details. Finally, we note that the case of one-dimensional domain while being relatively simple (contrary to the multidimensional case, smooth solutions exist) is important for physics of ferromagnetism and applications of ferromagnetic nanowires, see \cite{carbou1}.
\par\medskip\noindent
To the best of our knowledge \eqref{sllg}  has not been studied before. The existence of a weak martingale solution is proved for a similar equation in a three-dimensional domain in our earlier work \cite{Brz+Gold+Jeg_2013}. Kohn, Reznikoff and vanden-Eijnden~\cite{KRV} modelled the magnetisation $M$ in a thin film, assuming that $M$ is constant across the domain for all times and $\beta=0$. In this case \eqref{sllg} reduces to an ordinary stochastic differential equation in $\R^3$.
%We also note that the convergence of finite elements schemes for stochastic LLG equations similar to \eqref{sllg}  have been recently studied in \cite{sllg-book,goldys,hocquet}.
%\[
%dM=M\times(F+\sqrt{\eps}\circ dW)-\alpha M\times(M\times(F+\sqrt{\eps}\circ dW)).\]
%where $W$ is now a standard Wiener process in $\mathbb{R}^3$.
They\delc{Kohn, Reznikoff and Vanden-Eijnden} used the large deviations theory to make a detailed computational and theoretical study of the behaviour of the solution. They also remark that little is known about the behaviour of solutions to the stochastic Landau-Lifshitz equation when $M$ is not constant on the space domain.
\par\medskip\noindent
In this work we address the question raised in \cite{KRV}. We show first the existence and uniqueness of smooth pathwise solutions to \eqref{sllg}. Then we prove the Large Deviations Principle (LDP) for \eqref{sllg} and finally, we apply the LDP to the analysis of transitions between equilibria in the limit of vanishing noise.
 \begin{comment}
Then it is natural to model the dynamics of the magnetisation $M$ using the Stratonovitch stochastic differential equation
\begin{equation}\label{eq-ener}
dM=-M\times\la\nabla_{\mathbb L^2}\mathcal E(M)dt- \circ \sqrt{\eps}dW\ra+\alpha M\times(M\times\la\nabla_{\mathbb L^2}\mathcal E(M)dt-\circ \sqrt{\eps}dW,\ra)
\end{equation}
where $\circ$ denotes the Stratonovitch integral. Let us note that the Stratonovitch integral allows us to preserve the constraint $|M(t,x)|=1$ for all times.
\par\bigskip
In order to formulate the final version of the equation we are going to study, we need some additional notations. Let $\la f_i\ra$ be an orthonormal basis in $\R^3$,
\[f(y)=\la y\cdot f_2\ra f_2+\la y\cdot f_3\ra f_3\quad y\in\R^3,\]

\noindent  Note that since $G(y)h$ is orthogonal to $y\in\R^3$ for every $h\in\R^3$, the formal application of the Ito formula yields $|M(t,x)|=1$.
\par\noindent
\end{comment}
\par\medskip\noindent
 We will describe now the content and new results obtained in this paper.
 \par\medskip\noindent
 We start with Section \ref{prelim} containing some definitions and auxiliary fact that will are needed later.\\
  In Section {\ref{sec_exist} we prove the existence of a weak martingale solution stated in Theorem \ref{thm-weak solution}. The proof combines  the ideas of the proof of the existence theorem in ~\cite{Brz+Gold+Jeg_2013} with the application of the Girsanov theorem. We only sketch the steps that repeat almost verbatim the arguments from \cite{Brz+Gold+Jeg_2013} and concentrate on new arguments. \\
In Section \ref{unique} we consider the existence of strong solutions to the stochastic LLG Equation \eqref{sllg}.
\par\medskip\noindent
In Theorem~\ref{thm-uniqueness pathwise weak} we state a pathwise uniqueness result for solutions
of equation~(\ref{sllg}) with trajectories belonging to the space  $S_T=C([0,T];\LL^2)\cap L^4\la
0,T;H^{1,2}\la\domain,\R^3\ra\ra$.
\par\medskip\noindent
In Section \ref{sec:furtherreg}, we prove maximal regularity of solutions to \eqref{sllg}. Namely, we show that
\[\mathbb{E} \imp_0^T\imp_\domain |\D M(t,x)|^4\, dx\,dt+\imp_0^T\imp_\domain|\Delta M(t,x)|^2\,dx\,dt<\infty.\]
The proof of this result follows from the maximal regularity and ultracontractivity properties of the heat semigroup
generated by the  Laplace operator with the Neumann boundary conditions and the estimates for weak solutions of equation \eqref{sllg} obtained in Theorem
\ref{theoremWE}.
\par\medskip\noindent
The Large Deviations Principle for equation \eqref{sllg} is studied in Section  \ref{sec-LDP}. We first identify the
rate function and prove in Lemma \ref{lem19} that it has compact level sets in the space
\[\mathcal X_T=C\la[0,T]; H^{1,2}\la\domain;\R^3\ra\ra\cap L^2\la 0,T; H^{2,2}\la\domain;\R^3\ra\ra.\]
In particular, we show in Lemma \ref{lem19} certain compactness property of solutions to the deterministic LLG
equation. It seems that such a result is new in the deterministic theory and is of independent interest.\\
The Large Deviations Principle is proved in Theorem \ref{thm-large deviations}.
To prove this theorem, we use the weak convergence method of Budhiraja and
Dupuis~\cite[Theorem 4.4]{Budhiraja+Dupuis_2000}. Following their work we show that the two
conditions of Budhiraja and Dupuis, see Statements~\ref{sta1} and~\ref{sta2} in
Section~\ref{sec-LDP}, are satisfied and then Theorem
\ref{thm-large deviations} easily follows.
We note that our proof is simpler than the corresponding proofs in~\cite{CM}
and~\cite{Duan+Mil_2009} as we do not need to partition the time interval
$[0,T]$ into small subintervals. \\
In Section~\ref{sec:applyldp}
we apply the Large Deviations Principle to a simple stochastic model of
magnetisation in a needle-shaped domain.
we first obtain explicit estimates of the size of  domains of attraction of the
North and South Pole which are stationary solutions for the deterministic LLG
equation. Then we show that in the presence of small noise in equation
\eqref{sllg} there is a positive probability of transitions transitions between
the domains of attraction. Using the Freidlin-Ventzell estimates we obtain
explicit estimates for this probability. These estimate show also the importance
of the parameter $\beta$ (interpreted as the measure of shape anisotropy) for
reducing the disturbance of the magnetisation caused by small noise. The results
we obtain partially answer a question posed in \cite{KRV} and provide a
foundation for the computational study of stability of ferromagnetic nanowires
under the influence of small noise.
\subsection{Notations}
The inner product of vectors $x,y\in\R^3$ will be denoted by $x\cdot y$ and
$|x|$ will denote the Euclidean norm of $x$. We will use the standard notation
$x\times y$ for the vector product in $\R^3$.\\
For  a domain $\domain$ we will use the notation $\LL^p$ for the space
$L^p\la\domain;\R^3\ra$, $\mathbb W^{1,p}$
for the Sobolev space $W^{1,p}\la\domain;\R^3\ra$ and so on. For $p=2$ we will
often write  $\HH^k$ instead of $\mathbb W^{k,2}$.  We will
always emphasize the norm of the corresponding space writing
$|f|_{\mathbb{L}^2}$, $|f|_{\HH^1}$ and so on.
\par
We will also need the spaces $L^p(0,T;E)$ and $C([0,T];E)$ of Bochner
$p$-integrable, respectively continuous,
functions $f:[0,T]\to E$ with values in a Banach space $E$. If $E=\R$ then we
write simply $L^p(0,T)$ and $C([0,T])$.
For a Banach space $E$  we will denote by $\mathcal{L}(E)$ the space of all
linear and bounded maps from $E$ to itself.
\par
Throughout the paper $C$ stands for a positive real constant whose actual value
may vary from line to line. We include an argument list, $C(a_1,\ldots,a_m)$, if
we wish to emphasize that the constant depends only on the values of the
arguments $a_1$ to $a_m$.
\section{Preliminaries}\label{prelim}
We us assume that a $\alpha>0$. Let us denote by $g$   a map $g:\R^3\to\mathcal
L\la\R^3\ra$ defined by
\begin{equation}\label{eqndef-g}
g:\R^3 \ni y \mapsto \bigl\{ \R^3 \ni h \mapsto g(y)h:=y\times h-\alpha y\times (y\times h) \in \R^3\bigr\} \in  \mathcal
L\la\R^3\ra.
\end{equation}
The function $g$ \del{begin a polynomial function (od degree $2$)} is of class
$C^\infty$. In particular, we have
\begin{equation}\label{si2}
\left[g^\prime(y)h\right]z=\D\left[g(y)h\right]z=z\times
h-\alpha[z\times(y\times h)+y\times(z\times h)],\quad h,y,z\in\R^3,
\end{equation}
and \del{$G$ and $G^\prime$ are locally bounded, i.e.} for every $r>0$
\begin{equation}\label{gprime}
\sup_{|y|\le r} \;\Bigl[ \vert g(y)\vert_{\mathbb{R}^3} +\left\vert g^\prime(y)\right\vert_{\mathcal{L}\la \R^3,\mathcal
L(\R^3)\ra} \Bigr]<\infty.
\end{equation}
Clearly, we can define a map  $(u,h)\mapsto (g \circ u)h$, if  $u,h$ belong to some  function spaces of $\R^3$-valued functions on $\domain$.
For instance,  if $u\in\LL^\infty$ and $h\in\LL^2$ then $(g \circ u)h$ is a well
defined element of $\LL^2$.
We will denote by $G$ a Nemytski type map associated with the function $g$. To be
precise, we will
use the notation $G(u)$, if $u\in\LL^\infty$,  for a linear map  defined, for every $q\in [1,\infty]$,  by
\begin{equation}\label{eqndef_function G}
G(u): \LL^q \ni h \mapsto  u\times h-\alpha u\times(u\times h)\in \LL^q\,.
\end{equation}
\par\noindent
For fixed functions $e_i\in\LL^2$, $i=1,2,3$, let  $B:\R^3\to\LL^2$ be a linear
operator defined by
\begin{equation}\label{def_B}
B: \R^3 \ni k \mapsto \sum_{i=1}^3k_ie_i\in \LL^2.
\end{equation}
In the next lemma we use the notation $\newe=\la \newe_i\ra$ and
\[|\newe|_{\,\LL^\infty}=\max_{1 \leq i\le
3}\left|\newe_i\right|_{\,\LL^\infty}\,.\]
\begin{lemma}\label{lem_g-polynomial} Assume that  $q\in [1,\infty]$. Then the  map $G:\mathbb{L}^\infty \to \mathcal{L}(\mathbb{L}^q,\mathbb{L}^q)$ is a  polynomial map, hence of polynomal growth and locally Lipschitz, i.e. there exists $C_0>0$ such that
\begin{equation}
\label{eqn-G-polynomial growth}
\left|G\la u\ra h\right|_{\LL^2}\leq
C_0|h|_{\LL^q}\bigl[\left|u\right|_{\LL^\infty} +\left|u\right|_{\LL^\infty}^2\bigr],\quad h\in\LL^q\,.
\end{equation}
and, for every  $r>0$  there exists $C_r>0$ such that for all
$u_i\in\LL^\infty$, $i=1,2$ satisfying
$\left|u_i\right|_{\LL^\infty}\le r$, one has
\begin{equation}
\label{eqn-G-Lipschitz local}
\left|G\la u_1\ra h-G\la u_2\ra h\right|_{\LL^q}\leq
C_r|h|_{\LL^q}\left|u_1-u_2\right|_{\LL^\infty},\quad h\in\LL^2\,.
\end{equation}
Moreover, there exists $a>0$ such that
\begin{equation}
\label{eqn-G-polynomial growth-H^1-norm}
\left|G\la u\ra h\right|_{\mathbb{H}^1}\leq
a|h|_{\mathbb{H}^1}\bigl[1 +\left|u\right|_{\mathbb{H}^1}^2\bigr],\quad u,h \in \mathbb{H}^1.
\end{equation}
\end{lemma}
\begin{proof} The last part of the above Lemma is a consequenc of the fact that $\mathbb{H}^1$ is an algebra.\end{proof}

Given two vectors $f_2,f_3 \in \mathbb{R}^3$, the  function $f:\R^3\to\R^3$ is
defined by
\begin{equation}\label{eqn_function f}
f(y)=\la y\cdot f_2\ra f_2+\la y\cdot f_3\ra f_3,\quad y\in\R^3\,.
\end{equation}
As before, we will denote by $F$ the Nemytki  map associated with the function
$f$, i.e. for $q\in [1,\infty]$,
\begin{equation}\label{eqndef_function G}
F: \LL^q \ni u \mapsto f\circ u=\la u(\cdot)\cdot f_2\ra f_2+\la
u(\cdot)\cdot f_3\ra f_3 \in \LL^q\,.
\end{equation}
Note that $F: \LL^q \to \LL^q$ is a bounded linear map.  In conjuction with Lemma \ref{lem_g-polynomial}\ we get the following result.
\begin{lemma}\label{lem_G-F-polynomial}
For every $e\in\mathbb{R}^3$, the maps
\begin{eqnarray*}
GF:&& \LL^\infty \ni u \mapsto G(u)F(u) \in \LL^\infty
\\
G^\prime e G e:&& \LL^\infty \ni u \mapsto [G^\prime(u)e][G(u)e] \in \LL^\infty
\end{eqnarray*}
\end{lemma}
\par\noindent

Let us recall that $\domain\subset\R $ is a bounded interval.
We define the Laplacian with the Neumann boundary conditions by
$\rA:D(\rA)\subset \mathbb{L}^2\rightarrow \mathbb{L}^2$ by
\begin{eqnarray}
\left\{\begin{array}{rcl}
D(\rA)&:=&\{u\in\mathbb{H}^2:\, \rD u(x)=0\mbox{ for }x\in\partial\domain\},\\
\rA u&:=&-\Delta u\mbox{ for }u\in D(\rA).\end{array}\right.\label{Adefinition}
\end{eqnarray}
Let us recall that the operator $A$ is self-adjoint and nonnegative and
$D\left(A^{1/2}\right)$ when endowed with the
graph norm coincides with $\mathbb{H}^1$. Moreover, the operator $(A+I)^{-1}$ is
compact.

For any real number $\beta\geq 0$, we write $\mathbb X^\beta$ for the domain of
the fractional power
operator $D\la \rA^\beta\ra$ endowed with the \dela{graph} norm
$\vert x\vert_{\mathbb X^\beta}= \vert (I+\rA)^{\beta}x\vert$  and $\mathbb
X^{-\beta}$ denotes the dual
space of $\mathbb X^\beta$ so
that $\mathbb X^\beta \subset \mathbb{H}=\mathbb{H}^\prime\subset \mathbb
X^{-\beta}$ is a Gelfand triple.
Note that for $\beta\in [0,\frac34)$,
\[
\mathbb X^\beta=\mathbb{H}^{2\beta}.
\]

In what follows we will need the following, well known, interpolation
inequality:
\begin{equation}\label{ineq-interp}
|u|_{\LL^\infty}^2\le k^2|u|_{\mathbb{H}}|u|_{\mathbb{H}^1}\;\;\;\mbox{ }\forall
u\in\mathbb{H}^1,
\end{equation}
where the optimal value of the constant $k$ is
\[k=2\max\la 1,\frac{1}{\sqrt{|\domain|}}\ra.\]

For $v,w,z\in\mathbb{H}^1$ by the expressions $w\times\Delta v$ and
$z\times(w\times\Delta v)$ we understand the unique elements
 of the dual space $(\mathbb{H}^1)^\prime$ of $\mathbb{H}^1$ such that for any $\phi\in\HH^1$
\begin{equation}\label{indualspace1}
\fourIdx{}{(\mathbb{H}^1)^\prime}{}{\mathbb{H}^1}{\langle w\times\Delta v,\,\phi\rangle}
%\lsa{(\mathbb{H}^1)^\prime}w\times\Delta v,\,\phi\rangle_{\mathbb{H}^1}
=-\langle \D(\phi\times w),\,\D v\rangle_{\mathbb{L}^2}
\end{equation}
and
\begin{equation}\label{indualspace2}
%\lsa{(\mathbb{H}^1)^\prime}z\times(w\times\Delta v)\,,\phi\rangle_{\mathbb{H}^1}
\fourIdx{}{(\mathbb{H}^1)^\prime}{}{\mathbb{H}^1}{\langle z\times(w\times\Delta
v)\,,\phi\rangle}=-\langle \D((\phi\times z)\times w),\,\D
v\rangle_{\mathbb{L}^2},
\end{equation}
respectively. Note that the space $H^1\la \domain\ra$ is an algebra, hence for
$v,w,z\in\mathbb{H}^1$, linear functionals $\mathbb{H}^1\ni \phi \mapsto \mbox{ RHS of }$
\eqref{indualspace1} (or \eqref{indualspace2}) are continuous. In particular,
since $\langle a\times b,a\rangle=0$ for $a,b\in\mathbb{R}^3$, we obtain
\begin{eqnarray}
\label{indualspace1'}
%\lsa{(\mathbb{H}^1)^\prime}w\times\Delta v,\,v\rangle_{\mathbb{H}^1}
\fourIdx{}{(\mathbb{H}^1)^\prime}{}{\mathbb{H}^1}{\langle w\times\Delta v,\,v\rangle}
&=&-\langle  v\times \D w,\,\D v\rangle_{\mathbb{L}^2}
\\
\label{indualspace2'}
%\lsa{(\mathbb{H}^1)^\prime}z\times(v\times\Delta v)\,,\phi\rangle_{\mathbb{H}^1}
\fourIdx{}{(\mathbb{H}^1)^\prime}{}{\mathbb{H}^1}{\langle z\times(v\times\Delta
v)\,,\phi\rangle}
&=&-\langle \D(\phi\times z)\times v,\,\D v\rangle_{\mathbb{L}^2},
\end{eqnarray}
and since $a\times a=0$ for $a\in\mathbb{R}^3$, equation \eqref{indualspace1}
yields
\begin{equation}\label{indualspace1''}
%\lsa{(\mathbb{H}^1)^\prime}\phi\times\Delta v,\,\phi\rangle_{\mathbb{H}^1}
\fourIdx{}{(\mathbb{H}^1)^\prime}{}{\mathbb{H}^1}{\langle \phi\times\Delta v,\,\phi\rangle}
=-\langle \D(\phi\times \phi),\,\D v\rangle_{\mathbb{L}^2}=0.
\end{equation}
The maps $\mathbb{H}^1\ni y\mapsto y\times\Delta y\in(\mathbb{H}^1)^\prime$ and $\mathbb{H}^1\ni
y\mapsto y\times(y\times\Delta y)\in(\mathbb{H}^1)^\prime$ are continuous homogenous
polynomials of degree $2$, resp. $3$ hence they are locally Lipschitz
continuous.

\bigskip
\par\noindent

\par\noindent
 \section{The Existence of solutions}\label{sec_exist}
We will be concerned with the following stochastic integral equation form of
problem \eqref{sllg}
\begin{equation}\label{eqn-main}
\begin{aligned}
M(t)=M_0&+\imp_0^t[M(s)\times\Delta M(s))-\alpha M(s)\times(M(s)\times\Delta M(s))]\,ds \\
&+\sqrt{\eps}\imp_0^tG(M(s))B\,dW(s)+\frac{\eps}{2}\sum_{i=1}^3\imp_0^t\,
\big[G^\prime(M(s))\newe_i\big]\big(G(M(s))\newe_i)\,ds\\
&-\beta\imp_0^tG(M(s))F(M(s))\,ds,\quad\quad t\in [0,T].
\end{aligned}
\end{equation}
where $G$ and $F$ are the Nemytski maps associated with functions $g$ and $f$
defined in the previous section. For instance, since $\mathbb{H}^1$ is an algebra, the map
\[
G(\cdot)B:\mathbb{H}^1\ni M\mapsto \bigl\{ \mathbb{R}^3 \ni k \mapsto \sum_{i=1}^3 G(M)(k_ie_i)=\sum_{i=1}^3 k_iG(M)e_i \in \mathbb{H}^1 \bigr\} \in \mathcal{L}(\mathbb{R}^3,\mathbb{H}^1)
\]
is a continuous polynomial function (and hence of $C^\infty$-class and Lipschitz on balls).
\begin{comment}
Introducing the notations
\[\widehat F(s,W)=\sum_{i=1}^3F_i(s,W)\newe_i,\quad\mathrm{and}\quad
\widehat{W}(t)=\sum_{i=1}^3W_i(t)\newe_i\,,\]
we will sometimes use shorter notations
\[\imp_0^tG(M(s))\,d\widehat{W}(s)=\sum_{i=1}^3\imp_0^tG(M(s))\newe_i\,dW_i(s)\]
and
\[\imp_0^tG(M)\widehat
F(s,W)\,ds=\sum_{i=1}^3\imp_0^tF_i(s,W)G(M)\newe_i\,ds\,.\]
\end{comment}

Note, that the expresssion
\[\sqrt{\eps}\imp_0^tG(M(s))B\,dW(s)+\frac{\eps}{2}\sum_{i=1}^3\imp_0^t \,
\big[G^\prime(M(s))\newe_i\big]\big(G(M(s))\newe_i)\,ds\]
can be identified with the Stratonovich integral
\[\sqrt{\eps}\imp_0^tG(M(s))B\,\circ dW(s)\]
but we will not use this concept in the paper.
\par\noindent
We will now formulate the main result of this Section.
\del{Benie, I see a problem here. One of the given data is a process
$F\in\mathcal P_T$. It is defined with respect to the original filtration. When
in the proof we get a new probability space and a new filtration,
then the space $\mathcal P_T$ should also be modified. But if so, where to $F$
should then belong to?}
\begin{theorem}[Existence of a weak martingale solution]\label{thm-weak solution}
Assume that  $\newe=(\newe_i)_{i=1}^3\in\big(\HH^1\big)^3$, $\Vert e \Vert_{\mathbb{H}^1} \leq r$  and that function $f$
 defined
by \eqref{eqn_function f} is fixed. Assume also that  $M_0\in\HH^1$, $\Vert M_0 \Vert_{\mathbb{H}^1} \leq \rho$. Then  there
exists a system
\begin{equation}\label{eqndef-weak mart sol}
 \Bigl( \Omega,\mathcal{F},\mathbb{F},\mathbb{P},W, M\Bigr)
\end{equation}
consisting of a  probability space $\la\Omega,\mathcal F,\mathbb{P}\ra$,
of a filtration $\mathbb{F}=\la\mathcal{F}_t\ra$,
 of a canonical $\R^3$-valued   $\mathbb{F}$-Wiener process $W=\la{W(t)}\ra$
 and of an $\mathbb{F}$-progressively measurable process $M=\la{M(t)}\ra$ such
that
\begin{enumerate}
\item for each $\beta<\frac{1}{2}$ the paths of $M$ are continuous  $\mathbb
H^{2\beta}$-valued functions $\P$-a.s.;
\item For every $p\ge 1$ and every $T>0$,
\begin{equation}\label{ineqn-H^1 norm}\mathbb{E}\sup_{t\in
[0,T]}|M(t)|^{p}_{\mathbb{H}^1}\le C\la T,p,\alpha,\rho,r\ra;\end{equation}
\item For almost every $t\in [0,\infty)$, $M(t)\times\Delta M(t)\in\LL^2$  and
every $T>0$ we have
\begin{equation}\label{ineqn-projected Laplacian}
\mathbb{E}
\left(\imp_0^T|M(s)\times\Delta M(s)|_{\mathbb{L}^2}^2\,ds\right)^p\le
C\la T,p,\alpha,\rho,r\ra
\end{equation}
\item\label{theoremWE-04} $|M(t)(x)|_{\mathbb{R}^3}=1$ for all $x\in\domain$ and
for all $t\in [0,\infty)$, $\mathbb{P}$-a.s.;
\item For every $t\in [0,\infty)$ equation \eqref{eqn-main} holds $\mathbb P$-a.s.
\item  for every $\alpha\in\left(0,\frac{1}{2}\right)$, $ \mathbb P$-a.s.,
\begin{equation}\label{eqn-u_regularity_2}
 u(\cdot)\in C^{\alpha}\left([0,T],\mathbb L^{2}\right).\end{equation}

\end{enumerate}
\end{theorem}
\noindent Note that in Theorem~\ref{thm-weak solution}, $M$ is an $\mathbb{H}^1$-valued
process, hence the expressions $M(s)\times\Delta M(s)$ and
$M(s)\times(M(s)\times\Delta M(s))$  are interpreted in the sense of
\eqref{indualspace1} and \eqref{indualspace2} respectively.
\begin{proof}
The proof of Theorem~\ref{thm-weak solution} is very similar to the proof of Theorem 2.7
in \cite{Brz+Gold+Jeg_2013}.
Here we only sketch the main arguments. Full details can be found in
\cite{Brzezniak+Liang Li_2015}.
It is sufficient to prove the theorem for a bounded time interval $[0,T]$.
We start with some auxiliary definitions.
For each $n\in \mathbb{N}$, let $\mathbb{H}_n$ be the linear span of the first $n$
elements of the orthonormal basis of $\mathbb{L}^2$ composed of eigenvectors of $\rA$ and let
\begin{equation}\label{eqn-pi_n}
\pi_n:\mathbb{L}^2 \to \mathbb{H}_n
\end{equation}
be the corresponding orthogonal projection.
Let us define a map $G_n :\mathbb{H}_n \to \mathcal{L}(\mathbb{H}_n)$ by
\[G_n(u)=\pi_nG\la\pi_nu  \ra\pi_n ,\quad u\in \mathbb{H}_n.\]
and let $G_n^\prime : \mathbb{H}_n \to \mathcal{L}\bigl(\mathbb{H}_n,\mathcal{L}(\mathbb{H}_n)\bigr)$
be the
Fr{\'e}chet derivative of $G_n$. Since the space $\mathbb{H}_n $ is finite dimensional and contained in $\LL^\infty$, \\
\begin{comment}
 Let $\psi\in C_c^\infty(\R)$ be defined as follows: $\psi(r)=1$ for all
$r\in[0,|h|_{L^\infty}+1]$ and the support of $\psi$ is contained in
$(-1,|h|_{L^\infty}+2)$.
 \end{comment}
 For each $n\in \mathbb{N}$, we define a process
$M_n:[0,T]\times\Omega\rightarrow\mathbb{H}_n$ to be a solution of the following ordinary stochastic differential
equation on  $\mathbb{H}_n$:
\begin{eqnarray}\label{eqn-Galerkin}%\!\!\!\!\!\!\!\!\!\!\!\!\!\!
M_n(t)
&=&\pi_n M_0+\imp_0^t\pi_n(M_n\times\Delta M_n)\,ds\\
&-&\alpha\imp_0^t\pi_n(M_n\times(M_n\times\Delta M_n))\,ds \nonumber\\
&+&\sqrt{\eps}\imp_0^tG_n\la M_n\ra
BdW(s)+\frac{\eps}{2}\sum_{i=1}^3\imp_0^t\left[G_n^\prime\la M_n\ra
\newe_i\right]\la G_n\la M_n\ra \newe_i \ra\,ds\nonumber\\
&-&\beta\imp_0^tG_n\la M_n\ra F\la M_n\ra\, ds\,.
\nonumber
\end{eqnarray}
\begin{comment}
where for $f\in\LL^\infty$
\[\Psi_n(f)=\psi\la |f|_{\LL^\infty}\ra\psi\la\left|\pi_n(f\times
h)\right|_{\LL^\infty}\ra
\psi\la\left|\pi_n(f\times(f\times h))\right|_{\LL^\infty}\ra.\]
We note that in equation \eqref{eqn-Galerkin} the term containing $\frac{\eps}{2}$
is not a Stratonovitch correction term but we show that solutions converge to
the solution of equation \eqref{eqn-main}.
\end{comment}
Since the space $\mathbb{H}_n $ is finite dimensional and contained in $\LL^\infty$, by Lemmata \ref{lem_g-polynomial} and \ref{lem_G-F-polynomial}, the maps
the  maps $G_n$, $\left[G_n^\prime\la \cdot\ra
\newe_i\right]\la G_n\la \cdot\ra \newe_i \ra$ and  $G_n\la \cdot\ra F\la \cdot\ra$ are bounded polynomial maps on $\mathbb{H}_n$, hence locally Lipschitz and of polynomal growth. Since the coefficients in \eqref{eqn-Galerkin} are of one-sided linear growth,
by standard arguments we can prove, see e.g. \cite{Albeverio+Brzezniak+Wu_2010},  that for each
$n\in \mathbb{N}$,
equation \eqref{eqn-Galerkin} has a unique strong (in the probabilistic sense) solution.
\noindent Applying the It\^{o} formula and the Gronwall Lemma to the processes
$|M_n(\cdot)|^2_{\mathbb{H}}$ and $|M_n(\cdot)|^2_{\mathbb{H}^1}$, one can
obtain the following,
uniform in $n\in \mathbb{N}$, estimates.
\begin{lemma}\label{UB1}
Let the assumptions of Theorem \ref{thm-weak solution} be satisfyied. Then for each
$n\in \mathbb{N}$
\[|M_n(t)|_{\mathbb{L}^2}=|\pi_nu_0|_{\mathbb{L}^2},\quad\mathrm{for\,\,
all}\quad t\in [0,T]\quad \mathbb{P}-a.s.\]
Moreover, for each $p\in[1,\infty)$ there exists a constant $C\la T,p,\alpha,
\rho,r\ra$ such that, if $\left\|M_0\right\|_{\mathbb{H}^1}\leq \rho$ and $\|e\|_{\mathbb{H}^1}\leq r$, then  for every $n\in \mathbb{N}$
\begin{equation}\label{h1n}
\mathbb{E}\sup_{t\in [0,T]}|M_n(t)|^{p}_{\mathbb{H}^1}\le C\la T,p,\alpha,
\rho,r\ra,
\end{equation}
\[\mathbb{E}\la\imp_0^T|M_n(s)\times\Delta M_n(s)|^2_{\mathbb{L}^2}\,ds\ra^p\le
C\la T,p,\alpha, \rho,r\ra\]
and
\[\mathbb{E}\left(\imp_0^T|M_n(s)\times(M_n(s)\times\Delta
M_n(s))|^2_{\mathbb{L}^2}\,ds\right)^{p}\le C\la T,p,\alpha,
\rho,r\ra.\]
\end{lemma}

The above {\em a priori} estimates from Lemma \ref{UB1} on the sequence $(M_n)$
imply, by applying
two key results of Flandoli and G\c{a}tarek~\cite[Lemma 2.1 and Theorem
2.2]{FG},
that the corresponding sequence of laws of pairs
$(W,M_n)$ is tight on the space $C\la [0,T],\R^3\ra\times \bigl[ C([0,T];\mathbb
X^{-1/2})\cap
L^4(0,T;\mathbb{L}^4)\bigr]$ and hence by the Prokhorod Theorem, modulo
extracting a subsequnce,
these laws converge weakly to a Borel probability measure $\P^{W,M}$ on $C\la
[0,T],\R^3\ra\times \bigl[ C([0,T];\mathbb X^{-1/2})\cap
L^4(0,T;\mathbb{L}^4)\bigr]$. Next  we have the following result.
\begin{proposition}\label{skorohodconstruction}
There exists a probability space
$(\Omega^\prime,\mathcal{F}^\prime,\mathbb{P}^\prime)$
and there exists a sequence $(W^\prime_n,M^\prime_n)$ of $C\la
[0,T],\R^3\ra\times \bigl[ C([0,T];\mathbb X^{-1/2})\cap
L^4(0,T;\mathbb{L}^4)\bigr]$-valued random variables defined on
$(\Omega^\prime,\mathcal{F}^\prime,\mathbb{P}^\prime)$
such that the laws of $(W,M_n)$ and $(W^\prime_n,M_n^\prime)$ are
equal for each $n\in \mathbb{N}$ and $(W_n^\prime,M^\prime_n)$ converges
pointwise
in $C\la [0,T],\R^3\ra\times C([0,T];\mathbb X^{-\frac{1}{2}})$,
$\mathbb{P}^\prime$-a.s.,
to a limit $(W^\prime,M^\prime)$ whose law is equal to  $\P^{W,M}$.
\end{proposition}
\begin{proof}
 The proposition
follows from the Skorohod theorem (see~\cite[Theorem 4.30]{Kall}).
\end{proof}
It remains to show that  the pointwise limit $(W^\prime,M^\prime)$ defined on
the probability space $(\Omega^\prime,\mathcal{F}^\prime,\mathbb{P}^\prime)$
satisfies all the claims of Theorem~\ref{thm-weak solution}. For each $n\in \mathbb{N}$,
$(W^\prime_n,M^\prime_n)$ satisfies an equation obtained from (\ref{eqn-Galerkin}) by
replacing $W$ and $M_n$ by $W^\prime_n$ and $M^\prime_n$, respectively. Then the
processes $M^\prime_n$ satisfy the estimates of Lemma~\ref{UB1}. These estimates
together with the pointwise convergence of the sequence
$((W^\prime_n,M^\prime_n))_{n\in \mathbb{N}}$ imply that the
$(W^\prime,M^\prime)$ satisfies equation \eqref{eqn-main}. The proof of part (5) of
Theorem~\ref{thm-weak solution} is analogous to the proofs of Lemma 5.1 and Lemma 5.2 in
\cite{Brz+Gold+Jeg_2013}. \\
The proof of part (6) is similar to the proof of inequality (2.17) in Theorem 2.7(c) in \cite{Brz+Gold+Jeg_2013}. The only difference being the last two terms on the RHS of equation \eqref{eqn-main}. However, by part (5), the integrands in these terms are uniformly bounded and hence by Lemmata \eqref{lem_g-polynomial} and \ref{lem_G-F-polynomial} we infer that the expectation of the increments corresponding to these terms is Lipschitz with respect to the time parameters.
\end{proof}
\section{The pathwise uniqueness and the existence of a strong solution}\label{unique}

The main result in this section is Theorem~\ref{thm-uniqueness pathwise weak}, on pathwise
uniqueness of solutions of equation
(\ref{eqn-main}).  Although we could have formulated a theorem of  Yamada-Watanabe type
on the uniqueness in law and the existence of a strong solution to equation
(\ref{eqn-main}) we have decided to do so at the end of the next section after we had proved some further regularity properties of the solutions.

We start with a simple

\begin{lemma}\label{lem_harmonicmapsform}
Let $u$ be an element of $\mathbb{H}^1$ such that
\begin{equation}\label{eqn-sphere constraint}
|u(x)|=1\;\mbox{ for all }x\in \domain.
\end{equation}
Then, in $(\mathbb{H}^1)^\prime$, we have
\begin{equation}\label{harmonicmap}
 u\times(u\times\Delta u)= -|\D u|^2u-\Delta u.
\end{equation}
\end{lemma}
\begin{proof} Let us choose and fix $u, \phi\in \mathbb{H}^1$. Note that by \eqref{eqn-sphere constraint}, $\D \vert u\vert^2=0$.
By equality \eqref{indualspace2'} and the product rule we have
\[
-\fourIdx{}{(\mathbb{H}^1)^\prime}{}{\mathbb{H}^1}{\langle u\times(u\times\Delta
u)\,,\phi\rangle}
=\langle \D(\phi\times u)\times u,\,\D u\rangle_{\mathbb{L}^2}
=\langle (\D\phi\times u)\times u,\,\D u\rangle_{\mathbb{L}^2}
+\langle (\phi\times \D u)\times u,\,\D u\rangle_{\mathbb{L}^2}
.\]

Invoking a well known identity
\[a\times(b\times c)=(a\cdot c) b-(a\cdot b) c,\quad a,b,c\in\R^3,\]
we obtain
\begin{eqnarray*}
  \langle (\D\phi\times u)\times u,\,\D u\rangle_{\mathbb{L}^2} &=& \langle (\D\phi \cdot u) u,\,\D u\rangle_{\mathbb{L}^2}
  -\langle (u \cdot  u)\D\phi,\,\D u\rangle_{\mathbb{L}^2}
  \\
=\frac12 \int_{\domain}  (\D\phi(x) \cdot u(x)) \D \vert u(x) \vert^2\, dx
  -\langle \vert u \vert^2\D\phi,\,\D u\rangle_{\mathbb{L}^2}
  &=&
  -\langle \D\phi,\,\D u\rangle_{\mathbb{L}^2}
\end{eqnarray*}
and similarly
\begin{eqnarray*}
 \langle (\phi\times \D u)\times u,\,\D u\rangle_{\mathbb{L}^2}
  &=& \langle (\phi \cdot  u)\D  u,\,\D u\rangle_{\mathbb{L}^2}-\langle ( \D u \cdot u)\phi,\,\D u\rangle_{\mathbb{L}^2}
  \\
= \langle (\phi \cdot  u)\D  u,\,\D u\rangle_{\mathbb{L}^2}-\frac12 \langle ( \D \vert u\vert^2 )\phi,\,\D u\rangle_{\mathbb{L}^2}
&=& \langle (\phi \cdot  u)\D  u,\,\D u\rangle_{\mathbb{L}^2}.
\end{eqnarray*}

Therefore we obtain
\begin{eqnarray*}
-{}_{(\mathbb{H}^1)^\prime}\langle u\times(u\times\Delta
u),\phi\rangle_{\mathbb{H}^1}&=& -\langle \D\phi,\,\D u\rangle_{\mathbb{L}^2} + \langle (\phi \cdot  u)\D  u,\,\D u\rangle_{\mathbb{L}^2}\\
&=&{}_{(\mathbb{H}^1)^\prime}\!\las |\D u|^2u+\Delta u,\phi\ras_{\mathbb{H}^1}
\end{eqnarray*}
\end{proof}

The following uniqueness result applies to a more general problem than \eqref{eqn-main}. It will be used, in this generality, in the uniqueness part of the proof of Theorem \ref{thm-with h}.

\begin{theorem}[Pathwise uniqueness]\label{thm-uniqueness pathwise weak}
Assume that $(\Omega,\mathcal{F},\mathbb{F}, \mathbb{P})$, where
$\mathbb{F}=(\mathcal{F}_t)_{t\in[0,T]}$,
is a filtered probability space
 and    $W=(W(t))_{t\in[0,T]}$ is an $\mathbb{R}^3$-valued $\mathbb{F}$-Wiener
process.
 Assume that  $\newe=(\newe_i)_{i=1}^3\in\big(\HH^1\big)^3$.
\noindent Let $M_1,\;M_2:[0,T]\times\Omega\rightarrow\mathbb{H}$ be
$\mathbb{F}$-progressively measurable continuous
processes such that, for $i=1,2$, the paths of $M_i$ lie in $L^4\la
0,T;\mathbb{H}^1\ra$, satisfy
property \eqref{theoremWE-04} from Theorem \ref{thm-weak solution}
 and each $M_i$ satisfies the equation
\begin{eqnarray}
M_i(t)&=&M_0+\imp_0^tM_i\times\Delta
M_i\,ds-\alpha\imp_0^tM_i\times(M_i\times\Delta M_i)\,ds\nonumber\\
&&{}+\sqrt{\eps}\imp_0^tG\la M_i\ra
BdW(s)+\frac{\eps}{2}\sum_{j=1}^3\imp_0^t\left[G^\prime\la M_i\ra \newe_j\right]
G\la M_i\ra \newe_j\,ds\nonumber\\
&&{}-\beta\imp_0^tG\la M_i\ra  F\la M_i\ra\, ds+ \imp_0^tG\la M_i\ra Bh(s)\, ds
\label{eqbox12}
\end{eqnarray}
for all $t\in [0,T]$, $\mathbb{P}$-almost everywhere. Then
\[M_1(\cdot,\omega)=M_2(\cdot,\omega),\quad\mathrm{for}\quad \mathbb{P}-a.e.\,\,
\omega\in\Omega.\]
\end{theorem}

\begin{proof}
Let us fix $h\in\mathcal P_T$ and let $R>0$ be such that
\[
\imp_0^T \vert h(t)\vert^2 \, dt \leq R^2,\quad\mathbb{P}-\mbox{a.s.}
\]
Note that the above implies that
\begin{equation}\label{ineqn-F-bis}
\imp_0^T \vert h(t)\vert \, dt \leq R\sqrt{T},\quad\mathbb{P}-\mbox{a.s.}
\end{equation}
First, we  note that by Lemma~\ref{lem_harmonicmapsform} the following equality
holds in $\mathbb X^{-1/2}$.
\[M_i(s)\times(M_i(s)\times\Delta M_i(s))=-|\D M_i(s)|^2M_i(s)-\Delta
M_i(s).\]
\begin{comment}
Then using the It\^{o} formula and the Gronwall lemma we obtain
\[\mathbb{E}\sup_{t\in [0,T]}\left|M_1\la t\ra-M_2\la
t\ra\right|^4_{\mathbb{H}}=0.\]
\addf{\[\mathbb{E}\left|M_1\la t\ra-M_2\la t\ra\right|^2_{\mathbb{H}}=0, \mbox{
for } t \in [0,T].\]
}
\end{comment}
Let us assume that $M_1$ and $M_2$ are two solutions satisfying all
assumptions.
Because both $M_i$ satisfy \eqref{theoremWE-04} from Theorem \ref{thm-weak solution},
we infer that
$\vert M_i\vert $ are  bounded. Hence,  by the local  Lipschitz property of
maps
$G$, $G^\prime$ and $f$, as well by the assumptions that each $\newe_i\in\mathbb
L^\infty$,
there exists a constant $C_1>0$, such that for all $t \in [0,T]$,
\begin{eqnarray}\label{ineq-aux-001}
&&\hspace{-9truecm}\lefteqn{
\sum_{i=1}^3\vert G(M_2(t))\newe_i-G(M_1(t))\newe_i\vert_{\LL^2}^2 \leq C_1\vert
\newe\vert_{\LL^\infty}^2 \vert M_2(t)-M_1(t) \vert_{\LL^2}^2,}
\\
&&\hspace{-12.5truecm}\lefteqn{
\sum_{i=1}^3\vert
G^\prime(M_2(t))\newe_iG(M_2(t))\newe_i-G^\prime(M_1(t))\newe_iG(M_1(t))\newe_i
\vert_{\LL^2}^2 \leq C_1 \vert \newe\vert_{\LL^\infty}^2
\vert M_2(t)-M_1(t) \vert_{L^2}^2,}
\label{ineq-aux-002}
\\
\label{ineq-aux-003}&&\hspace{-12.5truecm}\lefteqn{\lb \big[ \big(G\la M_2(t)\ra - G\la M_1(t)\ra \big) \big] B h(s),
M_2(t)-M_1(t)\rb_{\LL^2}
\leq C  \vert h(t) \vert \vert
\newe\vert_{\LL^\infty} \vert  M_2(t) - M_1(t)  \vert_{\LL^2}^2
}
\\
&&\hspace{-12.5truecm}\lefteqn{
\lb G\la M_2(t)\ra   F\la M_2(t)\ra- G\la M_1(t)\ra   F\la M_1(t)\ra, M_2(t) -
M_1(t) \rb_{\LL^2}
\leq C_1 \vert M_2(t)-M_1(t) \vert_{\LL^2}^2}.
\label{ineq-aux-004}
\end{eqnarray}

Let $Z=M_2-M_1$. Then  the process $Z$ belongs to
$\mathbb{M}^2(0,T;V)\cap L^2(\Omega, C([0,T];\mathbb{H})$
and by Lemma~\ref{lem_harmonicmapsform} is a weak solution of the problem
\begin{eqnarray}
\label{eqbox14}
dZ(t)&=& \alpha A Z\, dt+ \Bigl[ \alpha \big( \vert \D M_2\vert^2 M_2 -
\vert \D M_1\vert^2 M_1\big)\Bigr]\,dt \\
&+&\Bigl[ M_2\times\Delta M_2 -M_1\times\Delta M_1 \Bigl[\,dt \nonumber\\
&+&\sqrt{\eps}\bigl( G\la M_2 \ra- G\la M_1 \ra \bigr) \,Bd W(s)\nonumber\\
&+&\frac{\eps}{2}\sum_{j=1}^3 \Big[ G^\prime\la M_2\ra \newe_j G\la M_2\ra
\newe_j- G^\prime\la M_1\ra \newe_j G\la M_1\ra \newe_j\Big] \,dt\nonumber\\
&-&\beta  \Big[ G\la M_2\ra   F\la M_2\ra- G\la M_1\ra   F\la M_1\ra\Big]\, dt.
\nonumber\\
&+& \Big[ \big(G\la M_2(t)\ra - G\la M_1(t)\ra \big) \Big] B h(s)\,dt
\nonumber
\end{eqnarray}
We can check that all assumptions of the It\^o Lemma from \cite{Pardoux_1979}
are satisfied and therefore
\begin{eqnarray}
\nonumber
\frac12 d\vert Z(t)\vert_{\rH}^2&=&- \lb A Z, Z\rb \, dt
\\
&+& \alpha  \lb \vert \D M_2(t)\vert^2 M_2(t)\, dt - \vert \D
M_1(t)\vert^2 M_1(t), Z\rb \, dt
\nonumber
\\ &+& \alpha  \lb \big( \D M_1(t)+\D M_2(t))M_1(t) \D Z, Z\rb\, dt
\nonumber \\
&+& \Big[ \lb  M_2(t)\times\Delta Z, Z\rb  -  \lb Z \times\Delta M_1(t), Z
\rb\Big] \,dt\nonumber\\
&+&\frac{\eps}{2}\sum_{j=1}^3 \lb G^\prime\la M_2(t)\ra \newe_j G\la M_2(t)\ra
\newe_j- G^\prime\la M_1(t)\ra \newe_j G\la M_1(t)\ra \newe_j, Z\rb
\,dt\nonumber\\
&&{}-\beta\lb G\la M_2(t)\ra   F\la M_2(t)\ra- G\la M_1(t)\ra   F\la M_1(t)\ra,
Z \rb \, dt
\nonumber\\
&+& \lb \big[ \big(G\la M_2(t)\ra - G\la M_1(t)\ra \big) \big] B h(s), Z\rb
\,dt
\nonumber
\\
&+& \frac12 \eps\sum_{j=1}^3 \vert \big( G\la M_2(t) \ra- G\la M_1(t) \ra \big)
\newe_j \vert_{\rH}^2 dt
\nonumber
\\
&&{}+\sqrt{\eps}\sum_{j=1}^3 \lb  G\la M_2(t) \ra- G\la M_1(t) \ra \big)
\newe_j, Z\rb \,dW_j(s)\nonumber\\
\nonumber
\\
&=&\sum_{i=1}^{8}I_i(t) \, dt +\sum_{j=1}^3 I_{9,j}(t)\,dW_j(t)
\label{eqbox15}
\end{eqnarray}
We will estimate all the terms in \eqref{eqbox15}. In what follows we will often
use inequality \eqref{ineq-interp} and $k$ is the constant from that inequality.
 \del{Let us recall that $\vert M_i(t)\vert_{\LL^\infty} = 1$ for $i=1,2$.} Let
us start with the 1$^{\rm st}$ term:
\begin{eqnarray*}
I_1(t)=-\lb A Z(t), Z(t)\rb &=&- \vert \D Z(t)\vert^2.
\end{eqnarray*}
As for the 2nd term we have
\begin{eqnarray*}
&&\hspace{-2truecm}\lefteqn{\lb \vert \D M_2\vert^2 M_2- \vert \D
M_1\vert^2 M_1, Z\rb }
\\
&=&  \lb \vert \D M_2\vert^2 Z , Z\rb  + \lb ( \D M_1+ \D M_2) M_1
\D Z,  Z\rb =:II_0+\sum_{i=1}^2 II_i\,.
\end{eqnarray*}
Next,
\begin{eqnarray*}
II_0 &\leq &   \vert \D M_2\vert_{\LL^2}^2 \vert Z\vert_{\LL^\infty}^2
\\
 &\leq &  k^2 \vert \D M_2\vert_{L^2}^2 \vert Z\vert_{\LL^2} \vert
Z\vert_{\mathbb{H}^1}
 \\
  &\leq &  k^2 \vert \D M_2\vert_{\LL^2}^2 \vert Z\vert_{\LL^2} \big( \vert
Z\vert_{\LL^2} +\vert \D Z\vert_{\LL^2}\big)
  \\ &\leq & k^2 \vert \D M_2\vert_{\LL^2}^2 \vert Z\vert_{\LL^2}^2 + k^2
\vert \D M_2 \vert_{\LL^2}^2 \vert Z\vert_{\LL^2} \vert \D
Z\vert_{\LL^2}
\\ &\leq & k^2 \vert \D M_2\vert_{\LL^2}^2 \vert Z\vert_{\LL^2}^2 +
\frac{k^4}{2 \eta^2} \vert \D M_2 \vert_{\LL^2}^4 \vert Z\vert_{\LL^2}^2
 + \eta^2 \vert \D Z\vert_{\LL^2}^2\,,
 \end{eqnarray*}
 and, for $i=1,2$,
 \begin{eqnarray*}
 II_i  &\leq &   \vert \D M_i\vert_{\LL^2} \vert M_1\vert_{\LL^\infty} \vert
\D Z\vert_{\LL^2} \vert  Z\vert_{\LL^\infty}
 \leq    \vert \D M_i\vert_{\LL^2}  \vert \D Z\vert_{\LL^2} \vert
Z\vert_{\LL^\infty}
 \\  &\leq &   k \vert \D M_i\vert_{\LL^2}  \vert \D Z\vert_{\LL^2}
\vert Z\vert_{\LL^2}^{\frac12} \big( \vert Z\vert_{\LL^2}^{\frac12} +\vert
\D Z\vert_{\LL^2}^{\frac12} \big)
 \\  &\leq &   k \vert \D M_i\vert_{\LL^2}  \vert \D Z\vert_{\LL^2}
\vert Z\vert_{\LL^2}
 +   k \vert \D M_i\vert_{\LL^2}    \vert Z\vert_{\LL^2}^{\frac12} \vert
\D Z\vert_{\LL^2}^{\frac32}
 \\
 &\leq & \frac{k^2}{\eta^2} \vert \D M_i\vert_{\LL^2}^2 \vert
Z\vert_{\LL^2}^2 +\eta^2  \vert \D Z\vert_{\LL^2}^2
 +\frac{k^4}{4\eta^6}\vert \D M_i\vert_{\LL^2}^4\vert Z\vert_{\LL^2}^2 +
\frac34 \eta^2  \vert \D Z\vert_{\LL^2}^2\,.
  \end{eqnarray*}
Hence,
\begin{eqnarray*}
&&\hspace{-2truecm}\lefteqn{I_2(t)=\lb \vert \D M_2\vert^2 M_2- \vert \D
M_1\vert^2 M_1, Z\rb  \leq
 k^2 \Big[  \vert \D M_2\vert_{\LL^2}^2  +  \frac{k^2}{2 \eta^2} \vert
\D M_2 \vert_{\LL^2}^4 }\\ &+&
\sum_{i=1}^2 \frac{1}{\eta^2} \vert \D M_i\vert_{\LL^2}^2
 +\frac{k^2}{4\eta^6}\sum_{i=1}^2 \vert \D M_i\vert_{\LL^2}^4\Big] \vert
Z\vert_{\LL^2}^2
+\frac52\eta^2 \vert \D Z\vert_{\LL^2}^2
  \end{eqnarray*}

Let us note now that by \eqref{indualspace1''}, the 2$^{\rm nd}$ part of the
4$^{\rm th}$ term, i.e. $\lb Z \times\Delta M_1, Z \rb$ is equal to $0$.
  Next, by definition \eqref{indualspace1'}, similarly as the estimate of $II_i$
above, we have the following estimates for the 1st part of the 4$^{\rm th}$ term
using the bound $\vert Z\vert_{L^\infty} \leq 2$, we get
\begin{eqnarray*}
\lb  M_2\times\Delta Z, Z\rb &=& -\lb Z\times \D M_2, \D Z\rb
\leq  \vert Z\vert_{\LL^\infty} \vert \D M_2\vert_{\LL^2} \vert \D
Z\vert_{\LL^2}
\\
&\leq& \frac{k^2}{\eta^2} \vert \D M_i\vert_{\LL^2}^2 \vert
Z\vert_{\LL^2}^2 +\eta^2  \vert \D Z\vert_{\LL^2}^2
 +\frac{k^4}{4\eta^6}\vert \D M_i\vert_{\LL^2}^4\vert Z\vert_{\LL^2}^2 +
\frac34 \eta^2  \vert \D Z\vert_{\LL^2}^2
  \end{eqnarray*}
Therefore, we get the following inequality for the   4$^{\rm th}$ term
\begin{eqnarray*}
&&\hspace{-2truecm}\lefteqn{I_4(t)=\Big[ \lb  M_2(t)\times\Delta Z, Z\rb  -  \lb
Z \times\Delta M_1(t), Z \rb\Big]}\\
&\leq& \frac{k^2}{\eta^2} \vert \D M_i\vert_{\LL^2}^2 \vert
Z\vert_{\LL^2}^2 +\eta^2  \vert \D Z\vert_{\LL^2}^2
 +\frac{k^4}{4\eta^6}\vert \D M_i\vert_{\LL^2}^4\vert Z\vert_{\LL^2}^2 +
\frac34 \eta^2  \vert \D Z\vert_{\LL^2}^2
  \end{eqnarray*}
Next, we will deal with the 3$^{\rm rd}$ term.
Since $\vert M_1\vert_{\LL^\infty} =1$, the H\"older inequality yields
\begin{eqnarray*}
\lb \D M_j M_1 \D Z, Z\rb &\leq & \vert  \D M_j \vert_{\LL^2} \vert
M_1\vert_{\LL^\infty} \vert  \D Z \vert_{\LL^2} \vert Z \vert_{\LL^\infty}
\leq  \D M_j \vert_{\LL^2}  \vert  \D Z \vert_{\LL^2} \vert Z
\vert_{\LL^\infty}
\\
&\leq& \frac{k^2}{\eta^2} \vert \D M_i\vert_{\LL^2}^2 \vert
Z\vert_{\LL^2}^2 +\eta^2  \vert \D Z\vert_{\LL^2}^2
 +\frac{k^4}{4\eta^6}\vert \D M_i\vert_{\LL^2}^4\vert Z\vert_{\LL^2}^2 +
\frac34 \eta^2  \vert \D Z\vert_{\LL^2}^2.
  \end{eqnarray*}
Therefore, we get the following inequality for the   3$^{\rm rd}$ term
\begin{eqnarray*}
I_3(t)=\lb \big( \D M_1+\D M_2)M_1 \D Z, Z\rb &=& \sum_{j=1}^2 \lb
\D M_j M_1 \D Z, Z\rb
\leq  \frac{k^2}{\eta^2} \Big(\sum_{i=1}^2\vert \D M_i\vert_{\LL^2}^2\Big)
\vert  Z\vert_{\LL^2}^2
\\&+&\eta^2  \vert \D Z\vert_{\LL^2}^2
 +\frac{k^4}{4\eta^6}\Big(\sum_{i=1}^2 \vert \D M_i\vert_{\LL^2}^4\Big)
\vert Z\vert_{\LL^2}^2 + \frac32 \eta^2  \vert \D Z\vert_{\LL^2}^2
  \end{eqnarray*}

By inequalities \eqref{ineq-aux-001}, \eqref{ineq-aux-003} and
\eqref{ineq-aux-004}
 we get the following bound for the 5$^{\rm th}$, 6$^{\rm th}$ and
7$^{\rm th}$ terms
\begin{eqnarray*}
 \sum_{i=5,6,8}I_i(t) \leq C_1 \vert Z(t) \vert_{\LL^2}^2.
\end{eqnarray*}
Finally, for the last term we get by \eqref{ineq-aux-003}
\begin{eqnarray*}
 I_7(t) \leq C_1 \vert F(t) \vert \vert Z(t) \vert_{\LL^2}^2.
\end{eqnarray*}
Finally, let us define an $\mathbb{R}$-valued process
\[
\xi_9(t) :=\imp_0^t\sum_{j=1}^3 I_{9,j}(s)\,dW_j(s), \;\; t \in [0,T].
\]
Obviously, $\xi_9$ is an $\LL^2$-valued martingale.
 Next we add together the terms containing $\eta^2 \vert \D
Z\vert_{\LL^2}^2$ to obtain
 \[\frac{19}4\eta^2\vert \D Z\vert_{\LL^2}^2\leq 5 \eta^2\vert \D
Z\vert_{\LL^2}^2\,.\]
Choosing $\eta$ in such a way that $5\eta^2=\frac12$, for a number $C>0$ we
introduce a process
\begin{eqnarray*}
\varphi(t)=\varphi_C(t) &=& C+
k^2 \Big[  \vert \D M_2\vert_{\LL^2}^2  +  \frac{k^2}{2 \eta^2} \vert \D
M_2(t) \vert_{\LL^2}^4
\\ &+&
\sum_{i=1}^2 \frac{1}{\eta^2} \vert \D M_i(t)\vert_{\LL^2}^2
 +\frac{k^2}{4\eta^6}\sum_{i=1}^2 \vert \D M_i(t)\vert_{\LL^2}^4\Big]
 \\
&+& \frac{k^2}{\eta^2} \vert \D M_i(t)\vert_{\LL^2}^2
 +\frac{k^4}{4\eta^6}\vert \D M_i(t)\vert_{\LL^2}^4
 + \frac{k^4}{4\eta^6}\sum_{i=1}^2 \vert \D M_i(t)\vert_{\LL^2}^4, \;\; t
\in [0,T].
  \end{eqnarray*}
From all our inequalities we infer that there exist a constant $C>0$ such that
\begin{eqnarray}\label{ineq-Z^2}
\vert Z(t)\vert_{\LL^2}^2 \leq \imp_0^t \varphi_C(s) \, \vert
Z(s)\vert_{\LL^2}^2\, ds + \xi_9(t), \;\ t\in [0,T]
  \end{eqnarray}
By the It\^o Lemma applied to the following an $\mathbb{R}$-valued
process,
\[ Y(t):=\vert Z(t)\vert_{\LL^2}^2e^{-\imp_0^t \varphi_C(s)  \, ds}, \;\; t \in
[0,T],
\]
 see
\cite{Schmalfuss_1997} for a similar idea, we infer that
\begin{eqnarray*}
Y(t) &\leq& \imp_0^{t} e^{-\imp_0^{s} \varphi_C(r)  \, dr} \, d\xi_9(s)
\\&=&
\sqrt{\eps} \sum_{j=1}^3 \imp_0^{t} e^{-\imp_0^{t} \varphi_C(s)  \, ds}
\lb  G\la M_2(s) \ra- G\la M_1(s) \ra \big) \newe_j, Z\rb \,dW_j(s), \;\; t \in
[0,T].
  \end{eqnarray*}
Since $M_1$, $M_2$ and $Z$ are uniformly bounded and $G$ is locally Lipschitz
the process defined by the RHS of the last inequality is an
$\mathbb{F}$-martingale.

Thus, we infer that
\[
\mathbb{E} Y(t) \leq 0, \;\; t \in [0,T]\,,
\]
and since $Y$ is  nonnegative, we deduce that
$Y(t) =0$, $\mathbb{P}$-a.s.,  for every $ t \in [0,T]$.
Finally, the definition of $Y$ yields
\[Z(t) =0\quad\mathbb{P}-a.s.,\quad\mathrm{for\,\,every}\quad t \in [0,T]\,.\]
This completes the proof.
\end{proof}
\begin{remark}\label{rem-unique} Let us note first that the processes $M_i$,
$i=1,2$ in
Theorem  \ref{thm-uniqueness pathwise weak} satisfy
weaker conditions that those guaranteed by the existence result from Theorem
\ref{thm-weak solution}.
Hence  our uniqueness result in Theorem  \ref{thm-uniqueness pathwise weak}  holds in the following
sense. \\
Suppose that $M_1$ is a solution satisfying assumptions of  Theorem
\ref{thm-uniqueness pathwise weak} and $M_2$
a solution in the sense of
Theorem \ref{thm-weak solution}, both defined on the same filtered probability space,
then $M_1=M_2$.
\end{remark}
\section{Further regularity}\label{sec:furtherreg}
In this section, we assume that
a system
$ \Bigl( \Omega,\mathcal{F},\mathbb{F},\mathbb{P},W, M\Bigr)
$
is a  weak martingale solution to problem
(\ref{eqn-main})
such that $M$ has paths in the space $S_T$ defined by
\begin{equation}\label{eqn-S}
S_T:=C([0,T];\mathbb{H})\cap L^4(0,T;\mathbb{H}^1).
\end{equation}
 Some regularity properties of $M$ are listed in
Theorem~\ref{thm-weak solution}.
The main result of this section is Theorem~\ref{thm-reg}, where we prove
stronger regularity of the solution.
In Proposition~\ref{contsinH1}, we use this estimate to show that paths of $M$
lie
in $C([0,T];\mathbb{H}^1)$, $\mathbb{P}$-almost everywhere; this improves upon
the continuity property in
Theorem~\ref{thm-weak solution}.
\par\noindent We start with a lemma that expresses $M$ in a mild-form which
allows us to exploit the
regularizing properties of the semigroup $(e^{-tA})$. The proof of this well
known fact is omitted, see for instance
\begin{lemma}\label{lem-mild}
For each $t\in [0,T]$, $\P$-a.s.
\begin{eqnarray}\label{eq-mild}
M(t)=e^{-\alpha tA}M_0&+&\imp_0^te^{-\alpha(t-s)A}(M(s)\times\Delta
M(s))ds+\alpha\imp_0^te^{-\alpha(t-s)A}\la|\D M(s)|^2M(s)\ra\,ds
\nonumber\\
&+&\eps^{\frac{1}{2}}\imp_0^te^{-\alpha(t-s)A}G(M(s))B\,d W(s)
\\
&-&\beta\imp_0^te^{-\alpha(t-s)A}G(M(s))F(M(s))\,ds
\nonumber\\
&+&\frac{\eps}{2}\sum_{i=1}^3\imp_0^t
e^{-\alpha(t-s)A}G^\prime(M(s))\newe_iG(M(s))\newe_i\,ds.
\nonumber
\end{eqnarray}

\end{lemma}
\begin{comment}
\begin{proof} The proof of this lemma is much the same as the first part of the
proof of Theorem 9.15 in~\cite{PZ}. We outline the idea: for each $t\in [0,T]$
and $j\ge 1$ and $n\geq j$, one can write down the It\^{o} formula for the
process $(\langle \mathrm{e}^{-\alpha(t-s)A}e_j, \pi_n
y(s)\rangle_{\mathbb{L}^2})_{s\in[0,t]}$ (here $\pi_n$ is orthogonal projection
onto $\mathbb{H}_n$). One sets $s=t$ in this formula and obtains
\[
\langle y(t) , e_j\rangle_{\mathbb{L}^2}=\langle\mbox{ right hand side
of~(\ref{eq-mild})}, e_j\rangle_{\mathbb{L}^2}
\]
$\mathbb{P}$-almost everywhere.  Since $j$ is arbitrary, this completes the
proof.
\end{proof}
\end{comment}
Before we state the main result of this section let us make the following
important remark.
\begin{remark}\label{rem-ident}
\deld{By estimate \eqref{sun1}} Suppose that  the vector $\Delta M(t,x)\in\R^3$ is
a.e. well defined  and that \[ |M(t,x)|^2=1 \mbox{ a.e.}.\] Then we infer that
\[M(t,x)\cdot\Delta M(t,x)=-|\D M(t,x)|^2,\quad \deld{t,x}\,a.e.\]
and therefore,  an elementary identity
\[|a\times b|^2+|a\cdot b|^2=|a|^2\cdot|b|^2,\quad a,b\in\R^3,\]
yields
\[|M(t,x)\times\Delta M(t,x)|^2+|\D M(t,x)|^4=|\Delta M(t,x)|^2,\quad
t,x-a.e.\]
\end{remark}

\begin{theorem}\label{thm-reg}
Assume that $p\in [1,\infty)$. Then for every  $M_0\in\HH^1$ \deld{,
$F\in\mathcal P_T$}
and $\newe=\big(\newe_i\big)_{i=1}^3\in(\mathbb H^1)^3$  Then there exists a
constant
$C_p=C_p\la\alpha,T,\deld{\Vert F\Vert ,}\Vert M_0\Vert_{\HH^1} ,\vert
\newe\vert_{\HH^1} \ra$
such that the unique solution $M$ of the problem  \eqref{eqn-main} satifies
\begin{equation}\label{ineq-maximal}
\E\la\imp_0^T |\D M(t)|^4_{\LL^4}\,dt  + \imp_0^T \left|\Delta
M(t)\right|^2_{\mathbb{L}^2}\,dt\ra^{p}\leq C_p.
\end{equation}
\end{theorem}
\begin{definition}\label{def-strong martingale solution}
A weak martingale solution
\begin{equation}\label{eqndef-weak mart sol}
 \Bigl( \Omega,\mathcal{F},\mathbb{F},\mathbb{P},W, M\Bigr)
\end{equation}
to problem \eqref{eqn-main} is called a martingale strong solution to to problem \eqref{eqn-main} iff it satisfies condition \eqref{ineq-maximal} for $p=1$.
\end{definition}
\begin{proof} By the uniqueness it is sufficient to prove the theorem for
the solution constructed in
Theorem \ref{thm-weak solution}.
Let us describe the structure of the proof. In Step 1 we will show the
first part of inequality
\eqref{ineq-maximal} for every $p\in [1, \infty)$. In Step 2 we will show the
second part of
inequality \eqref{ineq-maximal} for every $p=1$. In step 3 we will use Step 2
and Remark \ref{rem-ident}
to deduce the second part of inequality \eqref{ineq-maximal} for every $p\in [1,
\infty)$.

We will use repeatedly the following well known properties of the semigroup $\la
e^{-tA}\ra$. \par\noindent
 The semigroup $\la e^{-tA}\ra$, where $A$ is defined in \ref{Adefinition}, is
ultracontractive, see, for example,~\cite{Arendt}, that is, there exists $C>0$
such that if $1\le p\le q\le\infty$, then
\begin{equation}\label{ultra}
|\mathrm{e}^{-tA}f|_{\LL^q}\le
\frac{C}{t^{\frac{1}{2}\left(\frac{1}{p}-\frac{1}{q}\right)}}|f|_{\LL^p},\quad f
\in\LL^p,\,\, t>0.
\end{equation}
\begin{comment}
and
\begin{equation}\label{ultra-grad}
|\D\mathrm{e}^{-tA}f|_{\LL^q}\le
\frac{C}{t^{\frac{1}{2}+\frac{1}{2}\left(\frac{1}{p}-\frac{1}{q}\right)}}
|f|_{\LL^p},\quad f\in\LL^p,
\end{equation}
\end{comment}
It is also well known that $A$ has maximal regularity property, that is, there
exists $C>0$ such that for any $f\in L^2\la 0,T;\mathbb H\ra$ and
\[u(t)=\imp_0^te^{-(t-s)A}f(s)ds,\quad t\in [0,T],\]
we have
\begin{equation}\label{eqn-maxr}
\imp_0^T\left|Au(t)\right|_{\mathbb H}^2\,dt\le C\imp_0^T|f(t)|_{\mathbb
H}^2\,dt.
\end{equation}

Let us fix for the rest of the proof $T>0$, and an auxiliary number
$\delta\in\la\frac{5}{8},\frac{3}{4}\ra$. Let us also fix $\rho>0$ and $r>0$ such that $\Vert M_0\Vert_{\mathbb{H}^1} \leq \rho$ and  $\Vert e \Vert_{\mathbb{H}^1} \leq r$.
\begin{itemize}
\item{\textbf{Step 1}} Let us additionally choose and fix $ p \in [1,\infty)$.
 By Lemma \ref{lem-mild} $M$ can be written as a sum of six
terms:
\[M(t)=\sum_{i=0}^5m_i(t),\]
and we will consider each term separately.  In what follows, $C$ stands for a
generic constant that depends on $p$, $T$, $\alpha$, $\rho$ and $r$ only.
In order to simplify notation, we put, without loss of generality,
$\eps=\alpha=\beta=1$.\\
We will show first that
\begin{equation}\label{sun1}
\E \Bigl(\imp_0^T|M(t)|^4_{\mathbb W^{1,4}}\,dt\Bigr)^{p}\le
C\la p,T,\alpha,\rho,r\ra.
\end{equation}
Since $\delta >\frac58$ so that the Sobolev imbedding $\mathbb
X^\delta\hookrightarrow\mathbb W^{1,4}$, it is sufficient to prove
the following stronger estimate:
\begin{equation}\label{sun11}
\E \Bigl(\imp_0^T\left|A^\delta M(t)\right|^4_{\mathbb{L}^2}\,dt\Bigr)^{p}\le
C\la p,T,\alpha,\rho,r\ra.
\end{equation}

We start with $m_0$. For each $t\in(0,T]$, we have
\[\left|A^\delta e^{-
tA}M_0\right|^4_{\mathbb{L}^2}\le\frac{C}{t^{4\delta-2}}\left|M_0\right|^4_{
\mathbb{H}^1},\]
and therefore, since $\delta <\frac34$, we infer that
\begin{equation}\label{m0}
\imp_0^T\left|A^\delta m_0(t)\right|^4_{\mathbb{L}^2}\, dt\le
C\left|M_0\right|^4_{\mathbb H^1}.
\end{equation}
We will consider $m_1$. Putting $f=M\times\Delta M$ we have
\[
|A^\delta e^{-(t-s)A}f(s)|_{\mathbb{L}^2}\le
C(t-s)^{-\delta}|f(s)|_{\mathbb{L}^2},\quad 0<s<t<T,
\]
hence applying the  Young inequality we obtain
\begin{eqnarray*}
\imp_0^T|A^\delta m_1(t)|^4_{\mathbb{L}^2}\,dt&\le&
C\imp_0^T\left(\imp_0^t(t-s)^{-\delta}|f(s)|_{\mathbb{L}^2}\,ds\right)^4\,dt\\
&\le&
C\left(\imp_0^Ts^{-\frac{4\delta}{3}}\,ds\right)^3\left(\imp_0^T|f(s)|^2_{\mathbb{L}^2}\
,ds\right)^2.
\end{eqnarray*}
Thereby, since $\frac{4\delta}{3}<1$, part (3) of Theorem~\ref{thm-weak solution} yields
\begin{equation}\label{m1}
\E\Bigl( \imp_0^T|A^\delta m_1(t)|^4_{\mathbb{L}^2}\,dt\Bigr)^{p}\le
C(2p,T,\alpha,\rho,r).
\end{equation}
Since for every $t\in [0,T]$, $|M(t,x)|=1$ almost everywhere, and $e_i\in\mathbb
H^1$, $i=1,2,3$, the estimate \eqref{gprime} implies that there exists deterministic $c>0$
such that
\[\sum_{i=1}^3|G(M)e_i|_{\LL^2}+\sum_{i=1}^3|G^\prime(M)e_iG(M)e_i|_{\LL^2}\le c.\]
Therefore, the same arguments as for $m_1$ yield
\begin{equation}\label{m5}
\E\Bigl(\imp_0^T|A^\delta m_5(t)|^4_{\mathbb{L}^2}\,dt\Bigr)^{p} \le
C(p,T,\alpha,\rho,r).
\end{equation}
We will now consider the  term $m_2$ using the fact that  $f=|\D
M|^2M\in\LL^\infty(0,T;\LL^1)$. Invoking the semigroup property of $e^{-tA}$ and the
ultracontractive estimate \eqref{ultra} with $p=1$ and $q=2$ we find that there
exists $C>0$ such that  $\P$-a.s.
\[|A^\delta e^{-(t-s)A}f(s)|_{\mathbb{L}^2}\le
\frac{C}{(t-s)^{\delta+\frac{1}{4}}}\sup_{r \in [0,T]}|M\la r\ra|^2_{\HH^1},\quad
0<s<t\in [0,T].\]
Therefore,
\[\imp_0^T\left|\imp_0^tA^\delta
e^{-(t-s)A}f(s)\,ds\right|^4_{\mathbb{L}^2}\,dt\le C \vert f\vert^4_{L^\infty(0,T;\mathbb{L}^1)}
\imp_0^T\la\imp_0^t\frac{ds}{(t-s)^{\delta+\frac{1}{4}}}\,ds\ra^4\,dt.\]
Hence (since $\delta+\frac{1}{4}<1$) Theorem~\ref{thm-weak solution} yields
\begin{equation}\label{m2}
\E\imp_0^T\left|A^\delta m_2(t)\right|^4_{\mathbb{L}^2}\,dt\leq C \vert DM\vert^8_{L^\infty(0,T;\mathbb{L}^2)} \leq  C\la
T,\rho,r\ra.
\end{equation}
In order to estimate $m_3$ we recall that there exist $a_r>0$ such that
\begin{equation}\label{ds}
\Vert G(M)e_i\Vert_{\mathbb{H}^1}\le a_r(1+\Vert M\Vert_{\mathbb{H}^1}^2),\;\; i=1,2,3.
\end{equation} Invoking Lemma 7.2 in
\cite{DZ}
we find that for $i=1,2,3$ and any $t\in [0,T]$
\begin{eqnarray*}
&&\!\!\!\!\!\!\!\!\!\!\!\!\!\!\E\left|\imp_0^tA^\delta
e^{-(t-s)A}G(M(s))e_i\,dW(s)\right|^4_{\mathbb{L}^2}\\
&\le&C(T)\E\left(\imp_0^t|A^\delta
e^{-\alpha(t-s)A}G(M(s))e_i|^2_{\mathbb{L}^2}\,ds\right)^2\\
&=&C(T)\E\left(\imp_0^t|A^{\delta-\frac{1}{2}}
e^{-(t-s)A}A^\frac{1}{2}G(M(s))e_i|^2_{\mathbb{L}^2}\,ds\right)^2\\
&\le&C(T)\E\left(\imp_0^t\frac{|G(M(s))e_i|^2_{\mathbb{H}^1}}{(t-s)^{2\delta-1}}\,
ds\right)^2\\
&\le&C(T) \E\sup_{r \in [0,T]}\bigl[1+ |M\la r\ra|^8_{\mathbb{H}^1} \bigr].
\end{eqnarray*}
Thus, Theorem~\ref{thm-weak solution} now yields
\begin{equation}\label{m4}
\E\imp_0^T|A^\delta m_3(t)|^4_{\mathbb{L}^2}\,dt\le C(T,\rho,r).
\end{equation}
Because by inequality \eqref{eqn-G-polynomial growth-H^1-norm},
\[
\vert G(M)F(M)\vert_{\mathbb{H}^1} \leq aC  (1+\Vert M\Vert_{\mathbb{H}^1}^2)\Vert M\Vert_{\mathbb{H}^1}^2
\]
the case of $m_4$ can be treated very easily.

Finally, combining estimates \eqref{m0} to \eqref{m4} we obtain \eqref{sun11}
and \eqref{sun1} follows.

\item{\textbf{Step 2}}
We will prove that
\begin{equation}\label{eqn-maxr}
\E\imp_0^T|AM(t)|^2_{\mathbb{L}^2}\,dt\le C(T,\rho,r).
\end{equation}
To this end we note first that using the maximal inequality \eqref{eqn-maxr} and the
first part of the proof it is easy to see that
\begin{equation}\label{maxm1}
\E\imp_0^T\left|Am_i(t)\right|_{\mathbb{L}^2}^2\,dt\le
C(T,\rho,r),\quad i=1,2,4.
\end{equation}
The estimate
\begin{equation}\label{maxm0}
\imp_0^T\left|Am_0(t)\right|^2_{\mathbb{L}^2}\,dt\le C\la T,\rho\ra,
\end{equation}
is an immediate consequence of the fact that $M_0\in\HH^1=D\la A^{1/2}\ra$ .\\
We will consider now the stochastic term $m_3$. Using \eqref{ds}, a result of
Pardoux in \cite{Pardoux_1979} and part 1 of Theorem \ref{thm-weak solution} we find
that
\begin{equation}\label{maxm4}
\begin{aligned}
\E\imp_0^T|Am_3(t)|_{\mathbb{L}^2}^2dt&\leq C\E\imp_0^T\la
|M(t)|^4_{\mathbb{H}^1}+1\ra\,dt
\le C(T,\rho,r).
\end{aligned}
\end{equation}
Combining \eqref{maxm1}, \eqref{maxm0} and \eqref{maxm4} we obtain \eqref{eqn-maxr}.

\item{\textbf{Step 3}} Take $p \geq 1$. By Step 2 and Remark we infer that

\begin{eqnarray}\label{ineq-maximal-part2}
&&\E\la\imp_0^T  \left|\Delta M(t)\right|^2_{\mathbb{L}^2}\,dt\ra^{p} =
\E\la\imp_0^T  \int_{\domain}
\left|\Delta  M(t.x)\right|^2\,dx\,dt\ra^{p}\\
&=& \E\la\imp_0^T  \int_{\domain} |M(t,x)\times\Delta M(t,x)|^2\,dx\,dt +
\imp_0^T  \int_{\domain}|\D M(t,x)|^4\,dx\,dt\ra^{p}.
\nonumber
\end{eqnarray}
Hence the second part of inequality \eqref{ineq-maximal} in Theorem
\ref{thm-reg} folllows from the first part (proved above in Step 2) and
inequality \eqref{ineqn-projected Laplacian} from the Theorem \ref{thm-weak solution}
about the existence of weak solutions.

\end{itemize}
 The proof is complete.

\end{proof}
\begin{proposition}\label{contsinH1}
$\mathbb{P}$ almost surely, the  paths of $M$ lie in the space $C([0,T];\mathbb{H}^1)$.
\end{proposition}
\begin{proof}
The proposition follows easily from the results in \cite{Pardoux_1979}.
\end{proof}

\begin{corollary}\label{strong}
Let $\newe_i\in\mathbb \mathbb{H}^{1}$, $i=1,2,3$. Let $W$ be an
$\mathbb  F$ Wiener process defined on the probability space
$\la\Omega,\mathcal F,\mathbb F,\P\ra$. Then, for every $M_0\in\HH^1$ and $\eps>0$, there exists a unique pathwise
solution $M^\eps\in C\la [0,T];\HH^1\ra\cap L^2\la 0,T;D(A)\ra$ of the
problem \eqref{eqn-main}, i.e.
\begin{equation}
\label{eqn-ME1b}\begin{aligned}
M(t)&=M_0+\alpha\imp_0^t\Delta M(s)\,ds+\alpha\imp_0^t|\D M(s)|^2M(s)\,ds+\imp_0^t
M(s)\times\Delta M(s)\,ds\\
&+\frac{\eps}{2}\sum_{i=1}^3\imp_0^tG^\prime(M(s))\newe_iG(M(s))\newe_i\,ds+\sqrt{\eps}
\imp_0^tG(M(s))B\,dW(s)\\
&-\beta\imp_0^tG(M(s))F(M(s))\,ds,
\end{aligned}
\end{equation}
where all the integrals are the Bochner or the It\^o integrals in
$\LL^2$.
\end{corollary}
In what follows we will denote by $\mathcal X_T$ the Banach space
\begin{equation}\label{eqn-X_T}
\mathcal X_T=C\la [0,T];\mathbb H^1\ra\cap L^2(0,T;D(A)).
\end{equation}

By an infinite-dimensional version of the Yamada and Watanabe Theorem, see
\cite[Theorems 12.1 (part 3) and 13.2]{ondrejat}, the   pathwise uniqueness and
the existence of weak solutions
implies uniqueness in law and the existence of a strong solution. In
Theorem~\ref{thm-YW2} below,
we state such a result for equation (\ref{eqn-main}).

Using the additional regularity results proven in this section,  we have the following result.

\begin{theorem}\label{thm-YW2} Let assumptions of Theorem \ref{thm-uniqueness pathwise weak} be
satisfied. Then  uniqueness in law and the existence of a strong solution holds
for equation (\ref{eqn-main}) in the following sense:
\begin{enumerate}
\item if $\Bigl( \Omega,\mathcal{F},\mathbb{F},\mathbb{P},W, M\Bigr)$ and
$\Bigl(
\Omega^\prime,\mathcal{F}^\prime,\mathbb{F}^\prime,\mathbb{P}^\prime,W^\prime,
M^\prime\Bigr)$
are two martingale strong solutions to problem (\ref{eqn-main}) such that both
 $M$ and $M^\prime$ are $\mathcal X_T$-valued random variables,
then $M$ and $M^\prime$ have the same laws on $\mathcal X_T$;
\item for every $\eps>0$ there exists a Borel measurable function
\begin{equation}\label{eqn_map J^eps}
 J^\eps: \fourIdx{}{0}{}{}{C  \la[0,T];\R^3\ra}:=\bigl\{ \omega\in  C\la
[0,T];\R^3\ra: \omega(0)=0\bigr\}
 \to \mathcal X_T  \end{equation}
such that for any filtered probability space
$({\Omega},{\mathcal{F}},{\mathbb{F}},{\mathbb{P}})$,
where the filtration ${\mathbb{F}}=({\mathcal{F}}_t)$ is such that
${\mathcal{F}}_0$
contains all ${\mathbb{P}}$-null sets from ${\mathcal{F}}$,
and for any
  $\mathbb{R}^3$-valued ${\mathbb{F}}$-Wiener process
${W}=( W(t))_{t\in[0,T]}$, the system $\Bigl(
{\Omega},{\mathcal{F}},{\mathbb{F}},
 {\mathbb{P}},{W}, {M}^\eps\Bigr)$, where  ${M}^\eps=J^\eps\circ {W}$, i.e.
  \[
  {M}^\eps: {\Omega} \ni \omega \mapsto J^\eps\bigl({W}(\omega)\bigr)\in
\mathcal X_T,
 \]
 is a strong  martingale  solution\footnote{In particular,
${M}$ is ${\mathbb{F}}$-progressively measurable.} to problem
(\ref{eqn-main}).
\end{enumerate}
\end{theorem}

%%%************************************************************************************************

\section{The Large Deviations Principle}\label{sec-LDP}

In this section we will prove the large deviation principle for the family of
laws of the solutions $M^\eps$  of equation (\ref{eqn-main}) with the
parameter $\eps\in(0,1]$ approaching zero and fixed $M_0\in \mathbb{H}^1$.

In what follows we will denote by $M^\eps$ the unique strong martingale solution to the problem  \eqref{eqn-main}.

The main result in this section is as follows.

\begin{theorem}\label{thm-large deviations}
The family of laws $\{\mathcal{L}(M^\eps):\,\eps\in(0,1]\}$ on $\mathcal
X_T$ satisfies the large deviation principle with rate function $I$ defined below in equation \eqref{eqn-rate function}.
\end{theorem}

Before we embark on the proof of the above result we will present the necessary background. In particular we will formulate crucial Lemmata  \ref{lem19} and \ref{lem22}. Then we will present the proof of Theorem \ref{thm-large deviations}. This will be followed by the proof  of  Lemma \ref{lem19}. The proof of Lemma \ref{lem22} will be given at the very end of this section.

In order to prove the above result we will present some results due to Buhhiraja and Dupuis \cite{Budhiraja+Dupuis_2000}. Following that
paper we will  formulate some two general claims. These claims will be consequence of Lemmata  \ref{lem19} and
\ref{lem22} which we first only formulate. This preliminary material will   be followed by the proof of Theorem
\ref{thm-large deviations}.
This will then be followed by the proof  of  Lemma \ref{lem19}. The
proof of Lemma \ref{lem22} will be given at the very end of this section.

\subsection{Large Deviations Principle  according to Buhhiraja and Dupuis}
\label{subsec-LDP}

In order to prove the Large Deviations Principle formulated in Theorem \ref{thm-large deviations} holds we need to
consider an equation slightly more general than equation \eqref{eqn-main}.

Let $\la\Omega,\mathcal{F},\mathbb{F},\mathbb{P}\ra$, be the classical Wiener space, i.e.
\begin{eqnarray*}
\Omega&=&\fourIdx{}{0}{}{}{C}\la[0,T];\R^3\ra,
\\
\mathbb{P}&&\mbox{ is the Wiener measure on }\Omega,
\\
W&=&\bigl(W(t)=W_t\bigr)_{t\in [0,T]} \mbox{ is the canonical $\mathbb{R}^3$-valued Wiener process on $(\Omega,\mathbb{P})$},
\\
\mathbb{F}&=&\la\mathcal{F}_t\ra_{t\in
[0,T]} \mbox{ is the $\mathbb{P}$-completion of the  natural filtration
$\mathbb{F}^0=\la\mathcal F_t^0\ra_{t\in [0,T]}$ generated by
$W$.}
\end{eqnarray*}
Note that filtration ${\mathbb{F}}=({\mathcal{F}}_t)$ is such that
${\mathcal{F}}_0$
contains all ${\mathbb{P}}$-null sets from ${\mathcal{F}}$.

By Theorem \ref{thm-YW2} for every $\eps>0$ there exists a Borel map
\begin{equation}\label{eqn_map J^eps-2}
 J^\eps: \fourIdx{}{0}{}{}{C  \la[0,T];\R^3\ra}
 \to \mathcal X_T  \end{equation}
 the system $\Bigl(
{\Omega},{\mathcal{F}},{\mathbb{F}},
 {\mathbb{P}},{W}, {M}^\eps\Bigr)$, where    \[
  {M}^\eps: {\Omega} \ni \omega \mapsto J^\eps\bigl({W}(\omega)\bigr)\in
\mathcal X_T,
 \]
 is a  martingale strong solution\footnote{In particular,
${M}$ is ${\mathbb{F}}$-progressively measurable.} to problem
(\ref{eqn-main}).\\
By  $\E\deld{_\P}$ we will denote the  integration with respect to the measure $\P$.

 Suppose that $\rX$ is a separable Banach space. We  say that an
$\rX$-valued $\mathbb{F}$-predictable process
$h:[0,T]\times \fourIdx{}{0}{}{}{C\la [0,T];\R^3\ra}\to\rX$  belongs to
$\mathcal P_T(\rX)$ iff
\begin{equation}\label{norm-h}
\|h\|_T^2=\esssup_{\omega\in \Omega }
\imp_0^T|h(t,\omega)|_{\rX}^2\,dt<\infty.
\end{equation}

Given $h\in\mathcal P_T(\mathbb{R}^3)$ we can consider an equation
\begin{equation}\label{eqn-ME1-h-new}
\begin{aligned}
M(t)=M_0&+\imp_0^t[M\times\Delta M-\alpha M\times(M\times\Delta M)]\,ds \\
&+\sqrt{\eps}\imp_0^tG(M)B\,dW(s)+\frac{\eps}{2}\sum_{i=1}^3\imp_0^t\, \big[G^\prime(M)\newe_i\big]\big(G(M)\newe_i)\,ds\\
&-\beta\imp_0^tG(M)f(M)\,ds+\imp_0^t G(M)Bh(s)\,ds,\quad t\in [0,T].
\end{aligned}
\end{equation}

\begin{theorem}\label{thm-with h}
Assume that ${h}\in\mathcal P_T(\mathbb{R}^3)$ and $\eps \in (0,1]$.
Then there exists a process $\widetilde M=M^{h,\eps}$ such that the system
\[
\Bigl( \Omega,\mathcal{F},\mathbb{F},\mathbb{P},W,  \widetilde M \Bigr)
\]
is a strong martingale  solution of problem \eqref{eqn-ME1-h-new} such that
for every $p \geq 1$,
\begin{eqnarray}
&& \E \sup_{t\in [0, T]}|\widetilde M(t)|_{\mathbb H^1}^{p}<\infty,
\\
&& \E\la\imp_0^T |\D \widetilde  M(t)|^4_{\LL^4}\,dt  + \imp_0^T \left|\Delta
\widetilde M(t)\right|^2_{\mathbb{L}^2}\,dt\ra^{p}< \infty.
\end{eqnarray}

\end{theorem}
\begin{proof}\textbf{Part I: the existence} Let us fix $\eps>0$.  For any  ${h}\in\mathcal P_T(\mathbb{R}^3)$  let us put
\begin{equation}\label{eqn-Girsanov density}
\tilde \rho_{h}=
\exp\la -\frac{1}{\sqrt{\eps}}\imp_0^T{h}(s)\,dW(s)+\frac1{2\eps}\imp_0^T|{h}(s)|^2\,ds\ra\,.
\end{equation}
and
\begin{equation}\label{eqn-Girsanov Wiener process}
\widetilde W_{h}(t)=W(t)+\frac{1}{\sqrt{\eps}} \imp_0^t {h}(s)\,ds,\quad t\in[0,T].
\end{equation}
Since ${h}\in\mathcal P_T(\mathbb{R}^3)$ we infer that
\[\E \bigl(\tilde \rho_{h}\bigr)^2<\infty\]
and therefore  there exists a probability measure $\tilde \P_{h}$ on ${\mathcal F}_T$ such that
\[\frac{d\P_{h}}{d\P}=\tilde \rho_{h}.\]
Invoking the Girsanov Theorem we find that the process $\widetilde W_h$ is a Wiener process on probability space
$(\Omega,\P_{h})$. Note that now $\Omega=\fourIdx{}{0}{}{}{C}\la[0,T];\R^3\ra$.  Therefore, by part (ii) of Theorem \ref{thm-YW2}, if the process
$\widetilde M$ is defined by
\[
\widetilde M: \Omega \ni \omega \mapsto J^\eps\bigl(\widetilde W_{h} (\omega)\bigr) \in \mathcal{X}_T
\]
then the system
\[
\Bigl( \Omega,\mathcal{F},\mathbb{F},\mathbb{P}_{h},\widetilde{W}_h, \widetilde M\Bigr)
\]
is a strong martingale  solution of problem \eqref{eqn-main}. In particular, by Theorems \ref{thm-weak solution} and \ref{thm-reg},
\begin{eqnarray*}
&&\widetilde \E \sup_{t\in [0, T]}|\widetilde M(t)|_{\mathbb H^1}^{2p}<\infty,\;\; p \geq 1,
\\
&&\widetilde \E\la\imp_0^T | \D \widetilde M(t)|^4_{\LL^4}\,dt  + \imp_0^T \left|\Delta
\widetilde M(t)\right|^2_{\mathbb{L}^2}\,dt\ra^{2p}< \infty,\;\; p \geq 1.
\end{eqnarray*}

On the other hand, since ${h}\in\mathcal P_T(\mathbb{R}^3)$ we infer that
\begin{equation}
\label{ineq-Girsanov density}
\widetilde \E \bigl(\tilde \rho_{h}\bigr)^{-2}<\infty
\end{equation}
and therefore  $\mathbb{P}$ is absolutely continuous w.r.t. $\mathbb{P}_{h}$ and
\[\frac{d\P}{d\P_{h}}=\tilde \rho_{h}^{-1}.\]
Therefore, by applying the H\"older inequality we infer that for any $p\geq 1$.
\begin{eqnarray*}
&&\E \sup_{t\in [0, T]}|\widetilde M(t)|_{\mathbb H^1}^{p}<\infty,
\\
&& \E\la\imp_0^T | \D \widetilde M(t)|^4_{\LL^4}\,dt  + \imp_0^T \left|\Delta
\widetilde M(t)\right|^2_{\mathbb{L}^2}\,dt\ra^{p}< \infty.
\end{eqnarray*}
Therefore, by a standard argument, we infer that
the system
\[
\Bigl( \Omega,\mathcal{F},\mathbb{F},\mathbb{P},W, \widetilde M\Bigr)
\]
is a strong martingale  solution of problem \eqref{eqn-ME1-h-new}, see e.g. Appendix A in \cite{DaPrato+F+P+R_2013}.
This completes the existence proof. \\
\textbf{Part II: uniqueness.}  The uniqueness follows from Theorem \ref{thm-uniqueness pathwise weak}.

\end{proof}

Let us note that that we have used the  Girsanov Theorem only to prove the
existence of a solution to problem \eqref{eqn-ME1-h-new}. Having this done, we return to
our fixed probability space with a fixed Wiener process.
\par\medskip\noindent
Let now define a Borel map
\[
J^0: \fourIdx{}{0}{}{}{C}([0,T];\mathbb{R}^3) \to \mathcal{X}_T.
\]
If $x\in \fourIdx{}{0}{}{}{C}([0,T];\mathbb{R}^3)
\setminus \fourIdx{}{0}{1,2}{}{H}([0,T];\mathbb{R}^3)$, then we put $J^0(x)=0$. If
  $x^\prime=h$ for some $h\in L^2\la 0,T;\R^3\ra$, then by Corollary \ref{strong} there exist a unique  function $y_h \in \mathcal{X}_T$ that is
the unique solution of the equation
\begin{eqnarray}\label{eqn-deterministic}
y_h(t)=M_0&+&
\imp_0^t[y_h(s)\times\Delta y_h(s)-\alpha y_h(s)\times(y_h(s)\times\Delta y_h(s))]\,ds
\\
&-&\beta \imp_0^tG\la y_h\ra f\la y_h\ra ds+\imp_0^tG\la y_h\ra Bh(s)\,ds,
\nonumber
\end{eqnarray}
where the mapping $G$ has been defined  in \eqref{eqndef_function G}. We put
\[
J^0 (x):= y_h, \;\; x=\int_0^{\cdot}h(s)\,ds, h\in L^2\la0,T;\R^3\ra.\]
Let us note that in view of Lemma , problem can be written in the following equivalent form
\begin{eqnarray}\label{eqn-deterministic-2}%{e108}
y_h(t)=M_0&+&
\alpha\imp_0^t\Delta y_h\,ds+\alpha\imp_0^t\left|\D
y_h\right|^2y_h\, ds+\imp_0^ty_h\times\Delta y_h\,ds
\\
&-&\beta \imp_0^tG\la y_h\ra f\la y_h\ra ds+\imp_0^tG\la y_h\ra Bh(s)\,ds,
\nonumber
\end{eqnarray}

{We can easily prove that
 the map $J^0: \fourIdx{}{0}{1,2}{}{H}([0,T];\mathbb{R}^3) \to \mathcal{X}_T$ is continuous.
Since $\fourIdx{}{0}{1,2}{}{H}([0,T];\mathbb{R}^3)$ is a Borel subset of $ \fourIdx{}{0}{}{}{C}([0,T];\mathbb{R}^3)$,
 we infer that
 the map $J^0: \fourIdx{}{0}{}{}{C}([0,T];\mathbb{R}^3) \to \mathcal{X}_T$ is Borel measurable.
}
We  define now the rate function $I:\mathcal X_T\to[0,\infty]$ by the
formula
\begin{equation}\label{eqn-rate function}
I(u):=\inf\left\{\frac{1}{2}\imp_0^T |h(s)|^2\,ds:\,h\in L^2\la
0,T;\R^3\ra\mbox{ and }u= J^0(\int_0^{\cdot}h(s)\, ds)\right\},
\end{equation}
where   $\inf
\emptyset=\infty$.
\deld{the infimum of the empty set is equal to the infinity.}\\

In order to prove Theorem \ref{thm-large deviations}, i.e.  that the family of laws
$\{\mathcal{L}(J_0^\eps(W)):\,\eps\in(0,1]\}$ satisfies the large deviation
principle on $\mathcal{X}_T$ with the rate function
$I$ we will follow the weak convergence method of Budhiraja and
Dupuis~\cite{Budhiraja+Dupuis_2000}, see also Duan and
Millet~\cite{Duan+Mil_2009} and Chueshov and Millet~\cite{CM}. To this end we
need to show that the following two statements are true.
\begin{sta}\label{sta1} For each $R>0$, the set $\left\{\deld{J_h^0}y_h:\,h\in
B_R\right\}$
 is a compact subset of $\mathcal X_T$.
\end{sta}

In the above, for $R>0$ we denote by $B_R$ the closed call of radius $R$ in the Hilbert space
$L^2\la 0,T;\R^3\ra$ endowed with the weak topology.

\begin{sta}\label{sta2} Assume that $R>0$, that $(\eps_n)$  is an $(0,1]$-valued  sequence
convergent to $0$,  that $\la h_n\ra$  is  a
sequence of $\R^3$-valued $\mathbb{F}$-predictable $\mathbb{R}^3$-valued processes, indexed by $[0,T]$,  such that
$\|h\|_T\le R$ on $\Omega$ and  the laws $\mathcal{L}(h_n)$
converge weakly on $B_R$ to the law $\mathcal{L}(h)$. Then the processes
\[
\fourIdx{}{0}{}{}{C}([0,T],\mathbb{R}^3) \ni \omega \mapsto J^{\eps_n}( \omega+ \frac{1}{\sqrt{\eps_n}}\int_0^{\cdot} h_{\eps_m}(s) \,ds ) \in \mathcal{X}_T
\]
converge in law on $\mathcal X_T$ to $J^0(\int_0^{\cdot} h(s) \,ds )$.
\end{sta}
\smallskip\noindent
The remaining part of this section is devoted to the proof of these two
statements.
\begin{lemma}\label{lem19}
Suppose that $({h}_n)\subset L^2\la 0,T;\R^3\ra$ is a sequence converging weakly
to $h$. Then the sequence $y_{{h}_n}$ converges strongly to $y_h$ in $\mathcal
X_T$. In particular, for every $R>0$, the mapping
\[  B_R\ni h \mapsto J^0\bigl( \int_0^{\cdot} h(s)\, ds\bigr)\in\mathcal X_T\]
is Borel.\\ In particular, if $R>0$ and  $h$ and $\tilde{h}$ are two $B_R$-valued random
variables, possibly defined on different probability spaces, with the same laws,
then the laws of $\mathcal X_T$-valued random variables $\Omega \ni \omega
\mapsto J^0\bigl( \int_0^{\cdot} h(s,\omega)\, ds\bigr) \in \mathcal X_T$ and
$\tilde{\Omega} \ni \tilde{\omega} \mapsto J^0\bigl( \int_0^{\cdot} \tilde{h}(s,\tilde{\omega})\, ds\bigr) \in
\mathcal X_T$ are also equal.
\end{lemma}

\begin{lemma}\label{lem22}
Assume that $R>0$ and that  an $(0,1]$-valued sequence  $\la\eps_n\ra$ converges to  $ 0$
and  $\la h_n\ra$  is an $\mathcal P_T(\mathbb{R}^3)$-valued
 sequence such that
\begin{equation}\label{cond-h_n}\sup_{n\in \mathbb{N}}\int_0^T \vert h_n(t)\vert^2\, dt \le R,\quad\mathrm{for\,\,
every}\,\,\omega\in\Omega,
\end{equation}
and $\mathcal{L}(h_n)$ converges to $\mathcal{L}(h)$ weakly on $B_R$.
  Then the sequence of $\mathcal X_T$-valued random variables
\[
\fourIdx{}{0}{}{}{C}([0,T],\mathbb{R}^3) \ni \omega \mapsto J^{\eps_n}( \omega+ \frac{1}{\sqrt{\eps}}\int_0^{\cdot} h_n(s) \,ds )- J^0(\int_0^{\cdot} h_n(s) \,ds) \in \mathcal{X}_T
\]
 converges in probability to $0$.
\end{lemma}

It seems that it will be useful to introduce some temporary notation. The process (of function) $J^0(\int_0^{\cdot} h(s) \,ds)$ will be denoted by $\Phi^0(h)$ and the process $\fourIdx{}{0}{}{}{C}([0,T],\mathbb{R}^3) \ni \omega \mapsto J^{\eps}( \omega+ \frac{1}{\sqrt{\eps}}\int_0^{\cdot} h(s) \,ds )$ will be denoted by  $\Phi^\eps(h)$.

\begin{remark}\label{rem-lem22} The weak to strong continuity result is related to the weak to weak continuity results from \cite{Ball_1997,Rosa_1998,Brz+Li_2006,Daners_2005}. In the first three of these references the weak to weak continuity is an important tool in proving the existence of an attractor for 2D  Navier-Stokes Equations (stochastic with an additive noise in the latter) in unbounded domains. Recently, this notion has been generalised to the framework of stochastic PDEs (with multiplicative noise) and used in \cite{Brz+Motyl+Ondrejat_2016} in the  proof of the existence of an invariant measure for 2D  Navier-Stokes Equations in unbounded domains, where the compactness of the embedding from the Sobolev space $H^1$ to $L^2$ does not hold. We use the weak to weak continuity in a rather sketchy way at the beginning of the proof of Lemma \ref{lem22}.
It would be interesting to understand in a deeper way the relationship between these two notions.
\end{remark}

\subsection{Proof of the main result from this section}
\label{subsec-proof main}

\begin{proof}[Proof of Theorem \ref{thm-large deviations}]

Obviously Statement  \ref{sta1} follows from Lemma~\ref{lem19}.

The proof we
propose here seem to based on a new idea of using deterministic result from
Statement 1 and the Skorokhod embedding theorem on a separable metric space
$B_R$.
\par\noindent
Now we will occupy ourselves with a proof of Statement \ref{sta2}. For this aim let us  choose and fix  that $R>0$. Consider also
 an  $(0,1]$-valued  sequence $(\eps_n)$ that is convergent  to $0$ and
 a sequence $\la h_n\ra$ of $\la\mathbb{F}\ra$-predictable processes satisfying condition \eqref{cond-h_n}
 that converges to ${h}$ in law on $B_R$. Then, the following claims hold true.
 \begin{trivlist}
 \item[\textbf{(a)}]
 the $\mathcal X_T$-valued random
variables   $\Phi^{\eps_n}(h_n)-\Phi^0(h_n)$ of  converges in probability  to $0$,
  \item[\textbf{(b)}] $\Phi^0(h_n)$ converges in law on $\mathcal X_T$ to
$\Phi^0(h)$.
\end{trivlist}
Claim \textbf{(a)} follows from by  Lemma~\ref{lem22}. \\
To prove Claim \textbf{(b)} let us first recall that
 $B_R$ is a separable metric space.  By the assumptions, the laws on $B_R$ of the sequence of laws  $\la \mathcal{L}(h_n)\ra$ converges weakly to
the law  $\mathcal{L}(h)$. Hence, by the Skorodhod Theorem,  see for example,~\cite[Theorem
4.30]{Kall}, there exists a  probability space $\bigl( \tilde{\Omega},
\tilde{\mathcal {F}}, \tilde{ \mathbb{P}}\bigr)$,
  and, on that probability space,  there exist $ B_R$-valued random variables $\tilde{h}_n$ and
$\tilde{h}$, with the same laws as $h_n$ and $h$, such that
$\tilde{h}_n \to \tilde{h}$ in $B_R$,  pointwise on $\tilde \Omega$.
  By the main part of Lemma \ref{lem19} this implies that
 \[\Phi^0(\tilde{h}_n) \to \Phi^0(\tilde{h}) \mbox{ in } \mathcal{X}_T \mbox{
pointwise on  } \tilde\Omega.\]
Moreover, by the second part of Lemma \ref{lem19}, the laws of
$\Phi^0(\tilde{h}_n)$ and $ \Phi^0(\tilde{h})$  are equal, respectively, to the laws
 of  $\Phi^0(h_n)$ and $ \Phi^0(h)$.

Note that  we can choose a subsequence, without introducing a new
notation such that
\begin{trivlist}
 \item[\textbf{(a')}]
 the sequence   $\Phi^{\eps_n}(\tilde{h}_n)-\Phi^0(\tilde{h}_n)$ of $\mathcal
X_T$-valued random variables converges to $0$, $\tilde{\mathbb{P}}$- almost surely.
 \end{trivlist}

 These two convergence results imply that $\la \Phi^{\eps_n}(h_n)\ra$
converges in law on $\mathcal X_T$ to $\Phi^0(h)$. Indeed, for any globally
Lipschitz continuous and bounded function $f:\mathcal
X_T\rightarrow\mathbb{R}$, see Dudley \cite[Theorem 11.3.3]{Dudley_1989} we have
\begin{eqnarray*}
&&\hspace{-2.0cm}\lefteqn{\left| \imp_{\mathcal X_T} f(x) \,d\mathcal{L}
(\Phi^{\eps_n}(h_n)) -\imp_{\mathcal X_T}
f(x)\,d\mathcal{L}(\Phi^0(h))\right|
=\left| \imp_{\mathcal X_T} f(x) \,d\mathcal{L}( \Phi^{\eps_n}(\tilde{h}_n))
-\imp_{\mathcal X_T} f(x)\,d\mathcal{L}(\Phi^0(\tilde{h}))\right|
}\\
&=&\left|\imp_{\tilde
\Omega}f\left(\Phi^{\eps_n}(\tilde{h}_n)\right)\,d\tilde{\mathbb{P}}
-\imp_{\Omega} f(\Phi^0(\tilde{h}))\,d \tilde{\mathbb{P}}\right|
\\
&\le&
\imp_{\tilde
\Omega}|f(\Phi^{\eps_n}(\tilde{h}_n))-f(\Phi^0(\tilde{h}_n))|\,d\tilde{\mathbb{P}
} +
 \left\vert \imp_{\tilde \Omega} f(\Phi^{0}(\tilde{h}_n))\,d \mathbb{P}
-\imp_{{\tilde \Omega}} f(\Phi^0(\tilde{h}))\,d \tilde {\mathbb{P}}\right\vert
\end{eqnarray*}
Now, we observe that because $\Phi^{0}(\tilde{h}_n) \to \Phi^0(\tilde{h})$ a.s. and $f$ is
a bounded and continuous function, we infer that the 2nd term on the RHS
converges to $0$. The first term converges to $0$ because  it is bounded by
$\vert f\vert_{Lip} \imp_{\tilde \Omega}|\Phi^{\eps_n}(h_n)-\Phi^0(h_n)|\,d\tilde
{\mathbb{P}}$  and the sequence $\Phi^{\eps_n}(h_n)-\Phi^0(h_n)$ is $\tilde{\mathbb{P}}$-a.s.
convergent.

Therefore, Statement \ref{sta2} is true as well and thus we conclude the proof
of Theorem \ref{thm-large deviations}.
\end{proof}

\subsection{Proof of the auxiliary results}
\label{subsec-proof auxiliary}

\begin{proof}[Proof of Lemma \ref{lem19}] Let us assume that $h_n \to h$ weakly in $L^2(0,T;\mathbb{R}^3)$.
To simplify notation, we write $\yn$
for $y_{{h}_n}$, $y$ for $y_{h}$ and set $u_n=y_n-y$. We have to show that $u_n \to 0$ in $\mathcal{X}_T$.

Let us put
\begin{equation}\label{ineq-h}
R^2=\sup_{n\in\mathbb{N}}\imp_0^T\left|{h}_n\right|^2(s)\,ds.\end{equation}

 By Theorems\footnote{In fact, the corresponding result for equations \eqref{eqn-deterministic}, or \eqref{eqn-deterministic-2}, could be proven directly without invoking the Girsanov Theorem.}~\ref{thm-weak solution} and~\ref{thm-reg} and the uniqueness of solutions, there exists a finite
constant $C=C\la T,\alpha,R,\rho \ra$,  such that if  $\vert M_0\vert_{\mathbb{H}^1} \leq \rho$, then
\begin{equation}\label{e101}
\sup_{t\in [0,T]}|\yn(t)|_{\mathbb{H}^1}\le C,\;\;\; n\in
\overline{\mathbb{N}},
\end{equation}
\begin{equation}\label{e102}
\imp_0^T\la\left|\Delta \yn(s)\right|^2_{\mathbb{L}^2}+\left|\D
y_n\right|^4_{\LL^4}\ra\,ds\le C,\;\;\; n\in
\overline{\mathbb{N}}
\end{equation}
and
\begin{equation}\label{eqn-103}
\imp_0^T \left| \yn^\prime(s)\right|^2_{\mathbb{L}^2}\,ds\le C,\;\;\; n\in
\overline{\mathbb{N}}.
\end{equation}

Let us also recall that
\begin{equation}\label{e102a}
|\yn(t)(x)|=1,\quad x\in\domain,\,\,t\in [0,T],\;\; n\in
\overline{\mathbb{N}}.
\end{equation}
The same properties  hold for $y$. Hence, in particular,
\begin{equation}\label{e102b}
|u_n(t)(x)|\leq 2,\quad x\in\domain,\,\,t\in [0,T],\;\; n\in
\overline{\mathbb{N}}.
\end{equation}

\par\medskip\noindent
\emph{Step 1.}
The above estimates, together with standard compactness results, see e.g. \cite{Temam_2001}, imply that the sequence $y_n$ has a subsequence, for which we do not introduce a separate notation, which converges weakly in $L^2(0,T;\mathbb{H}^2)$, strongly in $L^q(0,Y;\mathbb{H}^1)\cap  C([0,T];\mathbb{L}^2) $, for any $q<\infty$ and in $C_{\rm w}([0,T];\mathbb{H}^1)$ to some $\bar{y} \in  L^2(0,T;\mathbb{H}^2) \cap  C([0,T];\mathbb{H}^1)$ such that $\bar{y}^\prime \in L^2(0,T;\LL^2)$. Standard argument, see e.g. section 7 of \cite{Brz+Li_2006}  imply that $\bar{y}$ is a  unique solution of the problem \eqref{eqn-deterministic-2}. A deterministic version of our uniqueness result Theorem \ref{thm-uniqueness pathwise weak} implies then that, recall that $y=y_h$,  $\bar{y}=y$.
Using the subsequence argument, we deduce that  the whole sequence $y_n$ converges to $y$  weakly in $L^2(0,T;\mathbb{H}^2)$, strongly in $L^q(0,Y;\mathbb{H}^1)\cap  C([0,T];\mathbb{L}^2)$, for any $q<\infty$ and in $C_{\rm w}([0,T];\mathbb{H}^1)$.

\par\medskip\noindent
\emph{Step 2.}
Let $q\in L^2(0,T;\LL^2)$.  We claim that
\begin{equation}\label{lem21}
\lim_{n\to\infty}\la\sup_{t\in [0,T]}\left|\imp_0^t\langle q(s),
u_n(s)\rangle_{\mathbb{L}^2}({h}_n(s)-{h}(s))\,ds\right|\ra=0.
\end{equation}
By \emph{Step 1} we can assume that there exists an element $u_\infty \in C([0,T];\LL^2)$ such that  $u_n\to u_\infty$ in $C([0,T];\LL^2)$. For $n
\in \mathbb{N}\cup \{\infty\}$ we define an operator $\mathcal K_n:L^2\la
0,T;\R^3\ra\to \fourIdx{}{0}{}{}{C\la[0,T];\R^3\ra}$ by the following formula
\[\mathcal K_nv(t)=\imp_0^t\left\langle q(s),u_n(s)\right\rangle_{\mathbb{L}^2}
v(s)\,ds,\quad t\in [0,T],\;\; v\in L^2\la 0,T;\R^3\ra.\]
Each operator $\mathcal K_n$,  is compact because the
function $\langle q(\cdot), u_n(\cdot)\rangle_{\mathbb{L}^2}$ belongs to $L^2\la
0,T;\R\ra$. Moreover,  since the sequence  $\langle q(\cdot),
u_n(\cdot)\rangle_{\mathbb{L}^2}$  converges strongly in $L^2\la 0,T;\R\ra$ to a
function $\langle q(\cdot), u_\infty(\cdot)\rangle_{\mathbb{L}^2}$ we infer that
\[\lim_{n\to\infty}\left\|\mathcal K_n-\mathcal K_\infty\right\|=0.\]

 Since \deli{${h}_n\to {h}$ weakly in $L^2\la [0,T];\R^3\ra$. Then}
\[\left|\mathcal K_n\la
{h}_n-{h}\ra\right|_{\fourIdx{}{0}{}{}{C\la[0,T];\R^3\ra})}\le\left\|\mathcal
K_n-\mathcal
K_\infty\right\|\cdot\left|{h}_n-{h}\right|_{L^2\la[0,T];\R^3\ra}+\left|\mathcal
K_\infty\la {h}_n-{h}\ra\right|_{\fourIdx{}{0}{}{}{C\la[0,T];\R^3\ra}},\]
 the claim \eqref{lem21} follows immediately by the compactness of $K_\infty$ because ${h}_n\to
{h}$ weakly in $L^2\la [0,T];\R^3\ra$.
\par\medskip\noindent
\emph{Step 3.} We will show that

 \begin{equation}\label{e302a}
\lim_{n \to \infty }\;\Bigl[ \;\sup_{t\in
[0,T]}|u_n(t)|^2_{\mathbb{H}^1}+\alpha\imp_0^T|\Delta
u_n|^2_{\mathbb{L}^2}\,ds\;\Bigr]=0.
\end{equation}

Without loss of generality we may assume that $\newe_2=\newe_3=0$ and we put
$\newe=\newe_1 \in \mathbb{H}^1$.  In particular, we can assume that all functions $h$ and  $h_n$ are $\mathbb{R}$-valued, i.e. $h,h_n \in L^2(0,T;\mathbb{R})$. Note that in this case, the last term in \eqref{eqn-deterministic-2} reads
\begin{equation}\label{eqn-simplification}
\imp_0^tG\la y_h\ra Bh(s)\,ds=\imp_0^t  h(s) G( y_h(s)) e\,ds.
\end{equation}

Let us recall, see \eqref{eqn-pi_n}, that for $N\in \mathbb{N}$,
$\pi_N:\mathbb{L}^2 \to \mathbb{H}_n$ is the orthogonal projection onto the finite dimensional subspace of $\LL^2$ spanned by the first $N$ eigenvectors of the Neumann Laplacian.

For the aim of proving \eqref{e302a} we will show that  there exist $C>0$ such that for every $N\in\mathbb{N}$,
\begin{equation}\label{e302b}
\limsup_{n\to\infty}\la\sup_{t\in [0,T]}|\D u_n(t)|^2_{\LL^2}+\alpha\imp_0^T|\Delta
u_n|^2_{\mathbb{L}^2}\,ds\ra\le
C\left|e-\pi_Ne\right|_{\mathbb{L}^2}^{2}.
\end{equation}

By \eqref{eqn-deterministic-2}, for each $n\in \mathbb{N}$, we have

\begin{eqnarray}\label{eqn-deterministic-3}%{e108}
u_n(t)&=&
\alpha\imp_0^t\Delta u_n\,ds
+\alpha\imp_0^t \Bigl\{ \left|\D y_{n}\right|^2y_{n}-\left|\D y\right|^2y\Bigr\} \, ds\\
&+&\imp_0^t \Bigl\{ y_{n}\times\Delta y_{n}-  y\times\Delta y\Bigr\}\,ds
-\beta \imp_0^t \Bigl\{ G\la y_{n}\ra F\la y_{n}\ra-G\la y\ra F\la y\ra \Bigr\}\,ds\nonumber\\
&+&\imp_0^t\Bigl\{G\la y_{n}\ra Bh_n(s)  -G\la y\ra Bh(s)\Bigr\} \,ds,
\nonumber
\end{eqnarray}
Therefore, by some simple algebraic manipulations, formula \eqref{eqn-simplification}  and the linearity\footnote{In fact, Lipschitz property of $f$ would be sufficient.} of the function $f$, we infer that
\begin{eqnarray}
\label{eqn-deterministic-4}
u_n(t)&=&
\alpha\imp_0^t\Delta u_n\,ds\\
&+&\alpha\imp_0^t(|\D y_n|-|\D \yw|)(|\D y_n|+|\D
\yw|)y_n\,ds+\alpha\imp_0^t|\D \yw|^2u_n\,ds
\nonumber\\
&+& \imp_0^tu_n\times\Delta \yn\,ds+\imp_0^t\yw\times\Delta u_n\,ds
\nonumber\\
&-&\beta \imp_0^t G\la y_n(s)\ra F\la u_n(s)\ra\, ds-\beta\imp_0^t\la
G\la y_n(s)\ra-G(y(s))\ra f(y(s))\, ds\nonumber\\
&+&\imp_0^t h_n(s) \bigl[ G\la y_n(s)\ra e - G\la y(s)\ra e \bigr] \, ds+\imp_0^t [{h}_n(s)-{h}(s)]\,G(y(s))e\,ds.\nonumber
\end{eqnarray}

In order to prove \eqref{e302a} we could follow a standard method of getting a'priori bounds by invoking the Gronwall Lemma. This would work easily but the last term  on the RHS of \eqref{eqn-deterministic-4}. In order to be able to deal with that term we could use Step 2. However this would work had the  function $e$ were more regular, at least from $\mathbb{H}^2$. For this purpose, we will introduce an approximation of $e$ by a sequence of more regular functions, see below, and then prove, instead of \eqref{e302a}, \eqref{e302b}.

Since $u_n$ is a strong solution of the above equation, by Lions-Magenes \cite{Lions+Magenes_1972}, we infer that

\begin{eqnarray}
\label{eqn-deterministic-5}
\frac12\frac{d}{dt}\vert \D u_n(t)\vert^2&=&-\alpha \vert \Delta  u_n(t)\vert^2
-\alpha \lb (|\D y_n(t)|-|\D \yw(t)|)(|\D y_n(t)|+|\D
\yw(t)|)y_n(t), \Delta u_n(t)\rb
 \nonumber\\
 &-&\alpha \lb |\D \yw(t)|^2u_n(t), \Delta u_n(t)\rb - \lb u_n(t)\times\Delta \yn(t), \Delta u_n(t)\rb
\\
 &+&\beta \lb G\la y_n(t)\ra F\la u_n(t)\ra, \Delta u_n\rb +\beta \lb \la
G\la y_n(t)\ra-G(y(t))\ra f(y(t)), \Delta u_n\rb \nonumber\\
&-&\lb h_n(t) \bigl[ G\la y_n(t)\ra e - G\la y(t)\ra e \bigr] , \Delta u_n\rb
-\lb  [{h}_n(t)-{h}(t)]\,G(y(t))e, \Delta u_n\rb
\nonumber
\end{eqnarray}

Let us now fix an auxiliary natural number $N$. Subtracting and adding $\pi_Ne$ in the last term of the above equality and using integration by parts   we get  \begin{eqnarray}
\label{eqn-deterministic-5}
\frac12\frac{d}{dt}\vert \D u_n(t)\vert^2&=&-\alpha \vert \Delta  u_n(t)\vert^2
-\alpha \lb (|\D y_n(t)|-|\D \yw(t)|)(|\D y_n(t)|+|\D
\yw(t)|)y_n(t), \Delta u_n(t)\rb
 \nonumber\\
 &-&\alpha \lb   |\D \yw(t)|^2u_n(t), \Delta u_n(t)\rb
 - \lb  u_n(t)\times\Delta \yn(t), \Delta u_n(t)\rb
\\
 &+&\beta \lb G\la y_n(t)\ra F\la u_n(t)\ra, \Delta u_n(t)\rb +\beta \lb \la
G\la y_n(t)\ra-G(y(t))\ra f(y(t)), \Delta u_n(t)\rb \nonumber\\
&-&\lb h_n(t) \bigl[ G\la y_n(t)\ra e - G\la y(t)\ra e \bigr] , \Delta u_n\rb
\nonumber\\
&-&\lb  [{h}_n(t)-{h}(t)] \bigl( G(y(t))e-G(y(t)) \pi_N e\bigr), \Delta u_n\rb
\nonumber\\
&-&\lb  [{h}_n(t)-{h}(t)]\,\Delta\bigl( G(y(t)) \pi_N e\bigr), u_n\rb
\nonumber
\end{eqnarray}
Let us show how we estimate each of the terms on the RHS above. All norm below, unless otherwise stated, are those in $\LL^2$. We fix $\eps>0$. We begin with the 1st term. By the Young inequality and inequality \eqref{e102a} we have, where $C$ is the constant from the GNI below,
\begin{eqnarray*}
- \alpha \lb (|\D y_n|-|\D \yw|)|\D
\yw|y_n, \Delta u_n\rb   &\leq&
 \frac{\eps}2  \vert \Delta u_n \vert^2+\frac1{2\eps} \vert \D u_n\vert_{\LL^4}^2   \vert \D
\yw\vert_{\LL^4}^2
\end{eqnarray*}
Note that by the Gagliardo-Nirenberg inequality (and again and inequality \eqref{e102a})
\begin{eqnarray*}
\vert \D u_n\vert_{\LL^4}^2 &\leq& C \bigl[ \vert u_n\vert^2+ \vert \Delta u_n\vert^2\bigr]^{1/2} \vert u_n\vert_{\LL^\infty} \leq C \bigl[ \vert u_n\vert+ \vert \Delta u_n\vert\bigr] \vert u_n\vert_{\mathbb{H}^1}
\end{eqnarray*}
Hence

\begin{eqnarray*}
\vert \D u_n\vert_{\LL^4}^2   \vert \D
\yw\vert_{\LL^4}^2 &\leq& \frac{\eps^2}{2}\vert \Delta u_n \vert^2+ (C+\frac{1}{2\eps^2})\vert \D
\yw\vert_{\LL^4}^4\vert u_n\vert_{\mathbb{H}^1}^2
\end{eqnarray*}

Therefore,
\begin{eqnarray}\label{ineq-201}
- \alpha \lb (|\D y_n|-|\D \yw|)|\D
\yw|y_n, \Delta u_n\rb   &\leq&
\frac{3\eps}{4}  \vert \Delta u_n \vert^2+ (\frac{C}{2\eps}+\frac{1}{4\eps^3})\vert \D
\yw\vert_{\LL^4}^4\vert u_n\vert_{\mathbb{H}^1}^2.
\end{eqnarray}
Similarly, we also get
\begin{eqnarray}\label{ineq-202}
- \alpha \lb (|\D y_n|-|\D \yw|)|\D
\yw_n|y_n, \Delta u_n\rb   &\leq&
\frac{3\eps}{4}  \vert \Delta u_n \vert^2+ \frac{1}{4\eps^3})\vert \D
\yw_n\vert_{\LL^4}^4\vert u_n\vert_{\mathbb{H}^1}^2.
\end{eqnarray}
In an almost identical way, where instead of inequality \eqref{e102a} we use \eqref{e102b}, we get
\begin{eqnarray}\label{ineq-203}
-\alpha \lb   |\D \yw|^2u_n, \Delta u_n\rb
&\leq&
\frac{\eps}{2}  \vert \Delta u_n \vert^2+  \frac{2\alpha^2}{\eps} \vert \D
\yw\vert_{\LL^4}^4\vert u_n\vert^2.
\end{eqnarray}
Next, with $C>0$ such that $\vert u \vert_{\LL^\infty}^2 \leq \vert u \vert_{\mathbb{H}^1}^2$ we have

\begin{eqnarray}
\nonumber
- \lb  u_n\times\Delta \yn, \Delta u_n\rb &\leq&
\frac{\eps}{2}  \vert \Delta u_n \vert^2+  \frac{1}{2\eps} \vert u_n\vert_{\LL^\infty}^2 \vert \Delta \yn \vert^2
\\
&\leq&
\frac{\eps}{2}  \vert \Delta u_n \vert^2+  \frac{C}{2\eps}  \vert \Delta \yn \vert^2 \vert u_n\vert_{\mathbb{H}^1}^2
\label{ineq-204}
\end{eqnarray}
The next two terms are easy. By inequalities \eqref{e102a} and  \eqref{e102b},  and  the Lipschitz continuity of functions  $f$ and $g$  on balls we infer that
\begin{eqnarray}
\label{ineq-205}
\beta \lb G\la y_n(t)\ra F\la u_n(t)\ra, \Delta u_n\rb
&\leq&  \frac{\eps}{2}  \vert \Delta u_n \vert^2+  \frac{C\beta^2}{2\eps} \vert u_n\vert^2
\end{eqnarray}
and
\begin{eqnarray}
\label{ineq-206}
\beta \lb \la
G\la y_n(t)\ra-G(y(t))\ra f(y(t)), \Delta u_n\rb  &\leq&  \frac{\eps}{2}  \vert \Delta u_n \vert^2+  \frac{C\beta^2}{2\eps} \vert u_n\vert^2
\end{eqnarray}
Next, by   inequality \eqref{e102a}   and  the Lipschitz continuity of function $g$  on balls we infer that
\begin{eqnarray}
\label{ineq-207}
\lb h_n(t) \bigl[ G\la y_n(t)\ra e - G\la y(t)\ra e \bigr] , \Delta u_n(t)\rb
&\leq & \frac{\eps}{2}  \vert \Delta u_n \vert^2+  \frac{C}{2\eps} \vert h_n(t)\vert^2 \vert u_n(t)\vert^2
\end{eqnarray}

\begin{equation}
\label{ineq-208}
\lb  [{h}_n(t)-{h}(t)] \bigl( G(y(t))e-G(y(t))\pi_N e\bigr), \Delta u_n\rb
\leq  \frac{\eps}{2}  \vert \Delta u_n \vert^2+  \frac{C}{2\eps} \vert {h}_n(t)-{h}(t)\vert^2 \vert e -\pi_N e\vert^2
\end{equation}
We leave  the last term unchanged. From all the inequalities above we deduce that
\begin{eqnarray*}
\label{eqn-deterministic-6}
\frac12\frac{d}{dt}\vert \D u_n(t)\vert^2&+&\alpha \vert \Delta  u_n(t)\vert^2
\\
&\leq&
\frac{3\eps}{4}  \vert \Delta u_n(t) \vert^2+(\frac{C}{2\eps} +\frac{1}{4\eps^3})\vert \D
\yw(t)\vert_{\LL^4}^4\vert u_n(t)\vert_{\mathbb{H}^1}^2
\\
&+&
\frac{3\eps}{4}  \vert \Delta u_n(t) \vert^2+ (\frac{C}{2\eps}+\frac{1}{4\eps^3})\vert \D
\yw_n(t)\vert_{\LL^4}^4\vert u_n(t)\vert_{\mathbb{H}^1}^2
\\
&+&
\frac{\eps}{2}  \vert \Delta u_n(t) \vert^2+  \frac{2\alpha^2}{\eps} \vert \D
\yw(t)\vert_{\LL^4}^4\vert u_n(t)\vert^2
\\
&+& \frac{\eps}{2}  \vert \Delta u_n(t) \vert^2+  \frac{C}{2\eps}  \vert \Delta \yn \vert^2 \vert u_n(t)\vert_{\mathbb{H}^1}^2
\\
&+&   \frac{\eps}{2}  \vert \Delta u_n(t) \vert^2+  \frac{C\beta^2}{2\eps} \vert u_n(t)\vert^2
\\
&+&  \frac{\eps}{2}  \vert \Delta u_n(t) \vert^2+  \frac{C\beta^2}{2\eps} \vert u_n(t)\vert^2
\\
&+& \frac{\eps}{2}  \vert \Delta u_n(t) \vert^2+  \frac{C}{2\eps} \vert h_n(t)\vert^2 \vert u_n(t)\vert^2
\\
&+& \frac{\eps}{2}  \vert \Delta u_n(t) \vert^2+  \frac{C}{2\eps} \vert {h}_n(t)-{h}(t)\vert^2 \vert e -\pi_N e\vert^2
\\
&-& \lb  [{h}_n(t)-{h}(t)]\,\Delta\bigl( G(y(t)) \pi_N e\bigr), u_n(t)\rb
\end{eqnarray*}
Let us now choose $\eps=\frac{\alpha}{9}>0$, i.e.  such that
\[
\frac\alpha2= \bigl(2 \times \frac34+6 \times \frac12 \bigl)\eps=\frac92\eps.
\]
Then we get
\begin{eqnarray}
\label{eqn-deterministic-6}
\frac{d}{dt}\vert \D u_n(t)\vert^2&+&\alpha \vert \Delta  u_n(t)\vert^2
\\
&\leq& \psi_n(t) \vert u_n(t)\vert_{\mathbb{H}^1}^2+ \chi_n(t) \vert u_n(t)\vert_{\LL^2}^2
+  \frac{9C}{2\alpha} \vert {h}_n(t)-{h}(t)\vert^2 \vert e -\pi_N e\vert^2
\\
&-& \lb  [{h}_n(t)-{h}(t)]\,\Delta\bigl( G(y(t)) \pi_N e\bigr), u_n(t)\rb
\nonumber
\end{eqnarray}
where
\begin{eqnarray}
\label{eqn-deterministic-7}
\psi_n(t)&=& (\frac{9C}{2\alpha} +\frac{9^3}{4\alpha^3})\bigl( \vert \D
\yw(t)\vert_{\LL^4}^4+ \vert \D
\yw_n(t)\vert_{\LL^4}^4\bigr) +\frac{9C}{2\alpha}  \vert \Delta \yn \vert^2
\\
\chi_n(t)&=&  \frac{18\alpha^2}{\alpha} \vert \D
\yw(t)\vert_{\LL^4}^4+\frac{9C\beta^2}{\alpha}+  \frac{9C}{2\alpha} \vert h_n(t)\vert^2
\end{eqnarray}
Therefore, with
\[
b_{n,N}:=\sup_{t \in [0,T]} \vert \lb  [{h}_n(t)-{h}(t)]\,\Delta\bigl( G(y(t)) \pi_N e\bigr), u_n(t)\rb\vert
\]
we infer that for any $t\in [0,T]$,
\begin{eqnarray}
\label{eqn-deterministic-8}
\vert \D u_n(t)\vert^2&+&\alpha \int_0^t \vert \Delta  u_n(s)\vert^2\, ds
\\
&\leq& \int_0^t \psi_n(s) \vert u_n(s)\vert_{\mathbb{H}^1}^2\,ds + \int_0^T \chi_n(s) \vert u_n(s)\vert_{\LL^2}^2 \,ds\\
&+&  \frac{9C}{\alpha} \vert e -\pi_N e\vert^2 \int_0^T \bigl( \vert {h}_n(s)\vert^2 +\vert {h}(s)\vert^2 \bigr) \, ds
\nonumber
+ b_{n,N}.
\nonumber
\end{eqnarray}
Therefore, by the Gronwall Lemma and our assumption \eqref{ineq-h} we get
\[
\vert \D u_n(t)\vert^2 \leq \Bigl[  \frac{18CR^2}{\alpha} \vert e -\pi_N e\vert^2+b_{n,N}
+ \int_0^T \chi_n(s) \vert u_n(s)\vert_{\LL^2}^2 \,ds\Bigr]\; e^{\int_0^t \psi_n(s)}, \;\; t \in [0,T].
\]

By estimates (\ref{e101}) and (\ref{e102})
\[\gamma:=\sup_{n\in \mathbb{N}}\imp_0^T\psi_n(s)\,ds<\infty\]
and $\gamma$ depends on $\alpha$, $T$, $R$, $\rho$ and $r$ only. Therefore,  we infer that there exists a constant $C_T>0$ such that
\begin{equation}\label{e302}
\sup_{t\in [0,T]}|u_n(t)|^2_{\mathbb{H}^1}+\alpha\imp_0^T|\Delta
u_n|^2_{\mathbb{L}^2}\,ds\leq C_Te^{\gamma T}\Bigl[  \frac{18CR^2}{\alpha} \vert e -\pi_N e\vert^2+b_{n,N}
+ \int_0^T \chi_n(s) \vert u_n(s)\vert_{\LL^2}^2 \,ds\Bigr].
\end{equation}
Therefore, since by Claim\eqref{lem21} $b_{n,N} \to0$ as $n\to \infty$, and, by Step 1, $\int_0^T \chi_n(s) \vert u_n(s)\vert_{\LL^2}^2 \,ds $ converges to 0, we conclude the proof of
\eqref{e302b} and so of \eqref{e302a} as well.

\par\medskip\noindent
\emph{Step 4.}
We complete the proof of Lemma \ref{lem19} by taking the limit as $N\to\infty$.
\end{proof}
Note, that Statement \ref{sta1} follows Lemma \ref{lem19}.

\medskip

Now we will occupy ourselves with the proof of that Statement \ref{sta2}. For
this purpose let us chhose and fix the following processes:
\[Y_n=\Phi^{\eps_n}({h}_n) \mbox{ and }y_n=\Phi^0({h}_n).\] \deli{and let}
Let $N>\left|M_0\right|_{\mathbb{H}^1}$ be fixed. For each $n\in \mathbb{N}$ we define
an $\la\mathcal F_t\ra$-stopping time
\begin{equation}\label{e114}
\tau_n=\inf\left\{t>0: |Y_n(t)|_{\mathbb{H}^1}\geq N\right\}\wedge T.
\end{equation}
\begin{lemma}\label{lem23}
For $\tau_n$ as defined in (\ref{e114}) we have
\[
\lim_{n\to\infty}\mathbb{E}\left(\sup_{t\in [0,T]}\left|Y_n\la
t\wedge\tau_n\ra-y_n\la
t\wedge\tau_n\ra\right|^2_{\mathbb{L}^2}+\imp_0^{\tau_N}\left|Y_n-y_n\right|^2_{
\mathbb{H}^1}\,ds\right)=0.
\]
\end{lemma}
\noindent\begin{proof} Let $X_n=Y_n-y_n$. We assume without loss of generality
that $\beta=0$, $\newe_2=\newe_3=0$ and $\newe_1=h$. Then for any $n\in
\mathbb{N}$ we have
\begin{equation}\label{yu}
\begin{aligned}
dX_n&=\alpha\Delta X_n dt\\
&\,\,\,\,\,\,+\alpha\la\D X_n\ra\cdot\la\D\la Y_n+y_n\ra\ra Y_n
dt+\alpha\left|\D y_n\right|^2X_n dt\\
&\,\,\,\,\,\,+X_n\times\Delta Y_n dt+y_n\times\Delta X_n dt\\
&\,\,\,\,\,\,+\la G\la Y_n\ra-G\la y_n\ra\ra h{h}_n dt\\
&\,\,\,\,\,\,+\sqrt{\eps_n}G\la Y_n\ra hdW+\frac{\eps_n}{2}G^\prime\la
Y_n\ra G\la Y_n\ra h dt
\end{aligned}
\end{equation}
Using a version of the It\^{o} formula given in \cite{Pardoux_1979} and
integration by parts we obtain
\[\begin{aligned}
\frac{1}{2}d\left|X_n\right|^2_{\mathbb{L}^2}&=-\alpha\left|X_n\right|^2_{\mathbb{H}^1}
dt
+\alpha\left|\left|\D y_n\right|X_n\right|^2_{\mathbb{L}^2}dt\\
&\,\,\,\,\,\,\,\,\,+\alpha\las X_n,\la\D X_n\ra\cdot\la\D\la
Y_n+y_n\ra\ra  Y_n\ras_{\mathbb{L}^2}dt\\
&\,\,\,\,\,\,\,\,\,\,\,-\las\D X_n,X_n\times\D
y_n\ras_{\mathbb{L}^2}dt\\
&\,\,\,\,\,\,\,\,\,\,\,+\las\la G\la Y_n\ra-G\la y_n\ra\ra
h,X_n\ras_{\mathbb{L}^2}{h}_ndt\\
&\,\,\,\,\,\,\,\,\,\,\,+\frac{\eps_n}{2}z_ndt+\sqrt{\eps_n}\las G\la
Y_n\ra h,X_n\ras_{\mathbb{L}^2} dW
\end{aligned}\]
where $z_n$ is a process defined by
\[z_n=\las G^\prime\la Y_n\ra G\la Y_n\ra h, X_n\ras_{\mathbb{L}^2}+\left|G\la
Y_n\ra h\right|^2_{\mathbb{L}^2}.\]
Therefore
\[\begin{aligned}
\left|X_n(t)\right|^2_{\mathbb{L}^2}+2\alpha\imp_0^t\left|X_n\right|^2_{\mathbb{H}^1}
ds&\le
C\imp_0^t\left|X_n\right|_{\mathbb{L}^2}\left|X_n\right|_{\mathbb{H}^1}
\left|y_n\right|^2_{\mathbb{H}^1}ds\\
&\,\,\,\,\,\,+C\imp_0^t\left|X_n\right|_{\mathbb{H}^1}^{3/2}\left|X_n\right|_{\mathbb{L
}^2}\la\left|y_n\right|_{\mathbb{H}^1}+\left|y_n\right|_{\mathbb{H}^1}\ra ds\\
&\,\,\,\,\,\,+C\imp_0^t\left|X_n\right|_{\mathbb{H}^1}^{3/2}\left|X_n\right|_{\mathbb{L
}^2}\left|y_n\right|_{\mathbb{H}^1} ds\\
&\,\,\,\,\,\,+C\imp_0^t\left|X_n\right|_{\mathbb{L}^2}^2\left|{h}_n\right|ds\\
&\,\,\,\,\,\,+C\eps_n+\sqrt{\eps_n}\left|\imp_0^t\las G\la Y_n\ra
h,X_n\ras_{\mathbb{L}^2} dW\right|.
\end{aligned}\]
By \eqref{e101} we have $\sup_n\left|y_n\right|_{\mathbb{H}^1}<\infty$ and therefore,
using repeatedly the Young inequality we find that there exists $C>0$ such that
for all $t\in [0,T]$
\[\begin{aligned}
\left|X_n(t)\right|^2_{\mathbb{L}^2}+\alpha\imp_0^t\left|X_n\right|^2_{\mathbb{H}^1}
ds&\le C\imp_0^t\left|X_n\right|_{\mathbb{L}^2}^2\la
1+\left|{h}_n\right|+\beta\left|y_n\right|^4_{\mathbb{H}^1}\ra ds\\
&\,\,\,\,\,\,+C\eps_n+\sqrt{\eps_n}\left|\imp_0^t\las G\la Y_n\ra
h,X_n\ras_{\mathbb{L}^2} dW\right|.
\end{aligned}\]
Denoting the left hand side of the above inequality by $L_t$ and using the
definition of $\tau_n$ we have
\[\begin{aligned}
L_{t\wedge\tau_n}&\le
C\imp_0^{t\wedge\tau_n}\left|X_n\right|_{\mathbb{L}^2}^2\la
1+\left|{h}_n\right|+\beta\left|y_n\right|^4_{\mathbb{H}^1}\ra ds\\
&\,\,\,\,\,\,+C\eps_n+\sqrt{\eps_n}\left|\imp_0^{t\wedge\tau_n}\las G\la
Y_n\ra h,X_n\ras_{\mathbb{L}^2} dW\right|\\
&\le
\imp_0^{t\wedge\tau_n}\left|X_n\right|^2_{\mathbb{L}^2}\psi_{n,N}
ds+C\eps_n\\
&\,\,\,\,\,\,+\sqrt{\eps_n}\left|\imp_0^{t\wedge\tau_n}\las G\la Y_n\ra
h,X_n\ras_{\mathbb{L}^2} dW\right|,
\end{aligned}\]
where
\[\psi_{n,N}(s)=1+\left|{h}_n(s)\right|+\beta N^4,\quad s\le T.\]
Since
\[\sup_n\sup_{t\in [0,T]}\las G\la Y_n(t)\ra,X_n(t)\ras^2\le C,\quad\P-a.s.,\]
the Burkoholder-Davis-Gundy inequality yields
\begin{equation}\label{sup}
\E\sup_{s\le t}L_{s\wedge\tau_n}\le C\sqrt{\eps_n}+\imp_0^t\E\sup_{r\le
s}\left|X_n\la r\wedge\tau_n\ra\right|^2_{\mathbb{L}^2}\psi_{n,N}ds,
\end{equation}
and therefore
\[\mathbb{E}\sup_{r\le t}|X_n(r\wedge\tau_n)|^2_{\mathbb{L}^2}\le
C\sqrt{\eps_n}+\imp_0^t\E\sup_{r\le s}\left|X_n\la
r\wedge\tau_n\ra\right|^2_{\mathbb{L}^2}\psi_{n,N}ds.
\]
Clearly,
\[\sup_{n\in \mathbb{N}}\imp_0^T\psi_{n,N}ds<\infty,\]
hence the Gronwall Lemma implies
\[
\mathbb{E}\sup_{r \in [0,T]}|X_n(r\wedge\tau_n)|^2_{\mathbb{L}^2}\le
C\sqrt{\eps_n}\mathrm{e}^{\imp_0^T\psi_{n,N}ds}\rightarrow 0\mbox{ as
}n\rightarrow\infty.
\]
Returning now to \eqref{sup}, we also have
\[
\mathbb{E}\imp_0^{\tau_n}|X_n(s)|^2_{\mathbb{H}^1}\,ds\le
C\sqrt{\eps_n}\mathrm{e}^{\imp_0^T\psi_{n,N}ds}\rightarrow 0\mbox{ as
}n\rightarrow\infty.
\]
This completes the proof of Lemma~\ref{lem23}.
\end{proof}
\begin{lemma}\label{lem24}
For the stopping time $\tau_n$ defined in (\ref{e114}) we have
\[
\lim_{n\to\infty}\mathbb{E}\left(\sup_{t\in
[0,T]}|\D(Y_n(t\wedge\tau_n)-y_n(t\wedge\tau_n))|^2_{\mathbb{L}^2}+\imp_0^{
\tau_n}|\Delta(Y_n-y_n)|^2_{\mathbb{L}^2}\,ds\right)=0.
\]
\end{lemma}
\begin{proof} By a version of the It\^{o} formula, see \cite{Pardoux_1979},
\[\frac{1}{2}d|\D(Y_n(t)-y_n(t))|^2_{\mathbb{L}^2}=-\las\Delta\la
Y_n-y_n\ra,d\la Y_n-y_n\ra\ras_{\mathbb{L}^2}
+\eps_n\left|\D G\la Y_n\ra h\right|^2_{\mathbb{L}^2}dt.\]
Therefore, putting $X_n=Y_n-y_n$ and invoking equality \eqref{yu} we obtain for
any $\eta>0$
\begin{equation}\label{d1}
\begin{aligned}
\frac{1}{2}d\left|\D X_n(t)\right|^2_{\mathbb{L}^2}&=-\alpha\left|\Delta
X_n\right|^2_{\mathbb{L}^2}\\
&\,\,\,\,\,\,\,-\alpha\las \Delta X_n,\D X_n\cdot\la\D Y_n+\D y_n\ra
Y_n\ras_{\mathbb{L}^2}dt\\
&\,\,\,\,\,\,\,-\alpha\las \Delta X_n,\left|\D
y_n\right|^2X_n\ras_{\mathbb{L}^2}dt\\
&\,\,\,\,\,\,\,-\las X_n\times\Delta y_n,\Delta X_n\ras_{\mathbb{L}^2}dt\\
&\,\,\,\,\,\,\,-\las \la G\la Y_n\ra-G\la y_n\ra\ra h,\Delta
X_n\ras_{\mathbb{L}^2} {h}_ndt\\
&\,\,\,\,\,\,\,-\sqrt{\eps_n}\las \D G\la Y_n\ra h,\Delta
X_n\ras_{\mathbb{L}^2}dW\\
&\,\,\,\,\,\,\,-\frac{\eps_n}{2}\las G^\prime\la Y_n\ra G\la Y_n\ra h,\Delta
X_n\ras_{\mathbb{L}^2}dt\\
&\,\,\,\,\,\,\,+\eps_n\left|\D G\la Y_n\ra h\right|^2_{\mathbb{L}^2}dt.
\end{aligned}
\end{equation}
We will estimate the terms in \eqref{d1}. First, noting that
\[\las X_n\times\Delta Y_n,\Delta X_n\ras_{\mathbb{L}^2}=\las X_n\times\Delta
u_n,\Delta X_n\ras_{\mathbb{L}^2}\]
we find that
\begin{equation}\label{da}
\left|\las X_n\times\Delta Y_n,\Delta X_n\ras_{\mathbb{L}^2}\right|\le
C\eta^2\left|\Delta
X_n\right|^2_{\mathbb{L}^2}+\frac{C}{\eta^2}\left|X_n\right|_{\mathbb{L}^2}
\left|X_n\right|_{\mathbb{H}^1}.
\end{equation}
Next, by the Young inequality and the interpolation inequality
\eqref{ineq-interp}
\[\begin{aligned}
\left|\las \Delta X_n,\D X_n\cdot\la\D Y_n+\D y_n\ra
Y_n\ras_{\mathbb{L}^2}\right|&\le C\eta^2\left|\Delta
X_n\right|^2_{\mathbb{L}^2}\\
&\,\,\,\,\,\,\,+\frac{C}{\eta^2}\imp_\domain\left|\D
X_n\right|^2\la\left|\D Y_n\right|^2+\left|\D y_n\right|^2\ra dx\\
&\le C\eta^2\left|\Delta X_n\right|^2_{\mathbb{L}^2}\\
&\,\,\,\,\,\,\,+\frac{C}{\eta^2}\left|\D
X_n\right|_\infty^2\imp_\domain\la\left|\D Y_n\right|^2+\left|\D
y_n\right|^2\ra dx\\
&\le C\eta^2\left|\Delta X_n\right|^2_{\mathbb{L}^2}\\
&\,\,\,\,\,\,\,+\frac{C}{\eta^2}\left|X_n\right|_{\mathbb{H}^1}\la
\left|X_n\right|_{\mathbb{H}^1}+\left|\Delta
X_n\right|_{\mathbb{L}^2}\ra\la\left|Y_n\right|_{\mathbb
H^1}^2+\left|y_n\right|_{\mathbb H^1}^2\ra
\end{aligned},\]
and thereby
\begin{equation}\label{db}
\begin{aligned}
\left|\las \Delta X_n,\D X_n\cdot\la\D Y_n+\D y_n\ra
Y_n\ras_{\mathbb{L}^2}\right|&\le C\eta^2\left|\Delta
X_n\right|^2_{\mathbb{L}^2}\\
&\,\,\,\,\,\,\,+\frac{C}{\eta^2}\left|X_n\right|^2_{\mathbb{H}^1}\la\left|Y_n\right|^2_
{\mathbb{H}^1}+\left|y_n\right|^2_{\mathbb{H}^1}\ra\\
&\,\,\,\,\,\,\,+\frac{C}{\eta^6}\left|X_n\right|^2_{\mathbb{H}^1}\la\left|Y_n\right|^4_
{\mathbb{H}^1}+\left|y_n\right|^4_{\mathbb{H}^1}\ra\\
&\le C\eta^2\left|\Delta
X_n\right|^2_{\mathbb{L}^2}+C_\eta\left|X_n\right|^2_{\mathbb{H}^1}\la
1+\left|Y_n\right|^4_{\mathbb{H}^1}\ra.
\end{aligned}
\end{equation}
Finally, using \eqref{ineq-interp}
we obtain
\begin{equation}\label{dc}
\begin{aligned}
|\las \Delta X_n, \left|\D
y_n\right|^2X_n\ras_{\mathbb{L}^2}|&\le\left|\Delta
X_n\right|_{\mathbb{L}^2}\left|\D y_n\right|_{\LL^\infty}\left|\D
y_n\right|_{\mathbb{L}^2}\left|X_n\right|_{\LL^\infty}\\
&\le C\eta^2\left|\Delta X_n\right|^2_{\mathbb{L}^2}\\
&\,\,\,\,\,\,\,+\frac{C}{\eta^2}\left|y_n\right|_{\mathbb{H}^1}\la\left|y_n\right|_{
\mathbb{H}^1}+\left|\Delta
y_n\right|_{\mathbb{L}^2}\ra\left|y_n\right|^2_{\mathbb{H}^1}\left|X_n\right|_{\mathbb{
L}^2}\left|X_n\right|_{\mathbb{H}^1}.
\end{aligned}
\end{equation}
Taking into account \eqref{da}, \eqref{db} and \eqref{dc} we obtain from
\eqref{d1}
\begin{equation}\label{est44}
\begin{aligned}
|\D X_n(t)|^2_{\mathbb{L}^2}+2\alpha\imp_0^t|\Delta
X_n|^2_{\mathbb{L}^2}\,ds&\le C\eta^2\imp_0^t\left|\Delta
X_n\right|^2_{\mathbb{L}^2}\,ds\\
&\,\,\,\,\,\,\,+C_\eta\sup_{r\le t}\la
1+\left|Y_n\right|^4_{\mathbb{H}^1}\ra\imp_0^t\left|X_n\right|^2_{\mathbb{H}^1}\,ds\\
&\,\,\,\,\,\,\,+C\eta^2\imp_0^t\left|\Delta
X_n\right|^2_{\mathbb{L}^2}\,ds+C_\eta\la\sup_{r\le
t}\left|X_n(r)\right|_{\mathbb{L}^2}\ra\la\sup_{r\le
t}\left|X_n(r)\right|_{\mathbb{H}^1}\ra\\
&\,\,\,\,\,\,\,+C\eta^2\imp_0^t\left|\Delta
X_n\right|^2_{\mathbb{L}^2}\,ds+C_\eta\la\sup_{r\le
t}\left|X_n(r)\right|_{\mathbb{L}^2}\ra\la\sup_{r\le
t}\left|X_n(r)\right|_{\mathbb{H}^1}\ra\\
&\,\,\,\,\,\,\,+C\eta^2\imp_0^t\left|\Delta
X_n\right|^2_{\mathbb{L}^2}\,ds+C_\eta\sup_{r\le
t}\left|X_n(r)\right|^2_{\mathbb{L}^2}\\
&\,\,\,\,\,\,\,+\sqrt{\eps_n}\left|\imp_0^t\las\D G\la Y_n\ra h,\Delta
X_n\ras_{\mathbb{L}^2}dW\right|\\
&\,\,\,\,\,\,\,+C\eps_n \imp_0^t\la 1+\left|\Delta
X_n\right|^2_{\mathbb{L}^2}\ra\,ds.
  \end{aligned}
\end{equation}
Choosing $\eta$ in such a way that $4C\eta^2=\alpha$ we obtain
\[\begin{aligned}
\E\la\sup_{t\in [0,T]}|\D X_n\la
t\wedge\tau_n\ra|^2_{\mathbb{L}^2}+\alpha\imp_0^{t\wedge\tau_n}|\Delta
X_n|^2_{\mathbb{L}^2}\,ds\ra&\le C_\eta\la
1+N^4\ra\E\imp_0^{\tau_n}\left|X_n\right|^2_{\mathbb{H}^1}\,ds\\
&\,\,\,\,\,\,\,+C_\eta\la 1+N\ra\E\sup_{t\in [0,T]}\left|X_n\la
t\wedge\tau_n\ra\right|_{\mathbb{L}^2}\\
&\,\,\,\,\,\,\,+\sqrt{\eps_n}\E\sup_{t\in
[0,T]}\left|\imp_0^{t\wedge\tau_n}\las\D G\la Y_n\ra h,\Delta
X_n\ras_{\mathbb{L}^2}dW\right|\\
&\,\,\,\,\,\,\,+C\eps_n \E\imp_0^T\la 1+\left|\Delta
X_n\right|^2_{\mathbb{L}^2}\ra\,ds.
  \end{aligned}
\]
By Theorem \ref{thm-reg} there exists a finite constant $C$, depending on $T$,
$\alpha$, $R$, $M_0$ and $h$ only, such that for each $n\in \mathbb{N}$
\[\E\imp_0^T|\Delta Y_n(s)|^2_{\mathbb{L}^2}\,ds\le C(T,\alpha,M,u_0,h),\]
hence invoking the Burkholder-Davis-Gundy inequality we find that
\[\begin{aligned}
\E\la\sup_{t\in [0,T]}|\D X_n\la
t\wedge\tau_n\ra|^2_{\mathbb{L}^2}+\alpha\imp_0^{t\wedge\tau_n}|\Delta
X_n|^2_{\mathbb{L}^2}\,ds\ra&\le C_\eta\la
1+N^4\ra\E\imp_0^{\tau_n}\left|X_n\right|^2_{\mathbb{H}^1}\,ds\\
&\,\,\,\,\,\,\,+C_\eta\la 1+N\ra\E\sup_{t\in [0,T]}\left|X_n\la
t\wedge\tau_n\ra\right|_{\mathbb{L}^2}\\
&\,\,\,\,\,\,\,+C(1+N)\sqrt{\eps_n}.
  \end{aligned}
\]
Finally, Lemma \ref{lem24} follows from Lemma \ref{lem23}.
\end{proof}

We will conclude this section with the promised proof of Lemma
\ref{lem22}.

\begin{proof}[Proof of Lemma \ref{lem22}] We will use the same notation as in
the proof of Lemma \eqref{lem24}. Let $\delta>0$ and $\nu>0$. Invoking part (2)
of Theorem (\ref{thm-weak solution}) we can find $N>|M_0|_{\mathbb{H}^1}$ such that
\[\frac{1}{N}\sup_{n\in \mathbb{N}}\mathbb{E}\sup_{t\in
[0,T]}|Y_n(t)|_{\mathbb{H}^1}<\frac{\nu}{2}.\]
Then invoking Lemma \ref{lem24} we find that for all $n$ sufficiently large
\begin{eqnarray*}
&&\!\!\!\!\!\!\!\!\!\!\!\!\!\!\mathbb{P}\left(\sup_{t\in
[0,T]}|Y_n(t)-u_n(t)|^2_{\mathbb{H}^1}+\imp_0^T|Y_n-u_n|^2_{\xone}\,
ds\geq\delta\right)\\
&\le&\mathbb{P}\left(\sup_{t\in
[0,T]}|Y_n(t\wedge\tau_n)-u_n(t\wedge\tau_n)|^2_{\mathbb{H}^1}+\imp_0^{\tau_n}
|Y_n-u_n|^2_{\xone}\,ds\geq\delta,\,\tau_n=T\right)\\
&&{}+\mathbb{P}\left(\sup_{t\in [0,T]}|Y_n(t)|_{\mathbb{H}^1}\geq N\right)\\
&\le&\frac{1}{\delta}\mathbb{E}\left(\sup_{t\in
[0,T]}|Y_n(t\wedge\tau_n)-u_n(t\wedge\tau_n)|^2_{\mathbb{H}^1}+\imp_0^{\tau_n}
|Y_n-u_n|^2_{\xone}\,ds\right)\\
&&{}+\frac{1}{N}\mathbb{E}\sup_{t\in [0,T]}|Y_n(t)|_{\mathbb{H}^1}\\
&<&\nu.
\end{eqnarray*}
\end{proof}

\section{Application to a model of a ferromagnetic needle} \label{sec:applyldp}

\noindent  In this section we will use the large deviation principle established
in the previous section to investigate the dynamics of a stochastic
Landau-Lifshitz model of magnetization in a needle-shaped particle. Here the
shape anisotropy energy is crucial. When there is no applied field and no noise
in the field, the shape anisotropy energy gives rise to two locally stable
stationary states of opposite magnetization. We add a small noise term to the
field and use the large deviation principle to show that noise induced
magnetization reversal occurs and to quantify the effect of material parameters
on sensitivity to noise.

\smallskip

\noindent The axis of the needle is represented by the interval $\domain$ and at
each $x\in\domain$ the magnetization $u(x)\in\mathbb S^2$ is assumed to be
constant over the cross-section of the needle. We define the total magnetic
energy of magnetization $u\in\HH^1$ of the needle by
\begin{equation}\label{exe1}
E_t(u)=\frac{1}{2}\imp_\domain |\D
u(x)|^2\,dx+\beta\imp_{\domain}\Phi(u(x))\,dx-\imp_\domain\mathcal K(t,x)\cdot
u(x)\, dx,
\end{equation}
where
\[\Phi(u)=\Phi\la u_1,u_2,u_3\ra=\frac{1}{2}\la u_2^2+u_3^2\ra,\]
$\beta$ is the positive real shape anisotropy parameter and $\mathcal{K}$ is the
externally applied magnetic field, such that  $\mathcal K(t)\in\mathbb{H}$ for
each $t$.
% The negative of the derivative of the energy is
%\begin{equation}\label{exe1a}
%-DE(y)=\Delta y-\beta\left(\begin{array}{c}0\\y_2\\y_3\end{array}\right)+\mathcal{K} \;\;\;\mbox{ for }y\in \xone,
%\end{equation}
%which is the total (deterministic) magnetic field.

With this magnetic energy, the deterministic Landau-Lifshitz equation becomes:
\begin{equation}\label{exe2}
\frac{ \partial y}{ \partial t}(t)= y\times \Delta y-\alpha
y\times(y\times\Delta y)+G(y)\la -\beta f(y)+\mathcal K(t)\ra
\end{equation}
where $f(y)= \D \Phi(y)$, $y\in\mathbb{R}^3$.
We assume, as before, that the initial state $u_0\in\mathbb{H}^1$ and
$|u_0(x)|_{\mathbb{R}^3}=1$ for all $x\in\domain$. We also assume that the
applied field $\mathcal{K}(t):\domain\rightarrow\mathbb{R}^3$ is constant on
$\domain$ at each time $t$. Equation~(\ref{exe2}) has nice features: the
dynamics of the solution can be studied using elementary techniques and, when
the externally applied field $\mathcal{K}$ is zero, the equation has two stable
stationary states, $\zeta_{+}=(1,0,0)$ and $\zeta_{-}=(-1,0,0)=-\zeta_{+}$. In
what follow we will abuse the notation as by $\zeta_{\pm}$. We will also denote
a constant function $\domain \ni x \mapsto \zeta_{\pm} \in \mathbb{S}^2$ which
obviously belongs to $\mathbb{H}^1$.

\smallskip

\noindent We now outline the structure of this example.
\begin{comment}
Firstly, recall, see \eqref{ineq-interp},  that $\mathbb{H}^1$ is continuously
embedded in $C(\domain;\mathbb{R}^3)$, and by \cite[Theorem 5.8 and its
proof]{AF}
\begin{equation}\label{interpol}
\sup_{x\in\domain}|u(x)|^2\leq k^2|u|_{\mathbb{H}}|u|_{\mathbb{H}^1}\;\;\;\mbox{
for all }u\in\mathbb{H}^1;
\end{equation}
with $k=2(1\vee 1/\sqrt{l(\domain)})$, where $|\domain|$ is the length of the
interval $\domain$.
\end{comment}
\noindent In Proposition \ref{stablestates}, we show that if the applied field
$\mathcal{K}$ is zero and the initial state $y_0$ satisfies
\[\left|y_0-\zeta_{\pm}\right|_{\mathbb{H}^1}<\frac{1}{2k^2\sqrt{|\domain|}}
\frac{\alpha}{1+2\alpha},\]
then the solution $y(t)$ of (\ref{exe2}) converges to $\zeta_{\pm}$ in
$\mathbb{H}^1$ as $t$ goes to $\infty$. In Lemma \ref{stablestates2}, we show
that if $\lambda$  exceeds a certain value (depending on $\alpha$ and $\beta$) and the
applied field is
$\mathcal{K}=\lambda \mathbf{m}+\beta f(\mathbf{m})$ and
$|y_0-\mathbf{m}|_{\mathbb{H}^1}<\frac{1}{k}$, then
$y(t)$ converges in $\mathbb{H}^1$ to $\mathbf{m}$ as $t$
goes to $\infty$. Lemma~\ref{stablestates2} is used to show that, given
$\delta\in(0,\infty)$ and $T\in (0,\infty)$, there is a piecewise constant (in
time) externally applied field, $\mathcal{K}$, which drives the magnetization
from the initial state $\zeta_{-}$ to the $\mathbb{H}^1$-ball centred at
$\zeta_{+}$ and of radius $\delta$ by time $T$; in short, in the
deterministic system, this applied field causes magnetization reversal by time
$T$(see Definition~\ref{defmagrev}). What we are really interested in is the
effect of adding a small noise term to the field. We will show that if
$\mathcal{K}$ is zero but a noise term multiplied by $\sqrt{\eps}$ is added to
the field, then the solution of the resulting stochastic equation exhibits
magnetization reversal by time $T$ with positive probability for all
sufficiently small positive $\eps$. This result, in Proposition~\ref{magrev}, is
obtained using the lower bound of the large deviation principle. Finally, in
Proposition~\ref{expuppbd}, the upper bound of the large deviation principle is
used: we obtain an exponential estimate of the probability that, in time
interval $[0,T]$, the stochastic magnetization leaves a given
$\mathbb{H}^1$-ball centred at the initial state $\zeta_{-}$ and of
radius less than or equal to
$\frac{1}{2k^2\sqrt{|\domain|}}\frac{\alpha}{1+2\alpha}$. This estimate
emphasizes the importance of a large value of $\beta$ for reducing the
disturbance in the magnetization caused by noise in the field.

\subsection{Stable stationary states of the deterministic equation}

\noindent In this subsection, we identify stable stationary states of the
deterministic equation (\ref{exe2}) when the applied field $\mathcal{K}$ does
not vary with time.

\smallskip

\noindent Let $\zeta\in\mathbb S^2$. Since the time derivative $\frac{dy}{dt}$
of the solution $y$ to $(\ref{exe2})$, belongs to $L^2(0,T;\mathbb{H})$ and $y$
belongs to $L^2(0,T;\xone)$, we have for all $t\geq 0$:
\begin{eqnarray} \label{e1a}
|y(t)-\zeta|^2_{\mathbb{H}}&=&|y_0-\zeta|^2_{\mathbb{H}}\\
&+&2\imp_0^t\las y-\zeta, y\times\Delta y-\alpha y\times(y\times\Delta
y)\right.\nonumber\\
&&\mbox{\hspace{30mm}}{}\left. +G(y)(-\beta f(y)+\mathcal
K\ras_{\mathbb{H}}\,ds\,\, \nonumber\\
\label{e1}
&=&|y_0-\zeta|^2_{\mathbb{H}}+2\imp_0^t\las -\zeta, y\times\Delta y-\alpha
y\times(y\times\Delta y)\right.\nonumber\\
&&\mbox{\hspace{30mm}}{}\left. +G(y)(-\beta f(y)+\mathcal
K\ras_{\mathbb{H}}\,ds\,\, \nonumber
\end{eqnarray}
\begin{eqnarray}\label{e2}
|\D y(t)|^2_{\mathbb{H}}
%&=&|\D y_0|^2_{\mathbb{H}}-2\imp_0^t\langle \Delta y(s), y(s)\times\Delta
y(s)+y(s)\times(-\beta
f(y(s))\dela{\vthree{0}{y_2(s)}{y_3(s)}}+\mathcal{K})\nonumber\\
%&&\mbox{\hspace{30mm}}{}-\alpha y(s)\times(y(s)\times\Delta y(s))\nonumber\\
%&&\mbox{\hspace{30mm}}{}-\alpha y(s)\times(y(s)\times(-\beta
f(y(s))\dela{\vthree{0}{y_2(s)}{y_3(s)}}+\mathcal{K}))\rangle_{\mathbb{H}}\,
ds\nonumber\\
&=&|\D y_0|^2_{\mathbb{H}}
-2\imp_0^t\langle \Delta y, G(y)\la-\beta f(y)+\mathcal K\ra\\
&&\mbox{\hspace{30mm}}{}-\alpha y\times(y\times\Delta
y)\rangle_{\mathbb{L}^2}\,ds\nonumber
\end{eqnarray}
\begin{lemma}\label{uandz}
 Let $u\in\HH^1$ be such that $u(x)\in\mathbb S^2$ and
\[|u-\zeta_{\pm}|_{\mathbb{H}^1}\le\frac{1}{2k^2\sqrt{|\domain|}}\frac{\alpha}{
1+2\alpha}.\]
 Then for all $x\in\domain$

\begin{enumerate}
\item $\frac{1-u_1^2(x)}{u_1^2(x)}+\alpha\frac{(1-u^2_1(x))^2}{u_1^2(x)}-\alpha
u_1^2(x)\le 0\;\;\;$ for all $x\in\domain$,
\item $\langle u(x),\zeta_{\pm}\rangle\geq\frac{3}{4}\;\;\;$ and
\item $\frac{7}{8}|u(x)\zeta_{\pm}|^2\le |u(x)\times \zeta_{\pm}|^2\;\;\;$.
\end{enumerate}
\end{lemma}
\begin{proof}
By \eqref{ineq-interp}
\begin{eqnarray}
\nonumber
\sup_{x\in\domain}|u(x)-\zeta_{\pm}|^2&\le&k^2|u-\zeta_{\pm}|_{\mathbb{H}}
|u-\zeta_{\pm}|_{\mathbb{H}^1},\\
\label{ineqL1}&\le&k^2\,2\sqrt{|\domain|}\,\frac{1}{2k^2\sqrt{|\domain|}}\frac{
\alpha}{1+2\alpha}=\frac{\alpha}{1+2\alpha}.
\end{eqnarray}
Invoking (\ref{ineqL1}), we find that
\begin{equation}\label{ineqL2}
u_1^2(x)=1-(u_2^2(x)+u_3^2(x))\geq
1-|u(x)-\zeta_{\pm}|^2\geq\frac{1+\alpha}{1+2\alpha},\quad x\in\domain.
\end{equation}
Hence one can use (\ref{ineqL2}) and straightforward algebraic manipulations to
verify that
\[
\frac{1-u^2_1(x)}{u_1^2(x)}+\alpha\frac{(1-u_1^2(x))^2}{u_1^2(x)}-\alpha
u_1^2(x)\le 0.
\]
Statements 2 and 3 of Lemma~\ref{uandz} follow easily from (\ref{ineqL1}).
\begin{comment}
 and the elementary equalities which hold for all $x\in\domain$:
\[2\langle
u(x),\addc{-\zeta_{\pm}}\rangle=2-|u(x)\addc{-\zeta_{\pm}}|^2\quad\mathrm{and}
\quad |u(x)\times(\addc{-\zeta_{\pm}})|^2=1-\langle
u(x),\addc{-\zeta_{\pm}}\rangle^2.\]
\end{comment}
\end{proof}
\begin{proposition}\label{stablestates} Let the applied field $\mathcal{K}$ be
zero and let $y_0\in\HH^1$ satisfy
\begin{equation}\label{estab1}
|y_0-\zeta_{\pm}|_{\mathbb{H}^1}<\frac{1}{2k^2\sqrt{|\domain|}}\frac{\alpha}{
1+2\alpha}.
\end{equation}
Let the process $y$ be the solution to \eqref{exe2}.
Then $y(t)$ converges to $\zeta_{\pm}$ in $\mathbb{H}^1$ as $t \to \infty$.
\end{proposition}
\begin{proof}
Using some algebraic manipulation and the fact that $\langle \D
y(s),y(s)\rangle=0$ a.e. on $\domain$ for each $s\geq 0$, one may simplify
equations (\ref{e1}) and (\ref{e2}).

\noindent We obtain from (\ref{e1}):
\begin{eqnarray}
\nonumber
|y(t)-\zeta_{\pm}|^2_{\mathbb{H}}&=&|y_0-\zeta_{\pm}|^2_{\mathbb{H}}
-2\alpha\imp_0^t\imp_\domain|\D y(s)|^2\langle
y(s),-\zeta_{\pm}\rangle\,dx\,ds\\
\label{uminuszbound}&&{}-2\alpha\beta\imp_0^t\imp_{\domain}\langle
y(s),-\zeta_{\pm}\rangle|y(s)\times \zeta_{+}|^2\,dx\,ds\;\;\;\mbox{ }\forall
t\geq 0,
\end{eqnarray}
and
\begin{equation}\label{Dubound}
|\D y(t)|^2_{\mathbb{H}}=|\D
y_0|^2_{\mathbb{H}}-2\alpha\imp_0^t|y(s)\times\Delta
y(s)|^2_{\mathbb{H}}\,ds+2\beta\imp_0^t\imp_{\domain}R(s)\,dx\,ds\;\;\;\mbox{
}\forall t\geq 0,
\end{equation}
where
\begin{eqnarray}
\label{Rsdef} R&=&\D y_1(y_3\D y_2-y_2\D y_3)+\alpha(\D
y_1)^2-\alpha y_1^2|\D y|^2\\
\nonumber &=&\frac{-y_2\D y_2-y_3\D y_3}{y_1}(y_3\D y_2-y_2\D
y_3)\\
\nonumber&&{}+\alpha(1-y_1^2)\left(\frac{y_2\D y_2+y_3\D
y_3}{y_1}\right)^2-\alpha y_1^2((\D y_2)^2+(\D y_3)^2.
\end{eqnarray}
Define
\[\tau=\inf\left\{t\geq 0:\,
|y(t)-\zeta_{\pm}|_{\mathbb{H}^1}\geq\frac{1}{2k^2\sqrt{|\domain|}}\frac{\alpha}
{1+2\alpha}\right\}.\]
 Then, by our choice of $y_0$, $\tau>0$. For each $s\in[0,\tau)$, $y(s)$
satisfies the hypotheses of Lemma~\ref{uandz}, hence
 and
 \[y(s)(x)\cdot\la-\zeta_{\pm}\ra\geq \frac{3}{4},\quad x\in\domain,\]
  \[|y(s)(x)\times\la-\zeta_{\pm}\ra|^2\geq
\frac{7}{8}|y(s)(x)-\zeta_{\pm}|^2,\quad x\in\domain.\]
 and, invoking the Cauchy-Schwartz inequality
\begin{equation}\label{Rsineq}
R\le\left(\frac{1-y_1^2}{y_1^2}+\alpha\frac{(1-y^2_1)^2}{y^2_1}-\alpha
y_1^2\right)((\D y_2)^2+(\D y_3)^2)\le 0,\quad x\in\domain.
\end{equation}
   Consequently, from (\ref{uminuszbound}) and (\ref{Dubound}) we deduce that
the functions $|y(\cdot)-\zeta_{\pm}|^2_{\mathbb{H}}$ and $|\D
y(\cdot)|^2_{\mathbb{H}}$ are nonincreasing on $[0,\tau)$. Furthermore, we have
\begin{equation}
\label{uminuszbd2}
|y(t)-\zeta_{\pm}|^2_{\mathbb{H}}\le|y_0-\zeta_{\pm}|^2_{\mathbb{H}}-\frac{3}{2}
\alpha\imp_0^t|\D
y(s)|^2_{\mathbb{H}}\,ds-\frac{21}{16}\alpha\beta\imp_0^t|y(s)-\zeta_{\pm}|^2_{
\mathbb{H}}\,ds,\quad t<\tau,
\end{equation}
and
\begin{equation}\label{Dubound2}
|\D y(t)|^2_{\mathbb{H}}\le |\D y_0|^2_{\mathbb{H}},\quad t<\tau.
\end{equation}
Suppose, to get a contradiction, that $\tau<\infty$. Then, from
(\ref{uminuszbd2}) and (\ref{Dubound2}), we have
\[
|y(\tau)-\zeta_{\pm}|_{\mathbb{H}^1}\le
|y_0-\zeta_{\pm}|_{\mathbb{H}^1}<\frac{1}{2k^2\sqrt{l(\domain)}}\frac{\alpha}{
1+2\alpha},
\]
which contradicts the definition of $\tau$. Therefore, $\tau=\infty$. Since
(\ref{uminuszbd2}) holds for all $t\geq 0$, we have
\[\imp_0^\infty|\D
y(s)|^2_{\mathbb{H}}\,ds+\imp_0^\infty|y(s)-\zeta_{\pm}|^2_{\mathbb{H}}\,
ds<\infty.\]
Since both integrands are nonincreasing
\[\lim_{t\to\infty}\la|\D
y(t)|_{\mathbb{H}}+|y(t)-\zeta_{\pm}|_{\mathbb{H}}\ra=0.\]
\end{proof}
%\medskip
%\noindent It is straightforward to show that when $\mathcal{K}=0$ the only
other uniform stationary states of equation (\ref{DetEq}) are points of the form
$(0,y_2,y_3)\in\Rthree$, where $y^2_2+y^2_3=1$; however, such a point,
$(0,y_2,y_3)$, is not a stable stationary state because any given
$\mathbb{H}^1$-ball centred at the point contains another uniform state
$(\underline{y}_1,\underline{y}_2,\underline{y}_3)\in\mathbb{R}^3$ with lower
energy (that is, $\underline{y}_2^2+\underline{y}_3^2<1$) and, by (\ref{e4}),
energy is nonincreasing.
We will show next, that if the applied field has sufficiently large magnitude,
then there exists a stable stationary state that is roughly in the direction of
the applied field.

\begin{lemma}\label{stablestates2} Assume that
$\mathbf{m}\in \mathbb S^2$ and
 a real number $\lambda$ satisfies
\begin{equation}\label{magcurlyH}
\lambda
>\left(\frac{4\beta+4\alpha\beta}{3\alpha}\vee\frac{2\beta+4\alpha\beta-\alpha}{
\alpha}\right).
\end{equation}
Let the applied field be\footnote{Note that a constant function $\mathbf{m}$ is
a stationary solution to the problem \eqref{exe2}.} \[
\mathcal{K}:=\lambda \mathbf{m}+\beta  f(\mathbf{m}).
\]
Let $y$ be a solution to the problem \eqref{exe2}  with initial data $y_0$
satisfying $|y_0-\mathbf{m}|_{\mathbb{H}^1}<\frac{1}{k}$. Then
\begin{equation}\label{expdecay}
\left|y(t)-\mathbf{m}\right|_{\mathbb{H}^1}\le
\left|y_0-\mathbf{m}\right|_{\mathbb{H}^1}\mathrm{e}^{-\frac{1}{2}\gamma
t}\;\;\mbox{ }\forall t\geq 0,
\end{equation}
where
\[\gamma:=(\alpha \lambda +\alpha-2\beta-4\alpha\beta)\wedge\la\frac{3}{2}\alpha
\lambda -2\beta-2\alpha\beta\ra>0\]
 is positive, by condition (\ref{magcurlyH}).
\end{lemma}
\begin{proof}
\noindent We have, from (\ref{e1a}) and (\ref{e1}) with $\zeta$ replaced by
$\mathbf{m}$:
\begin{eqnarray}\label{yminusH}
&&\!\!\!\!\!\!\!\!\!\!\!\!|y(t)-\textstyle{\mathbf{m}}|^2_\mathbb{H}\\
&=&|y_0-\textstyle{\mathbf{m}}|^2_\mathbb{H}+2\displaystyle{\imp_0^t}\langle
y(s)-\textstyle{\mathbf{m}}, y(s)\times(\mathbf{m}-\beta f(y- \mathbf{m} )
)\rangle_\mathbb{H}\,ds\nonumber\\
&&\mbox{\hspace{30mm}}{}+2\alpha\imp_0^t\langle\Delta y(s),
y(s)\times(y(s)\times\textstyle{\mathbf{m}})\rangle_\mathbb{H}\,ds\nonumber\\
&&\mbox{\hspace{30mm}}-2\alpha\imp_0^t\langle
y(s)\times\textstyle{\mathbf{m}},y(s)\times(\mathbf{m}-\beta f(y- \mathbf{m}
) )\rangle_\mathbb{H}\,ds\nonumber\\
&=&|y_0-\textstyle{\mathbf{m}}|^2_\mathbb{H}-2\beta\displaystyle{\imp_0^t}
\langle y-\textstyle{\mathbf{m}}, y\times f(y- \mathbf{m} )
\rangle_{\mathbb{H}}\,ds\nonumber\\
&&\mbox{\hspace{30mm}}{}-2\alpha\imp_0^t\imp_\domain|\D y|^2(y\cdot
\textstyle{\mathbf{m}})\,dx\,ds\nonumber\\
&&\mbox{\hspace{30mm}}{}-2\alpha\imp_0^t|y\times\textstyle{\mathbf{m}}
|^2_\mathbb{H}\,ds\nonumber\\
&&\mbox{\hspace{30mm}}{}+2\alpha\beta\imp_0^t\langle
y\times\textstyle{\mathbf{m}},y\times f(y- \mathbf{m} )
\rangle_\mathbb{H}\,ds\;\;\mbox{ }\forall t\geq 0.\nonumber
\end{eqnarray}
From (\ref{e2}) we have:
\begin{eqnarray}\label{Dery}
|\D y(t)|^2_{\mathbb{H}}&=&|\D
y_0|^2_{\mathbb{H}}-2\imp_0^t\langle\Delta y(s), y(s)\times(\mathbf{m}-\beta
f(y- \mathbf{m} ) )\\
&&\mbox{\hspace{30mm}}{}-\alpha y(s)\times(y(s)\times\Delta y(s))\nonumber\\
&&\mbox{\hspace{30mm}}{}-\alpha y(s)\times(y(s)\times(\mathbf{m}-\beta f(y-
\mathbf{m} ) ))\rangle_\mathbb{H}\,ds\nonumber\\
\end{eqnarray}

\begin{eqnarray}\label{Dery2}
&=&|\D y_0|^2_{\mathbb{H}}+2\beta\imp_0^t\langle\Delta y, y\times f(y-
\mathbf{m} ) \rangle_\mathbb{H}\,ds\nonumber\\
&&\mbox{\hspace{0mm}}{}-2\alpha\imp_0^t|y\times\Delta
y|^2_\mathbb{H}\,ds\nonumber\\
&&\mbox{\hspace{0mm}}{}-2\alpha\imp_0^t\imp_\domain|\D y|^2(y\cdot
\textstyle{\mathbf{m}})\,dx\,ds\nonumber\\
&&\mbox{\hspace{0mm}}{}-2\alpha\beta\imp_0^t\langle\Delta y, y\times(y\times
f(y- \mathbf{m} ) )\rangle_\mathbb{H}\,ds\;\;\mbox{ }\forall t\geq 0. \nonumber
\end{eqnarray}
Define
\begin{equation}\label{tau1def}
\tau_1:=\inf\{t\geq 0: |y(t)-\textstyle{\mathbf{m}}|_{\mathbb{H}^1}\geq
\frac{1}{k}\}.
\end{equation}
By our choice of $y_0$, $\tau_1$ is greater than zero. Observe that
\begin{equation}\label{e18a}
\sup_{x\in\domain}|y(t)(x)-\textstyle{\mathbf{m}}|_\Rthree<1\;\;\mbox{ for all
}t<\tau_1.
\end{equation}
It is easy to check that for every $t<\tau_1$
%\begin{enumerate}
%\item
\begin{equation}\label{e18}
\frac{3}{4}|y(t)-\textstyle{\mathbf{m}}|^2_\mathbb{H}\le
|y(t)\times\textstyle{\mathbf{m}}|^2_\mathbb{H}\le
|y(t)-\textstyle{\mathbf{m}}|^2_\mathbb{H},
\end{equation}
%\item and
and
\begin{equation}\label{e19}
y(t,x)\cdot\textstyle{\mathbf{m}}\geq\frac{1}{2},\quad x\in\domain.
\end{equation}
%\end{enumerate}
Adding equalities (\ref{yminusH}) and (\ref{Dery}) we obtain for $t>0$
\begin{eqnarray}
&&\!\!\!\!\!\!\!\!\!\!\!\!|y(t)-\textstyle{\mathbf{m}}|^2_{\mathbb{H}^1}
\nonumber\\
&=&|y_0-\textstyle{\mathbf{m}}|^2_{\mathbb{H}^1}-4\alpha\displaystyle{
\imp_0^t\imp_\domain}|\D y|^2\la
y\cdot\textstyle{\mathbf{m}}\ra\,dx\,ds\nonumber\\
&&\mbox{\hspace{0mm}}{}+2\beta\imp_0^t\langle\Delta y, y\times f(y- \mathbf{m}
) \rangle_\mathbb{H}\,ds\nonumber\\
&&\mbox{\hspace{0mm}}{}-2\alpha\beta\imp_0^t\langle\Delta y, y\times(y\times
f(y- \mathbf{m} ) )\rangle_\mathbb{H}\,ds\nonumber\\
&&\mbox{\hspace{0mm}}{}-2\alpha\imp_0^t|y\times\textstyle{\mathbf{m}}
|^2_\mathbb{H}\,ds\nonumber\\
&&\mbox{\hspace{0mm}}{}-2\beta\imp_0^t\langle y-\textstyle{\mathbf{m}}, y\times
f(y- \mathbf{m} ) \rangle_\mathbb{H}\,ds\nonumber\\
&&\mbox{\hspace{0mm}}{}+2\alpha\beta\imp_0^t\langle
y\times\textstyle{\mathbf{m}}, y\times f(y- \mathbf{m} )
\rangle_\mathbb{H}\,ds\nonumber\\
&&\mbox{\hspace{0mm}}{}-2\alpha\imp_0^t|y\times\Delta y|^2_\mathbb{H}\,ds\,.
\label{e20}
\end{eqnarray}
Therefore for every $t<\tau_1$
\begin{eqnarray}
\!\!\!\!\!\!\!\!\!\!\!\!|y(t)-\textstyle{\mathbf{m}}|^2_{\mathbb{H}^1}
\nonumber&\le&
|y_0-\textstyle{\mathbf{m}}|^2_{\mathbb{H}^1}
-(2\alpha-2\beta-4\alpha\beta)\displaystyle{\imp_0^t}|\D
y|^2_\mathbb{H}\,ds\nonumber\\
&&\mbox{\hspace{0mm}}{}-(\textstyle{\frac{3}{2}}
\alpha-2\beta-2\alpha\beta)\displaystyle{\imp_0^t}|y-\textstyle{\mathbf{m}}
|^2_\mathbb{H}\,ds\nonumber\\
&&\mbox{\hspace{0mm}}{}-2\alpha\imp_0^t|y\times\Delta y|^2_\mathbb{H}\,ds\,,
\label{e21}
\end{eqnarray}
where we used (\ref{e18a}), (\ref{e18}) and (\ref{e19}). Because of hypothesis
(\ref{magcurlyH}), the two expressions $(2\alpha-2\beta-4\alpha\beta)$ and
$(\frac{3}{2}\alpha-2\beta-2\alpha\beta)$ on the right hand side of (\ref{e21})
are positive numbers.

\noindent Suppose, to get a contradiction, that $\tau_1<\infty$. Then, from
(\ref{e21}), we have
\[
|y(\tau_1)-\textstyle{\mathbf{m}}|_{\mathbb{H}^1}\le
|y_0-\textstyle{\mathbf{m}}|_{\mathbb{H}^1}<\frac{1}{k},
\]
which contradicts the definition of $\tau_1$ in (\ref{tau1def}).
Hence $\tau_1=\infty$. It now follows from (\ref{e21}) that
\begin{eqnarray}
\imp_0^\infty|\D y(s)|^2_{\mathbb{H}}\,ds&<&\infty, \label{e23}\\
\imp_0^\infty|y(s)-\textstyle{\mathbf{m}}|^2_\mathbb{H}\,ds&<&\infty
\label{e24}\\
\mbox{ and }\;\;\imp_0^\infty|y(s)\times\Delta
y(s)|^2_{\mathbb{H}}\,ds&<&\infty. \label{e25}
\end{eqnarray}
From (\ref{e20}) and these three inequalities, the function
$t\in[0,\infty)\mapsto |y(t)-\textstyle{\mathbf{m}}|^2_{\mathbb{H}^1}$ is
absolutely continuous and, for almost every $t\ge 0$, its derivative is:
\begin{eqnarray}
\frac{d}{dt}|y-\textstyle{\mathbf{m}}|^2_{\mathbb{H}^1}(t)
&=&-4\alpha\imp_\domain|\D y(t)|^2(y(t)\cdot
\textstyle{\mathbf{m}})\,dx\nonumber\\
&&{}+2\beta\langle\Delta y(t), y(t)\times
f(y-\mathbf{m})\rangle_\mathbb{H}\nonumber\\
&&{}-2\alpha\beta\langle\Delta y(t), y(t)\times(y(t)\times
f(y-\mathbf{m}))\rangle_\mathbb{H}\nonumber\\
&&{}-2\alpha|y(t)\times\textstyle{\mathbf{m}}|^2_\mathbb{H}\nonumber\\
&&{}-2\beta\big\langle y(t)-\textstyle{\mathbf{m}}, y(t)\times
f(y-\mathbf{m})\big\rangle_\mathbb{H}\nonumber\\
&&{}+2\alpha\beta\langle y(t)\times\textstyle{\mathbf{m}}, y(t)\times
f(y-\mathbf{m})\rangle_\mathbb{H}\nonumber\\
&&{}-2\alpha|y(t)\times\Delta y(t)|^2_\mathbb{H}\nonumber\\
&\le&-(2\alpha-2\beta-4\alpha\beta)|\D y(t)|^2_\mathbb{H}\nonumber\\
&&{}-(\textstyle{\frac{3}{2}}\alpha-2\beta-2\alpha\beta)|y(t)-\textstyle{
\mathbf{m}}|^2_\mathbb{H}\nonumber\\
&&{}-2\alpha|y(t)\times\Delta y(t)|^2_\mathbb{H}\nonumber\\
&\le&-\gamma|y(t)-\textstyle{\mathbf{m}}|^2_{\mathbb{H}^1}, \label{e26}
\end{eqnarray}
where
\[\gamma:=(\alpha+\alpha \lambda -2\beta-4\alpha\beta)\wedge\la\frac{3}{2}\alpha
\lambda -2\beta-2\alpha\beta\ra>0.\]
Now the lemma follows by a standard argument.
\end{proof}

\subsection{Noise induced instability and magnetization reversal}
In Proposition~\ref{stablestates} we showed that the states $\zeta_{+}$
and $\zeta_{-}$ are stable stationary states of the deterministic
Landau-Lifshitz equation $(\ref{exe2})$ when the externally applied field
$\mathcal{K}$ is zero. In this section we show that a small noise term in the
field may drive the magnetization from the initial state $\zeta_{-}$ to
any given $\mathbb{H}^1$-ball centred at $\zeta_{+}$ in any given time
interval $[0,T]$. We also find an exponential upper bound for the probability
that small noise in the field drives the magnetization outside a given
$\mathbb{H}^1$-ball centred at the initial state $\zeta_{-}$ in time
interval $[0,T]$. Firstly we need a definition.
\begin{definition}\label{defmagrev} Let $\delta$ be a given small positive real
number. Suppose that the initial magnetization is $\zeta_{-}$ and that at
some time $T$ the magnetization lies in the open $\mathbb{H}^1$-ball centred at
$\zeta_{+}$ and of radius $\delta$. Then we say that \emph{magnetization
reversal} has occurred by time $T$.
\end{definition}
\noindent
We consider a stochastic equation for the magnetization, obtained by setting
$\mathcal{K}$ to zero and adding a three dimensional noise term to the field.
Denoting the  magnetization by $Y$, the equation is:
\begin{equation}\label{e27}
\left.\begin{array}{rcl}dY&=&\la Y\times\Delta Y-\alpha Y\times(Y\times\Delta
Y)+\beta G(Y)f(Y)\ra\, dt\\
\\
&&{}+\sqrt{\eps}G(Y)B\circ dW(t)\\
Y(0)&=&\zeta_{-}.\end{array}\right\}
\end{equation}
In (\ref{e27}), we assume that the vectors $\newe_1,\newe_2,\newe_3\in\R^3$ are
linearly independent.
The parameter $\eps>0$ corresponds to the `dimensionless temperature' parameter
appearing in the following stochastic differential equation (\ref{krvstochmodel}) of Kohn,
Reznikoff and Vanden-Eijnden \cite{KRV}:
\begin{equation}\label{krvstochmodel}
\dot{m}=m\times(g+\varepsilon^\frac{1}{2}\sqrt{\textstyle{\frac{2\alpha}{1+\alpha^2}}}\dot{W})-\alpha m\times(m\times(g+\varepsilon^\frac{1}{2}\sqrt{\textstyle{\frac{2\alpha}{1+\alpha^2}}}\dot{W})),
\end{equation}

\smallskip

\noindent Fix $T>0$. There is no deterministic applied field in (\ref{e27}) but,
as we will see, the lower bound of the large deviation principle satisfied by
the solutions $Y^\eps$ ($\eps\in(0,1)$) of (\ref{e27}) implies that, for all
sufficiently small positive $\eps$, the probability of magnetization reversal by
time $T$ is positive.

\smallskip

\noindent Firstly, we shall use Lemma~\ref{stablestates2} to construct a
piecewise constant (in time) deterministic applied field, $\mathcal{K}$, such
that the solution $y$ of (\ref{exe2}), with initial state
$\zeta_{-}=(-1,0,0)$, undergoes magnetization reversal by time $T$.

\noindent Take points $u^i\in\mathbb S^2$, $i=0,1,\ldots,N$, such that
$u^0=\zeta_{-}$ and $u^N=\zeta_{+}$ and
\[
|u^i-u^{i+1}|_{\mathbb{H}^1}=|u^i-u^{i+1}|_\Rthree\sqrt{|\domain|}<\frac{1}{k}
\;\;\mbox{ for }i=0,1,\ldots,N-1.
\]
Let
\[\eta:=\min\left\{\frac{1}{k}-|u^i-u^{i+1}|_{\mathbb{H}^1}:
i=1,\ldots,N-1\right\}\wedge\frac{\delta}{2}.\] Using Lemma~\ref{stablestates2},
we can take the applied field to be
\begin{equation}\label{e28}
\mathcal{K}(t):=\sum_{i=0}^{N-1}\,1_{(i\frac{T}{N}, (i+1)\frac{T}{N}]}(t)\left(R
u^{i+1}+\beta f(u^{i+1})\right),\;\;\mbox{ }t\geq 0,
\end{equation}
with the positive real number $R$ chosen to ensure that, as $t$ varies from
$i\frac{T}{N}$ to $(i+1)\frac{T}{N}$, $y(t)$ starts at a distance of less than
$\eta$ from $u^i$ (i.e. $|y(i\frac{T}{N})-u^i|_{\mathbb{H}^1}<\eta$) and moves
to a distance of less than $\eta$ from $u^{i+1}$ (i.e.
$|y((i+1)\frac{T}{N})-u^{i+1}|_{\mathbb{H}^1}<\eta$). Specifically, we take
$R\in(0,\infty)$ such that
\[
\frac{1}{k}\mathrm{e}^{-\frac{1}{2}[(\alpha
R+\alpha-2\beta-4\alpha\beta)\wedge(\frac{3}{2}\alpha
R-2\beta-2\alpha\beta)]\frac{T}{N}}<\eta.
\]
For each $i=0,1,\ldots,N-1$, let $\phi^{i+1}=(\phi^{i+1}_1, \phi^{i+1}_2,
\phi^{i+1}_3)\in\mathbb{R}^3$ be the vector of scalar coefficients satisfying
the equality
\[
\phi^{i+1}_1 a^1+\phi^{i+1}_2 a^2+\phi^{i+1}_3 a^3 = Ru^{i+1}+\beta f(u^{i+1}),
\]
and define
\begin{equation}\label{e28b}
\phi(t):=\sum_{i=0}^{N-1}\,1_{(i\frac{T}{N}, (i+1)\frac{T}{N}]}(t)\;
\phi^{i+1},\;\;\mbox{ }t\in[0,T].
\end{equation}
We remark that the function $\phi$ depends on the chosen values of $\delta$ and
$T$, the material parameters $\domain$, $\alpha$ and $\beta$ and the noise
parameters $a^1$, $a^2$ and $a^3$.

\smallskip

\noindent Recall that $Y^\eps$ denotes the solution of (\ref{e27}). By an
argument very much like that leading to Theorem~\ref{thm-large deviations}, the
family of laws $\{\mathcal{L}(Y^\eps): \eps\in(0,1)\}$ on $\mathcal X_T$
satisfies a large deviation principle. In order to define the rate function, we
introduce an equation
\begin{eqnarray}
\nonumber y_\psi(t)&=&\zeta_{-}+\imp_0^ty_\psi\times\Delta
y_\psi\,ds-\alpha\imp_0^ty_\psi\times(y_\psi\times\Delta y_\psi)\,ds\\
\label{upsi}&&{}-\beta\imp_0^tG\la y_\psi\ra f\la y_\psi\ra\,ds+\imp_0^t G\la y_\psi\ra B\psi\,ds.
\end{eqnarray}
By Corollary \ref{strong} this equation has unique solution $y_\psi\in\mathcal X_T$ for every $\psi\in L^2\la
0,T;\R^3\ra$. The rate function $I:\mathcal X_T\rightarrow [0,\infty]$, is defined by:
\begin{equation}\label{ratefn}
I_T(v)
:=\inf\left\{\frac{1}{2}\imp_0^T|\psi(s)|^2\,ds : \psi\in L^2(0,T;\R^3)\mbox{ and }\;v=y_\psi\right\},
\end{equation}
where the infimum of the empty set is taken to be $\infty$.

\noindent Let $y$ be the solution of equation (\ref{exe2}) with $y_0=\zeta_{-}$ and $\mathcal{K}$ as defined in
(\ref{e28}). Using the notation in (\ref{upsi}), we have $y=y_\phi$, for $\phi$ defined in (\ref{e28b}).  Therefore
\[I_T(y)\le\frac{1}{2}\imp_0^T|\phi(s)|^2\,ds<\infty.\]
 Since $y$ undergoes magnetization reversal by time $T$, paths of $Y^\eps$ which lie close to $y$ also undergo
magnetization reversal by time $T$. In particular, by the Freidlin-Wentzell formulation of the lower bound of the large
deviation principle (see, for example, \cite[Proposition 12.2]{DZ}), given $\xi>0$, there exists an $\eps_0>0$ such that
for all $\eps\in(0,\eps_0)$ we have
\begin{eqnarray} &&\!\!\!\!\!\!\!\!\!\!\!\!\!\!\mathbb{P}\la\textstyle{\sup_{t\in
[0,T]}}|Y^\eps(t)-y(t)|_{\mathbb{H}^1}+\la\textstyle{\imp_0^T}|Y^\eps(s)-y(s)|^2_{\xone}\,ds\ra^\frac{1}{2}<\frac{
\delta}{2}\ra\nonumber\\
&\geq&\exp\left(\textstyle{\frac{-I_T(y)-\xi}{\eps}}\right)\nonumber\\
&\geq&\exp\left(\textstyle{\frac{-\frac{1}{2}\imp_0^T|\phi(s)|^2\,ds-\xi}{\eps}}\right). \label{e29}
\end{eqnarray}
Since we have $|y(T)-\zeta_{+}|_{\mathbb{H}^1}<\frac{\delta}{2}$, the right hand side of (\ref{e29}) provides a
lower bound for the probability that $Y^\eps$ undergoes magnetization reversal by time $T$. We summarize our conclusions
in the following proposition.
\begin{proposition}\label{magrev}
For all sufficiently small $\eps>0$, the probability that the solution $Y^\eps$ of (\ref{e27}) undergoes magnetization
reversal by time $T$ is bounded below by the expression on the right hand side of (\ref{e29}); in particular, it is
positive.
\end{proposition}
\noindent We shall now use the upper bound of the large deviation principle satisfied by $\{\mathcal{L}(Y^\eps):
\eps\in(0,1)\}$ to find an exponential upper bound for the probability that small noise in the field drives the
magnetization outside a given $\mathbb{H}^1$-ball centred at the initial state $\zeta_{-}$ in time interval
$[0,T]$. This is done in Proposition~\ref{expuppbd} below; the proof of the proposition uses Lemma~\ref{Dupsi0}. In
Lemma~\ref{Dupsi0} and Proposition~\ref{expuppbd}, for $h$ an arbitrary element of $L^2(0,T;\mathbb{R}^3)$, $y_h$
denotes the function in $\mathcal X_T$ which satisfies equality (\ref{upsi}) and $\tau_h$ is defined by
\[\tau_h:=\inf\left\{t\in [0,T]:
|y_h(t)+\zeta_{+}|_{\mathbb{H}^1}\geq\frac{1}{2k^2\sqrt{|\domain|}}\frac{\alpha}{1+2\alpha}\right\}.\]
\begin{lemma}\label{Dupsi0}
For each $h\in L^2(0,T;\mathbb{R}^3)$, we have $|\D y_h(t)|_{\mathbb{H}}=0$ for all $t\in[0,\tau_h\wedge
T)$.
\end{lemma}
\begin{proof} Let $h\in L^2(0,T;\mathbb{R}^3)$. To simplify notation in this proof, we write $y$ instead of $y_h$.
Proceeding as in the derivation of (\ref{Dubound}), we obtain
\begin{eqnarray}
\nonumber |\D y(t)|^2_\mathbb{H}&=&-2\alpha\imp_0^t|y\times\Delta y|^2_{\mathbb{H}}\,ds+2\beta\imp_0^t\imp_\domain
R\,dx\,ds\\
\label{Dupsibound} &&{}-2\alpha\sum_{i=1}^3\imp_0^t\langle\D y, y\times(\D y\times
a^i)\rangle_{\mathbb{H}}h_i\,ds,\quad t\in [0,T],
\end{eqnarray}
where $R(s)$ defined in (\ref{Rsdef}) satisfies inequality (\ref{Rsineq}). For each $s\in[0,\tau_h\wedge T)$, $y(s)$
satisfies the hypotheses of Lemma~\ref{uandz}, thus we have $R(s)(x)\le 0$ for all $x\in\domain$. It follows from
(\ref{Dupsibound}) that for all $t\in[0,\tau_h\wedge T)$:
\begin{equation}\label{Dygronineq}
|\D y(t)|^2_{\mathbb{H}}\le 2\alpha\imp_0^t|\D y|^2_{\mathbb{H}}\sum_{i=1}^3|a^i|\cdot|h_i|\,ds.
\end{equation}
By the Gronwall lemma applied to (\ref{Dygronineq}), $|\D y(t)|_{\mathbb{H}}^2=0$ for all $t\in[0,\tau_h\wedge
T)$.
\end{proof}
\begin{proposition}\label{expuppbd}
Let
\[0<r<\rho\le\frac{1}{2k^2\sqrt{|\domain|}}\frac{\alpha}{1+2\alpha}.\]
The for any $\xi>0$, there exists $\eps_0>0$ such that for all $\eps\in(0,\eps_0)$:
\begin{equation}\label{FWuppbd}
\mathbb{P}\la\sup_{t\in [0,T]}|Y^\eps(t)+\zeta_{+}|_{\mathbb{H}^1}\geq\rho\ra\le\exp\left(\frac{-\kappa
r^2+\xi}{\eps}\right),
\end{equation}
where
\[\kappa=\frac{\alpha\beta}{8\max_{1\le i\le 3}|a^i|^2\,|\domain|(1+\alpha^2)}\,.\]
\end{proposition}
\noindent
\begin{proof} We shall use the Freidlin-Wentzell formulation of the upper bound of the large deviation principle (see,
for example, \cite[Proposition 12.2]{DZ}) satisfied by $\{\mathcal{L}(Y^\eps): \eps\in(0,1)\}$. Recall that
$\mathcal{I}$, defined in (\ref{ratefn}), is the rate function of the large deviation principle. Our main task is to
show that
\[\left\{v\in\mathcal X_T:\,I_T(v)\le\kappa r^2\right\}\subset\left\{v\in C([0,T];\mathbb{H}^1): \sup_{t\in
[0,T]}|v(t)+\zeta_{+}|_{\mathbb{H}^1}\le r\right\}.\]
Take $h\in L^2(0,T;\mathbb{R}^3)$ such that \begin{equation}\label{psibound}
\frac{1}{2}\imp_0^T|h(s)|^2\,ds\le\kappa r^2.
\end{equation}
For simplicity of notation, in this proof we write $y$ in place of $y_h$.  By Lemma~\ref{Dupsi0} we have for all
$t\in [0,T]$,
\begin{equation}\label{utaupsi1}
\begin{aligned}
|y(t\wedge\tau_h)+\zeta_{+}|^2_{\mathbb{H}^1}&=2\alpha\imp_0^{t\wedge\tau_h}\imp_\domain|\D y|^2\la y\cdot
\zeta_{+}\ra\,dx\,ds\\
&\,\,\,\,\,\,\,+2\alpha\beta\imp_0^{t\wedge\tau_h}\imp_\domain\la y\cdot \zeta_{+}\ra |y\times \zeta_{+}|^2\,dx\,ds\\
&\,\,\,\,\,\,\,-2\alpha\beta\sum_{i=1}^3\imp_0^{t\wedge\tau_h}\left\langle\frac{1}{2}(y\times \zeta_{+}),
\frac{2}{\alpha\beta}a^i\right\rangle_{\mathbb{H}}h_i\,ds\\
&\,\,\,\,\,\,\,+2\alpha\beta\sum_{i=1}^3\imp_0^{t\wedge\tau_h}\left\langle\frac{1}{2}(y\times
\zeta_{+}),\frac{2}{\beta}(y\times a^i)\right\rangle_{\mathbb{H}}h_i\,ds\\
&\le-\frac{3}{2}\alpha\beta\imp_0^{t\wedge\tau_h}|y\times
\zeta_{+}|^2_{\mathbb{H}}\,ds+\frac{3}{2}\alpha\beta\imp_0^{t\wedge\tau_h}|y\times \zeta_{+}|^2_{\mathbb{H}}\,ds\\
&\,\,\,\,\,\,\,+\frac{4}{\beta}\la\frac{1}{\alpha}+\alpha\ra
|\domain|\sum_{i=1}^3|a^i|^2\imp_0^{t\wedge\tau_h}h_i^2\,ds,
\end{aligned}
\end{equation}
where we estimated the integrals on the right hand side of the second equality as follows: the first integral vanished
thanks to Lemma~\ref{Dupsi0}, Lemma~\ref{uandz} was used for the integrand of the second integral and the Cauchy-Schwarz
inequality and Young's inequality were used for the integrands of the other integrals.
Using (\ref{psibound}) in (\ref{utaupsi1}), we obtain
\begin{equation}\label{utaupsi2}
|y(t\wedge\tau_h)+\zeta_{+}|_{\mathbb{H}^1}\le r<\frac{1}{2k^2\sqrt{|\domain|}}\frac{\alpha}{1+2\alpha}\;\;\;\mbox{
}\forall t\in[0,T].
\end{equation}
From (\ref{utaupsi2}) and the definition of $\tau_h$, we conclude that $\tau_h>T$. Hence we have
\[\sup_{t\in[0,T]}|y(t)+\zeta_{+}|_{\mathbb{H}^1}\le r.\]
\noindent By the Freidlin-Wentzell formulation of the upper bound of the large deviation principle, since $r<\rho$,
given $\xi\in(0,\infty)$, there exists $\eps_0>0$ such that for all $\eps\in(0,\eps_0)$, inequality (\ref{FWuppbd})
holds.
\end{proof}

\smallskip

\begin{remark} Our use of Lemma~\ref{Dupsi0} in the proof of Proposition~\ref{expuppbd} means that, in this proposition,
we did not need to allow for the spatial variation of magnetization on $\domain$.
\end{remark}

\appendix

\section{Budhiraja-Dupuis result}

Let us recall Theorem 3.6 from \cite{Budhiraja+Dupuis_2000}.

\begin{theorem}\label{thm-3.6_BD}
Assume that $\rK,\rH$ be  separable Hilbert spaces such that the embedding
\begin{equation}\label{eqn-K-H}
\rK \embed \rH
\end{equation}
is $\gamma$-radonifying,
and
\[
f: \fourIdx{}{0}{}{}{C}([0,T];\rH) \to \mathbb{R}
\]
be a bounded (or bounded from below?) Borel measurable function. Then
\begin{equation}\label{eqn-3.14_BD}
- \log \mathbb{E} e^{-f (W)} =
\inf_{h \in \mathcal{A}} \mathbb{E} \Bigl( \frac12 \int_0^T \vert h(s)\vert^2_{\rK}+ f\bigl(W+ \int_0^{\cdot} h(s) \,ds\bigr) \Bigr),
\end{equation}
where $\mathcal{A}$ consist of all $\rK$-valued predictable processes $h$ such that
\begin{equation}\label{eqn-2.1_BD}
\mathbb{P} \Bigl\{ \int_0^T  \vert h(s)\vert^2_{\rK}< \infty \Bigr\}=1.
\end{equation}
\end{theorem}

For $R>0$ we denote by $\mathcal{A}_R$ the subset of $\mathcal{A}$ consisting of
of all $\rK$-valued predictable processes $h$ satisfying
\begin{equation}\label{eqn-B-04}
\mathbb{P} \Bigl\{ \int_0^T  \vert h(s)\vert^2_{\rK}\leq R^2 \Bigr\}=1.
\end{equation}
Note that $\bigcup_{R>0}\mathcal{A}_R $ is a proper subset of $\mathcal{A}$.
Let us also denote by $B_R$ the closed ball of radius $R$ in the set $L^2(0,T;\rK)$, i.e.
\begin{equation}\label{eqn-B-05}
B_R:=\bigl\{ h \in  L^2(0,T;\rK) :\int_0^T  \vert h(s)\vert^2_{\rK}\leq R^2\bigr\}.
\end{equation}
We endow $B_R$ with the weak topology induced by $L^2(0,T;\rK)$.

Let now $E$ be a Polish space and consider a family, indexed by $\eps \in (0,1]$, of Borel measurable maps
\[
J^\eps: \fourIdx{}{0}{}{}{C}([0,T];\rH) \to E.
\]

On the space $\fourIdx{}{0}{}{}C([0,T];\rH)$ we consider a Wiener measure $\mathbb{P}$ corresponding to the embedding  \eqref{eqn-K-H} (and the integration w.r.t. $\mathbb{P}$ we denote by $\mathbb{E}$). Note that the RKHS of $\mu$ is not the space $L^2(0,T;\rK)$ but the space $\fourIdx{}{0}{1,2}{}{H}(0,T;\rK)$, where
\[
\fourIdx{}{0}{1,2}{}{H}(0,T;\rK)=\bigl\{ \omega\in \fourIdx{}{0}{}{}{C}([0,T];\rK):  \omega^\prime\in L^2(0,T;\rK) \bigr\}.
\]
Note that the map
\[
 L^2(0,T;\rK) \ni h \mapsto \int_0^{\cdot} h(s)\, ds \in  \fourIdx{}{0}{1,2}{}{H}(0,T;\rK)
\]
is an isometric isomorphism.

We denote by $\mu^\eps$ the "image" measure on $E$ of $\mathbb{P}$ by $J^\eps$, i.e.
\begin{equation}\label{eqn-4.1_BD}
\mu^\eps= J^\eps(\mathbb{P}), \mbox{ i.e. } \mu^\eps(A)=\mathbb{P}\bigl( (J^\eps)^{-1}(A)\bigr), \;\; A\in \mathcal{B}(E).
\end{equation}

Assume
\begin{ass}\label{ass-4.3_BD}
There exists a Borel measurable map
\[
J^0: \fourIdx{}{0}{}{}{C}([0,T];\rH) \to E
\]
such that
\begin{trivlist}
\item[\textbf{(BD1)}]
 if $R>0$ and a family $h_\eps \subset \mathcal{A}_R$ converges in law on $B_R$ to $h \in \mathcal{A}_R$, then the processes
\[
\fourIdx{}{0}{}{}{C}([0,T],\rH) \ni \omega \mapsto J^\eps( \omega+ \frac{1}{\sqrt{\eps}}\int_0^{\cdot} h_\eps(s) \,ds ) \in E \]
converge in law, as $\eps \todown 0$,  to the process $ J^0(\int_0^{\cdot} h(s) \,ds)
$,
\item[] and
\item[\textbf{(BD2)}]
 the set
\[
\Bigl\{  J^0(\int_0^{\cdot} h(s) \,ds) : h\in B_R \Bigr\}
\]
is compact in $E$.
\end{trivlist}
\end{ass}

We have the following result.

\begin{theorem}\cite[Theorem 4.4]{Budhiraja+Dupuis_2000}\label{thm-4.4_BD} If the assumptions listed above, in particular Assumption \ref{ass-4.3_BD},  are satisfies, then
the family of measures $\mu_\eps$ satisfies the LDP with the rate function $I$ defined by
\begin{equation}\label{eqn-4.3_BD}
I(u):=\inf\bigl\{ \frac12 \int_0^T \vert h(s)\vert^2_{\rK}: h\in L^2(0,T;\rK) \mbox{ and } u=J^0(\int_0^{\cdot}h(s)\,ds), \bigr\}.
\end{equation}
\end{theorem}
Note that
\begin{equation}\label{eqn-B08}
I(u):=\inf\bigl\{ \frac12 \int_0^T \vert y^\prime(s)\vert^2_{\rK}: y\in  \fourIdx{}{0}{1,2}{}{H}(0,T;\rK)  \mbox{ and } u=J^0(y), \bigr\}.
\end{equation}
Obviously, we put, as always, $\inf \emptyset=\infty$.

\end{document}